\documentclass{amsart}
\usepackage{amssymb,amsmath,epsfig,amscd}
\usepackage[mathscr]{eucal}
\usepackage[all]{xy}

\theoremstyle{plain}
\newtheorem{thm}[equation]{Theorem}
\newtheorem{prop}[equation]{Proposition}
\newtheorem{lemma}[equation]{Lemma}
\newtheorem{cor}[equation]{Corollary}
\newtheorem{hyp}[equation]{Hypotheses}
\newtheorem{conj}[equation]{Conjecture}
\theoremstyle{definition}
\newtheorem{remark}[equation]{Remark}
\newtheorem{example}[equation]{Example}
\newtheorem{defn}[equation]{Definition}

\numberwithin{equation}{section}

\newcommand{\F}{{\mathbb F}}

\newcommand{\boldF}{{\mathbb F}}

\newcommand{\Fp}{{\mathbb F}_p}
\newcommand{\Fpp}{{\mathbb F}_{p^2}}
\newcommand{\Fpbar}{\overline{\boldF}_p}

\newcommand{\OO}{\mathcal O}
\newcommand{\Z}{{\mathbb Z}}

\newcommand{\Zp}{\Z_p}

\newcommand{\Q}{{\mathbb Q}}

\newcommand{\Qp}{{\mathbb Q}_p}
\newcommand{\Qpp}{{\mathbb Q}_{p^2}}
\newcommand{\Qpbar}{\overline{\Q}_p}

\newcommand{\kk}{{\mathbf k}}

\newcommand{\rhobar}{\overline{\rho}}
\newcommand{\om}{\tilde{\omega}}
\newcommand{\omt}{\om_2}

\newcommand{\Gal}{{\rm Gal}}
\newcommand{\Hom}{{\rm Hom}}

\newcommand{\GL}{{\rm GL}}

\newcommand{\Frob}{{\rm Frob}}

\DeclareMathOperator{\Spec}{{\rm Spec}}
\DeclareMathOperator{\End}{{\rm End \ }}
\newcommand{\WD}{{\rm WD}}
\newcommand{\Fil}{{\rm Fil}}
\newcommand{\val}{{\rm val}}

\newcommand{\G}{\mathcal G}

\newcommand{\DstF}{D_{\mathrm{st}}^{F}}
\newcommand{\VstF}{V_{\mathrm{st}}^{F'}}
\newcommand{\DstkF}{D_{\mathrm{st},k}^{F}}
\newcommand{\DstFtwo}{D_{\mathrm{st},2}^{F}}
\newcommand{\VstkF}{V_{\mathrm{st},k}^{F'}}
\newcommand{\VstFtwo}{V_{\mathrm{st},2}^{F'}}
\newcommand{\TstFtwo}{T_{\mathrm{st},2}^{F'}}
\newcommand{\Tst}{T_{\mathrm{st},k}}
\newcommand{\Rmod}{R-\Mod^{k-1}_{\mathrm{cris},\mathrm{dd}}}
\newcommand{\mE}{\mathfrak{m}_E}
\newcommand{\M}{{\mathscr M}}
\newcommand{\NN}{{\mathscr N}}
\newcommand{\D}{{\mathcal D}}
\newcommand{\invlim}{\varprojlim}
\newcommand{\X}{{\mathfrak X}}
\newcommand{\N}{{\mathbb N}}
\newcommand{\Ocris}{{\mathscr O}^{\mathrm{cris}}_{n,\pi}}
\newcommand{\Jcris}{{\mathscr J}^{\mathrm{cris}}_{n,\pi}}
\newcommand{\Ocrisi}{{\mathscr O}^{\mathrm{cris}}_{\infty,\pi}}
\newcommand{\Jcrisi}{{\mathscr J}^{\mathrm{cris}}_{\infty,\pi}}
\newcommand{\Ocrisg}{{\mathscr O}^{\mathrm{cris},(g)}_{n,\pi}}
\newcommand{\Jcrisg}{{\mathscr J}^{\mathrm{cris},(g)}_{n,\pi}}
\newcommand{\Ocrisgi}{{\mathscr O}^{\mathrm{cris},(g^{-1})}_{n,\pi}}
\newcommand{\AbR}{({\rm Ab}/R)}
\newcommand{\syn}{{\rm Spf}(R)_{\rm syn}}
\newcommand{\OOO}{{\mathscr O}}
\newcommand{\p}{\phi_1}
\newcommand{\ghat}{\widehat{g}}
\newcommand{\Acrish}{\widehat{A}_{\mathrm{cris}}}
\newcommand{\Asth}{\widehat{A}_{\mathrm{st}}}
\newcommand{\Acris}{{A_{\mathrm{cris}}}}
\newcommand{\Asthi}{\widehat{A}_{\mathrm{st},\infty}}
\newcommand{\Acrishi}{\widehat{A}_{\mathrm{cris},\infty}}
\newcommand{\Acrisi}{{A_{\mathrm{cris},\infty}}}
\newcommand{\upi}{\underline{\pi}}
\newcommand{\Mod}{{\rm Mod}}
\newcommand{\BrMod}{{\rm BrMod}}
\newcommand{\uMod}{\underline{\Mod}}
\newcommand{\mm}{\mathfrak{m}}
\newcommand{\e}{\mathbf{e}}
\newcommand{\ox}{\overline{x}}
\newcommand{\ob}{\overline{b}}
\newcommand{\ow}{\overline{w}}

\begin{document}
\title{On a Conjecture of Conrad, Diamond, and Taylor}
\author{David Savitt }
\date{April 19, 2004. \textit{Revised:} September 3, 2004. 
\\ Author's work supported by a National Science Foundation (NSF)
International Research Fellowship and by the Office of
Multidisciplinary Activities under the NSF Mathematics and Physical
Sciences Directorate Program.  This work also supported by a Natural
Sciences and Engineering Research Council of Canada (NSERC)
postdoctoral fellowship.}



\begin{abstract}  We prove a conjecture of Conrad, Diamond, and
Taylor on the size of certain deformation rings parametrizing
potentially Barsotti-Tate Galois representations.  To achieve
this, we extend results of Breuil and M\'ezard (classifying Galois
lattices in semistable representations in terms of ``strongly
divisible modules'') to the potentially crystalline case in
Hodge-Tate weights $(0,1)$. We then use these strongly divisible
modules to compute the desired deformation rings.  As a corollary,
we obtain new results on the modularity of potentially
Barsotti-Tate representations.
\end{abstract}
 \maketitle
\section{Introduction}

In their paper \cite{CDT}, Conrad, Diamond, and Taylor
conjectured that certain deformation rings parametrizing
potentially Barsotti-Tate Galois representations are sufficiently
small for the methods of Taylor-Wiles to yield a modularity
result. Breuil and M\'ezard \cite{BreuilMezard} reformulated and
vastly generalized these conjectures, and proved their new
conjectures for semistable Galois representations in even weight.
In this article, essentially a sequel to \cite{BreuilMezard}, we
prove the conjectures of Breuil and M\'ezard in the cases
originally conjectured by Conrad, Diamond, and Taylor.

We now describe these conjectures.  Fix $p$ an odd prime, and let
$E$ be a finite extension of $\Qp$ with residue field $\kk_E$. To
each potentially crystalline Galois representation $\rho : G_{\Qp}
\rightarrow \GL_2(E)$, we attach a representation $\WD(\rho)$ of the
Weil group $W_{\Qp}$ (see Def. \ref{galtype}), and hence a
Galois type $\tau(\rho) = \WD(\rho) \,|_{I_p}$.

Suppose that $\rhobar : G_{\Qp} \rightarrow \GL_2(\kk_E)$ is such
that $\End_{\kk_E[G_{\Qp}]} \rhobar = \kk_E$; we shall say that
$\rhobar$ has trivial endomorphisms.  Let
$R^{\mathrm{univ}}_{\OO_E}(\rhobar)$ be the universal deformation
ring parametrizing deformations of $\rhobar$ over complete local
noetherian $\OO_E$-algebras.  If $2 \le k < p$ and if $\OO_{E'}$ are
the integers in a finite extension of $E$, we say that a deformation
$\rho : G_{\Qp} \rightarrow \GL_2(\OO_{E'})$ of $\rhobar$ has type
$(k,\tau)$ if
\begin{itemize}
\item $\rho$ is potentially semi-stable and $\tau(\rho) \cong
\tau$,
\item $\rho$ has Hodge-Tate weights $(0,k-1)$, and
\item $\det(\rho)$ is a fixed lift of $\det(\rhobar)$ of the following form:
the $(k-1)$st power of the $p$-adic
cyclotomic character times a  finite character of order prime
to $p$.
\end{itemize}
The kernel $\mathfrak{p}$ of the corresponding map
$R^{\mathrm{univ}}_{\OO_E}(\rhobar) \rightarrow \OO_{E'}$ is also
said to have type $(k,\tau)$, and we define
$$ R(k, \tau, \rhobar)_{\OO_E} = R^{\mathrm{univ}}_{\OO_E}(\rhobar) /
\underset{{\mathfrak p} \ \text{type} \ (k,\tau)}{\bigcap}
{\mathfrak p} . $$

The first part of the conjectures of Breuil and M\'ezard
(see \cite[Conj. 2.2.2.4]{BreuilMezard}) posits that
$R(k,\tau,\rhobar)_{\OO_E}$ should be equidimensional of Krull
dimension~$2$, and that $R(k,\tau,\rhobar)_{\OO_E} \otimes E$ should
be regular.  Let $\mu_{\mathrm{gal}}(k,\tau,\rhobar)$ be the Samuel
multiplicity of $\overline{R} = R(k,\tau,\rhobar)_{\OO_E}
\otimes_{\OO_E} \kk_E$; so conjecturally, this is $\dim_{\kk_E}
\mm_{\overline{R}}^n / \mm_{\overline{R}}^{n+1}$ for $n$
sufficiently large.  Via a recipe on the automorphic side, Breuil
and M\'ezard also define an integer
$\mu_{\mathrm{aut}}(k,\tau,\rhobar)$ (see
\cite[Sec 2.1]{BreuilMezard} for the details).  We then have the following.

\begin{conj}[\cite{BreuilMezard}, Conj. 2.3.1.1] \label{conj:bm}   If
 \ $\det(\tau)$ is tame, then $$\mu_{\mathrm{gal}}(k,\tau,\rhobar) =
\mu_{\mathrm{aut}}(k,\tau,\rhobar).$$
\end{conj}

The conjectures of Conrad, Diamond, and Taylor to which we have
referred (see \cite[Conjs. 1.2.2, 1.2.3]{CDT}) are, more or
less, the case $k=2$ and $\tau$ tamely ramified in Conjecture
\ref{conj:bm}. Our main theorem, then, is the following.

\begin{thm} Conjecture \ref{conj:bm} holds when $k=2$ and $\tau$
is tamely ramified.
\end{thm}

Indeed, we show the following (see Exams. \ref{characters},
\ref{char2}, \ref{rednotation} for notation and Ths.
\ref{main-prin} and \ref{main-super} for more precise statements).

\begin{thm} \label{thm:main1} Suppose that $\rhobar : G_{\Qp} \rightarrow \GL_2(\kk_E)$
has trivial endomorphisms.  Suppose that $\tau \cong \om^i \oplus
\om^j$ with $i \not\equiv j \pmod{p-1}$.  Then we have the
following.
\begin{enumerate}
\item $\mu_{\mathrm{gal}}(2,\tau,\rhobar) = 0$ if $\rhobar \,|_{I_p}
\otimes_{\kk_E} \Fpbar \not\in \left\{
\begin{pmatrix}
\omega^{1+i} & * \\
0 & \omega^j
\end{pmatrix},
\begin{pmatrix}
\omega^{1+j} & * \\
0 & \omega^i
\end{pmatrix},
\omega_2^k \oplus \omega_2^{pk}\right\}$ with $k =
1+\{j-i\}+(p+1)i$, where $\{a\}$ is the unique integer in
$\{0,\ldots,p-2\}$ which is congruent to $a \pmod{p-1}$;

\item $\mu_{\mathrm{gal}}(2,\tau,\rhobar) = 1$ if $\rhobar \,|_{I_p} \otimes_{\kk_E} \Fpbar
\in \left\{
\begin{pmatrix}
\omega^{1+i} & * \\
0 & \omega^j
\end{pmatrix},
\begin{pmatrix}
\omega^{1+j} & * \\
0 & \omega^i
\end{pmatrix}\right\}$,

\item
$\mu_{\mathrm{gal}}(2,\tau,\rhobar) = 2$ if $\rhobar \,|_{I_p}
\otimes_{\kk_E} \Fpbar \cong \omega_2^k \oplus \omega_2^{pk}$ with
$k = 1+\{j-i\}+(p+1)i$.
\end{enumerate}
\end{thm}

\begin{thm} \label{thm:main2} Suppose that $\rhobar : G_{\Qp} \rightarrow \GL_2(\kk_E)$
has trivial endomorphisms.  Suppose that $\tau \cong \omt^m \oplus
\omt^{pm}$ with $p+1 \nmid \, m$.  Write $m = i + (p+1)j$ with $i
\in \{1,\ldots,p\}$ and $j \in \Z/(p-1)\Z$.  Then we have the
following.
\begin{enumerate}
\item $\mu_{\mathrm{gal}}(2,\tau,\rhobar)= 1$ if $\rhobar
\,|_{I_p} \otimes_{\kk_E} \Fpbar \in \left\{
\begin{pmatrix}
\omega^{i+j} & * \\
0 & \omega^{1+j}
\end{pmatrix},
\begin{pmatrix}
\omega^{1+j} & * \\
0 & \omega^{i+j}
\end{pmatrix}\right\}$, the first
$*$ peu ramifi\'e when $i=2$ and the second when $i=p-1$;

\item $\mu_{\mathrm{gal}}(2,\tau,\rhobar) = 1$ if $\rhobar
\,|_{I_p} \otimes_{\kk_E} \Fpbar \in \left\{ \omega_2^{p+m} \oplus
\omega_2^{1+pm}, \omega_2^{1+m} \oplus \omega_2^{p(1+m)}
\right\}$;

\item $\mu_{\mathrm{gal}}(2,\tau,\rhobar)= 0$ otherwise.
\end{enumerate}
\end{thm}

We note an important consequence of these results.  The method of
Taylor-Wiles, as utilized in \cite{BCDT}, may be reformulated as
follows.

\begin{thm}[\cite{BCDT}, Th. 1.4.1] \label{tw}   Let $\rho : G_{\Q}
\rightarrow \GL_2(E)$ be an odd continuous representation ramified
at only finitely many primes.  Assume that its reduction $\rhobar :
G_{\Q} \rightarrow \GL_2(\kk_E)$ is modular and is absolutely
irreducible after restriction to $\Q(\sqrt{(-1)^{(p-1)/2} p})$.
Further, suppose that
\begin{itemize}
\item $\rhobar \,|_{G_{\Qp}}$ has trivial endomorphisms,
\item $\rho_p = \rho \,|_{G_{\Qp}}$ is potentially Barsotti-Tate,
and
\item $\mu_{\mathrm{gal}}(2,\tau(\rho_p),\rhobar) \le 1 \le
\mu_{\mathrm{aut}}(2,\tau(\rho_p),\rhobar)$.
\end{itemize}
Then $\rho$ is modular.
\end{thm}

The import of Conjecture \ref{conj:bm} is that if it is true, then the
last condition of Theorem \ref{tw} may be replaced with
$\mu_{\mathrm{aut}}(2,\tau(\rho_p),\rhobar) \le 1$, removing the
irksome hypothesis involving $\mu_{\mathrm{gal}}$ from the theorem.
In particular, we obtain the following immediate corollary of
Theorems \ref{thm:main1} and \ref{thm:main2} (together with the
$k=2$, $\tau$ scalar case of Conj. \ref{conj:bm}, proved in
\cite{BreuilMezard}).

\begin{thm} \label{modularity} Let $\rho : G_{\Q}
\rightarrow \GL_2(E)$ be an odd continuous representation ramified
at only finitely many primes.  Assume that its reduction $\rhobar
: G_{\Q} \rightarrow \GL_2(\kk_E)$ is modular and is absolutely
irreducible after restriction to $\Q(\sqrt{(-1)^{(p-1)/2} p})$.
Further, suppose that
\begin{itemize}
\item $\rhobar \,|_{G_{\Qp}}$ has trivial endomorphisms,
\item $\rho_p = \rho \,|_{G_{\Qp}}$ is potentially Barsotti-Tate,
\item $\tau(\rho)$ is tamely ramified, and
\item if $\tau(\rho) \cong \om^i \oplus \om^j$ with $i \not\cong
j \pmod{p-1}$ then $\rhobar \,|_{G_{\Qp}} \otimes_{\kk_E} \Fpbar$
is reducible.
\end{itemize}
Then $\rho$ is modular.
\end{thm}

This is a significant improvement on the main results in
\cite{SavittCompositio}, where $\rhobar \,|_{G_{\Qp}}$ had to be
reducible and defined over $\Fp$.  It would be of interest to know
whether the methods of Taylor-Wiles could be extended to handle
cases where
$\mu_{\mathrm{gal}}(2,\tau,\rhobar)=\mu_{\mathrm{aut}}(2,\tau,\rhobar)
= 2$, in order to remove the last hypothesis from this
theorem.\footnote{This question is resolved in considerable
generality in a new preprint of Mark Kisin \cite{KisinModularity}.}

We give a brief outline of this paper.  We follow the same
strategy established by Breuil and M\'ezard to prove Conjecture
\ref{conj:bm} in the case $\tau$ scalar, $k$ even.  To achieve this,
we must provide (as best we can) ``potential'' versions, when $k=2$,
of the machinery of \cite{BreuilMezard} classifying lattices in
semi-stable Galois representations by means of strongly divisible
modules.  We begin in Section \ref{filmod} by recalling Fontaine's
filtered modules with coefficients and descent data and computing
the particular filtered modules that arise in the proofs of our main
theorems.

Sections \ref{strtame} and \ref{coef} contain the bulk of the
technicalities: in the former, we use the equivalence between
$p$-divisible groups and lattices in potentially Barsotti-Tate
representations to add (tame) descent data to the strongly divisible
modules of \cite{BreuilMezard} when $k=2$; in the latter, we
introduce coefficients into the mix. Since we are working over a
base ring that may be highly ramified, the results of
\cite{BreuilMezard} do not entirely go over to our situation, and so
in some cases we must scrape by with weaker results.

Finally, we perform the calculations using strongly divisible
modules (with coefficients and descent data) necessary to prove
our main theorems.  In Section \ref{sec:chars} we perform
calculations with characters, and use these results repeatedly in
Section \ref{calcs}, which contains the bulk of our calculations.

We remark that, in the course of our work, we completely determine
(Ths. \ref{red-prin} and \ref{red-super}) the reductions (mod
$p$) of $2$-dimensional potentially Barsotti-Tate Galois
representations that become crystalline over a tamely ramified
extension of $\Qp$.  In Section \ref{modforms}, we apply these
results to re-prove an old result on the (mod $p$) representations
attached to modular forms, and to suggest a first step towards a new
one.

\begin{remark}
  The current version of this paper fixes an error
  in the published version, as a consequence of which there is one
  more family of strongly divisible modules that
we must study by the methods of this paper than was studied in the
published version.  The main results of this paper are unaffected.

The mistake is in the statement and proof of Theorem 6.12(4) of the
published version.  Let  $T := T_{\mathrm{st},2}^{\Qp}(\M)$ be the $p$-adic Galois
representation considered there.   In
the notation of that item, if $m = 1 + (p+1)j$ --- i.e., if
$i=1$ --- then the two characters $\omega_2^{m+p}$ and
$\omega_2^{pm+1}$ are both characters of niveau one, and are equal.
Hence the proof of Theorem 6.12(4) does \emph{not} show
that the reduction mod~$p$ of~$T$ has niveau two in this case, and
indeed that proof can be modified
to check that the reduction mod~$p$ of~$T$ is split.  
This means that when $i=1$ we 
still need to construct a strongly divisible lattice in
$T[1/p]$ whose reduction mod $p$ has trivial endomorphisms.  That
construction is contained in this version of the paper.

For further details, see the corrigendum on the
author's website.   Numbering in this version is consistent with
the published version.
\end{remark}

\section{Filtered modules with coefficients and descent data}
\label{filmod}

The purpose of this section is to provide ``potential'' versions
of the results in \cite[Sec. 3.1]{BreuilMezard}.

\subsection{Weil-Deligne representations}

Suppose that $E/K$ is an extension of fields, and suppose that $F/K$ is a
finite Galois extension. Endow $F \otimes_K E$ with an action of
$G=\Gal(F/K)$ by letting $G$ act naturally on the first factor and
trivially on the second. Let $g$ denote an element of $G$.  In this
section, we examine the structure of $(F \otimes_K E)$-modules with
equivariant $G$-actions, which we dub $(F,E,G)$-modules, for short.
By a map of $(F,E,G)$-modules, we mean an $(F \otimes_K E)$-module
homomorphism that is also a $G$-homomorphism.

\begin{lemma} \label{freeness} Every $(F,E,G)$-module is free.
\end{lemma}

\begin{proof} Let $M$ be an $(F,E,G)$-module.
Let $V = M^{G}$, the $G$-invariants of $M$.  By Galois descent, we
have $M = F \otimes_{K} V$ as $F$-vector spaces with an action of
$G$. But since $G$ acts trivially on $E$ we find that $V$ is
actually an $E$-vector space; since the actions of $F$ and $E$ on
$M$ commute, $M$ is a free $(F,E,G)$-module.
\end{proof}

For the remainder of this section, we consider what happens
when $F$ is actually contained inside $E$.

\begin{lemma}
\label{structA} If $E$ contains $F$, then the map $\theta$ taking $
f \otimes e \mapsto (\sigma(f)e)_{\sigma}$, and extended by
linearity, is an isomorphism
\begin{equation}
\label{eqA} \theta : F \otimes_K E \rightarrow \coprod_{ \sigma :
F \underset{K}{\hookrightarrow} E } E
\end{equation}
of $(F,E,G)$-modules where, on the right-hand side, if
$(e_{\sigma})_{\sigma}$ denotes the vector that has $e_{\sigma}$ in
the $\sigma$-component then $(f \otimes e) \cdot
(e_{\sigma})_{\sigma} = (\sigma(f) e_{\sigma} e)_{\sigma}$ and $g
\cdot (e_{\sigma})_{\sigma} = (e_{\sigma})_{\sigma \circ g^{-1}}$.
\end{lemma}

\begin{proof}
To begin, we note that right-hand side of \eqref{eqA} is indeed an
$(F,E,G)$-module and that the map $\theta$ is well-defined, after
which it is easy to see that $\theta$ is a map of
$(F,E,G)$-modules. But $\theta$ is surjective, since the elements
of $G$ are linearly independent over $E$, and so by a dimension
count $\theta$ is an isomorphism.
\end{proof}

\begin{prop}
\label{structB} If $E$ contains $F$, any $(F,E,G)$-module $M$ is
isomorphic to one of the form
$$ M \cong \coprod_{ \sigma : F
\underset{K}{\hookrightarrow} E } V $$
for some $E$-vector space $V$, with the $(F,E,G)$-module structure on the
right-hand side defined as in
Lemma \ref{structA}.
\end{prop}

\begin{proof} Let $E_{\sigma}$ be the $(F
\otimes_{K} E)$-submodule of $\coprod_{\sigma} E$ consisting of
elements that are nonzero at most in the position corresponding to
$\sigma$.  Let $I_{\sigma}$ be the ideal $\theta^{-1}(E_{\sigma})$
in $F \otimes_K E$, and put $M_{\sigma} = I_{\sigma} M$; if $\tau =
\sigma \circ g^{-1}$, then $g$ induces $E$-linear maps $E_{\sigma}
\rightarrow E_{\tau}$, $I_{\sigma} \rightarrow I_{\tau}$, and
$\mu_{\sigma,\tau} : M_{\sigma} \rightarrow M_{\tau}$.  By
definition, $\mu_{\sigma,\tau}$ and $\mu_{\tau,\sigma}$ must be
inverses of one another, and hence they are isomorphisms of
$E$-vector spaces.

Now, the summation map $\coprod M_{\sigma} \rightarrow M$ is
evidently surjective.  To prove injectivity, suppose that we have a
relation $\sum_{\sigma} m_{\sigma} = 0$ with each $m_{\sigma} \in
M_{\sigma}$. Note that $(f \otimes 1) m_{\sigma} = (1 \otimes \sigma f)
m_{\sigma}$ follows from the analogous relation in $I_{\sigma}$, and so
$$ \sum_{\sigma} (1 \otimes \sigma f) m_{\sigma} = 0 $$
for all $f \in F$.  It follows from the linear independence of
the elements of $G$ that $m_{\sigma}=0$ for all $\sigma$, and so
$M = \coprod_{\sigma} M_{\sigma}$.

Fix any $\tau : F \hookrightarrow E$.  We map $M$ bijectively to
$\coprod_{\sigma} M_{\tau}$ via the map $\coprod \mu_{\sigma,\tau}$.
One checks without difficulty that, with the desired
$(F,E,G)$-module structure on $\coprod_{\sigma} M_{\tau}$, this map
is an isomorphism of $(F,E,G)$-modules.  For example, $g m_{\sigma}
= \mu_{\sigma, \sigma \circ g^{-1}} m_{\sigma}$ is mapped to the
element that is equal to
$$\mu_{\sigma \circ g^{-1},\tau} \mu_{\sigma, \sigma \circ g^{-1}}
m_{\sigma} = \mu_{\sigma,\tau} m_{\tau}$$ in the $(\sigma \circ
g^{-1})$-position and zero elsewhere. \end{proof}

\begin{remark}  Essentially the same argument shows that each
$(F,E,G)$-submodule of $\coprod_{\sigma} V$ is equal to $\coprod_{\sigma}
W$ for some sub-$E$-vector space $W \subset V$.
\end{remark}

Now fix a group $H$ and a surjection $\phi : H \twoheadrightarrow G$.
Suppose that $M$ is an $F \otimes_K E$-module endowed with two
$\phi$-semilinear,
$E$-linear actions $\cdot_{1}$ and $\cdot_{2}$ of $H$: that is, if $m \in
M$, $f \otimes e \in F
\otimes_K E$, and $h \in H$, we ask that $h \cdot_{i} (f \otimes e)m =
(\phi(h)f \otimes e)(h \cdot_{i} m)$ for $i=1,2$.  Moreover, assume that
the two actions of $H$ commute with one another, and that the second
action factors through an abelian quotient of $H$.

As in the proof of Proposition \ref{structB}, $M$ decomposes as a
coproduct $\coprod M_{\sigma}$ of
$E$-vector spaces, where $M_{\sigma} = I_{\sigma}M$.  The preceeding
hypotheses allow us to define a representation of $H$ on each
$M_{\sigma}$.  Indeed, both actions of an element $h \in H$ induce a
map $M_{\sigma} \rightarrow M_{\sigma \circ \phi(h^{-1})}$, and so we
obtain an $E$-linear map $\rho_{\sigma}(h) : M_{\sigma} \rightarrow
M_{\sigma}$ by setting
$$\WD_{\sigma}(h) (m_{\sigma}) = h^{-1} \cdot_{2} h \cdot_{1}
m_{\sigma}.$$
The commutativity hypotheses on the two actions guarantee that
$\WD_{\sigma}$ is a
representation.  Moreover, each $h$ induces an
isomorphism $\WD_{\sigma} \rightarrow \WD_{\sigma \circ \phi(h)^{-1}}$
via the second action; since $\phi$ is surjective, all of the
$\WD_{\sigma}$ are isomorphic.

\begin{defn} The isomorphism class of the $\WD_{\sigma}$ is called the
\textit{Weil-Deligne representation of $H$ attached to $M$}, and is
denoted $\WD(M)$.
\end{defn}

\subsection{Weakly admissible filtered modules} \label{sub:weakly}

Let $p$ be an odd prime.  Choose an algebraic closure $\Qpbar$ of
$\Qp$, let $E$ and $F$ be finite extensions of $\Qp$ inside
$\Qpbar$, and let $F'$ be a field lying between $\Qp$ and $F$ such
that $F/F'$ is Galois. Fix the uniformizer $p \in \Qp$, thereby
fixing an inclusion $B_{\mathrm{st}} \rightarrow B_{\mathrm{dR}}$.
Let $F_0$ denote the maximal unramified extension of $\Qp$ contained
in $F$. We retain this notation for the remainder of the
paper.

\begin{defn}  A \textit{filtered $(\varphi,N,F/F',E)$-module} of rank
$n$ is
a free $(F_0 \otimes_{\Qp} E)$-module $D$ of rank $n$ equipped with
\begin{itemize}
\item an $F_0$-semilinear, $E$-linear automorphism $\varphi$,
\item a nilpotent $(F_0 \otimes_{\Qp} E)$-linear endomorphism $N$ such that
$N\varphi = p\varphi N$,
\item a decreasing filtration on $D_F = F \otimes_{F_0} D$ such that $\Fil^i D_F$
is zero if $i \gg 0$ and is equal to $D_F$ if $i \ll 0$, and
\item an $F_0$-semilinear, $E$-linear action of $\Gal(F/F')$ commuting
with $\varphi$ and $N$ and preserving the filtration.
\end{itemize}
\end{defn}

Suppose that $\rho : G_{F'} \rightarrow \GL(V)$ is a potentially
semistable representation of $G_{F'}$ on an $n$-dimensional $E$-vector
space $V$, such
that $\rho \,|_{G_F}$ is semistable.  Then
$$ \DstF(V) = (B_{\mathrm{st}} \otimes V)^{G_F} $$
is an example of a filtered $(\varphi,N,F/F',E)$-module of rank $n$.  For
instance, to
see that the action of $\Gal(F/F')$ preserves the filtration, we note that
the filtration is induced from the map
$$ F \otimes_{F_0} \DstF(V) \rightarrow B_{\mathrm{dR}} \otimes V .$$
This map is Galois-equivariant because the inclusion
$B_{\mathrm{st}} \rightarrow B_{\mathrm{dR}}$ is; since the action
of Galois preserves the filtration on $B_{\mathrm{dR}}$, it also
preserves the filtration on $F \otimes_{F_0} \DstF(V)$. Note that in
the special case $F' = \Qp$, Lemma \ref{freeness} implies that the
filtration consists of free $(F \otimes_{\Qp} E)$-modules; this is not
always the case (see \cite[Rem 3.1.1.4]{BreuilMezard}).

\begin{defn} A filtered $(\varphi, N, F/F', E)$-module is said to be
\textit{weakly admissible} if the underlying
$(\varphi,N,F,E)$-module is weakly admissible in the sense of
\cite[D\'ef 3.1.1.1(ii)]{BreuilMezard} (i.e., if one
forgets about the $\Gal(F/F')$-action).
\end{defn}

Therefore $\DstF$ is a functor from the category of
$E$-representations of $G_{F'}$ which become semistable when
restricted to $G_{F}$, to the category of weakly admissible
$(\varphi,N,F/F',E)$-modules.  Conversely, if $D$ is a weakly
admissible $(\varphi,N,F/F',E)$-module, define
$$ \VstF(D) = (B_{\mathrm{st}} \otimes_{F_0} D)^{\varphi=1}_{N=0} \cap \Fil^0
(B_{\mathrm{dR}} \otimes_F (F \otimes_{F_0} D )) .$$ This is an
$E$-representation of $G_{F'}$, where $G_{F'}$ acts as usual on
$B_{\mathrm{st}}$ and through $\Gal(F/F')$ on $D$; moreover, by the
results of \cite{ColmezFontaine} we know that the restriction of
$\VstF(D)$ to $G_{F}$ is a semi-stable representation of dimension
(over $E$) equal to $\text{rk}_{F_0 \otimes_{\Qp} E} D$.

\begin{prop} $\DstF$ and $\VstF$ are quasi-inverses.
\end{prop}

\begin{proof}  Consider the natural $G_{F'}$-homomorphism
$$ B_{\mathrm{st}} \otimes \VstF(D) \rightarrow B_{\mathrm{st}} \otimes D .$$
Taking $G_F$-invariants yields a map
$$ \DstF(\VstF(D)) \rightarrow (B_{\mathrm{st}} \otimes D)^{G_F} = D,$$
which we know must be an isomorphism of underlying
$(\varphi,N,F,E)$-modules. Since our first map was actually a
$G_{F'}$-homomorphism, this isomorphism respects the action of
$\Gal(F/F')$  and so is an isomorphism of
$(\varphi,N,F/F',E)$-modules as well.

The argument for the map $V \rightarrow \VstF(\DstF(V))$ is analogous.
\end{proof}

\begin{cor}
The category of $E$-representations of $G_{F'}$ which become
semistable when restricted to $G_{F}$ and the category of weakly
admissible $(\varphi,N,F/F',E)$-modules are equivalent.
\end{cor}

Following \cite{BreuilMezard}, we make use of functors $\DstkF$ and
$\VstkF$, defined as follows:
$$ \VstkF(D) = (B_{\mathrm{st}} \otimes_{F_0} D)^{\varphi=p^{k-1}}_{N=0} \cap
\Fil^{k-1} (B_{\mathrm{dR}} \otimes_F (F \otimes_{F_0} D ))$$ and
$$ \DstkF(V) = \DstF(V(1-k)) ,$$
where $V(1-k)$ denotes the Tate twist $V \otimes_{\Zp} \Zp(1-k)$.
That these functors are quasi-inverse to one another follows from
the next lemma.

\begin{lemma} For all filtered $(\varphi,N,F/F',E)$-modules $D$, there is
an isomorphism $\VstkF(D) \cong \VstF(D)(k-1)$.
\end{lemma}

The proof of this lemma is exactly the same as the proof of
\cite[Lem. 3.1.1.2]{BreuilMezard}, and similarly we have the following
immediate corollary.

\begin{cor} The functor $\VstkF$ is an equivalence of
categories between the category of weakly admissible filtered
$(\varphi,N,F/F',E)$-modules $D$ such that $\Fil^0(F \otimes_{F_0}
D) = F \otimes_{F_0} D$ and $\Fil^k(F \otimes_{F_0} D)=0$, and the
category of $E$-representations of $G_{F'}$ which are semistable
when restricted to $G_F$ and have Hodge-Tate weights in the range
$\{0,\ldots,k-1\}$.
\end{cor}

\begin{example}  \label{characters} Let $\epsilon$ denote the
cyclotomic character of $G_{\Qp}$, let $\om$ denote the
Teichm\"uller lift of the mod $p$ reduction of $\epsilon$, and for
$a \in \OO_E^{\times}$ let $\lambda_a$ denote the unramified
character of $G_{\Qp}$ sending arithmetic Frobenius to $a$.   Set
$F_1=\Qp(\zeta_p)$. Then $\epsilon^i \om^j \lambda_a$ becomes
semistable when restricted to $G_{F_1}$, and
$D_{\mathrm{st},k}^{F_1}(\epsilon^i \om^j \lambda_a)$ is a
$1$-dimensional filtered module $E\cdot \e$ satisfying $N=0$,
$$\varphi(\e) = p^{k-i-1} a^{-1} \e ,$$ and, for $g \in
\Gal(F_1/\Qp)$, $$ g(\e) = \om^j(g)(\e).$$  Indeed, this
follows directly from the result in the special cases $\epsilon$,
$\om$, and $\lambda_a$. For the first two, use that the element $t
\in B_{\mathrm{st}}$ is a period for $\epsilon$, and that $\om
\,|_{G_{F_1}}$ is trivial, respectively. For $\lambda_a$, one can
use Hilbert's Theorem 90 for $\Fpbar/\Fp$ and a Hensel-like
approximation argument to show that the p-adic completion
$\widehat{\Qp^{\mathrm{un}}}$ of the maximal unramified extension of
$\Qp$ is a period ring for unramified representations. Then if $\e
=\sum x_i \otimes e_i \in (B_{\mathrm{st}} \otimes E)^{G_{F_1}}$
with the $x_i \in \widehat{\Qp^{\mathrm{un}}}$ we have
$$\varphi(\e) = \sum \Frob(x_i) \otimes e_i
 = a^{-1} \sum \Frob(x_i) \otimes a e_i = a^{-1} \Frob(\e) = a^{-1}
 \e,$$ where $\Frob$ is any representative of arithmetic Frobenius
 in $G_{\Qp}$.
\end{example}

\begin{example} \label{char2} Similarly, let $\varpi$ be a choice
of $(-p)^{1/(p^2-1)}$, set $F_2 = \Qpp(\varpi)$, and suppose $E$ is
a finite extension of $\Qpp$.  Let $\omt : G_{\Qpp} \rightarrow
\OO_E^{\times}$ be the character $\omt(g) = (g\varpi)/\varpi$. Then
the character $\omt^m (\epsilon^i \lambda_a) |_{G_{\Qpp}}$ becomes
semistable when restricted to $G_{F_2}$, and
$D_{\mathrm{st},k}^{F_2}(\omt^m (\epsilon^i \lambda_a)|_{G_{\Qpp}})$
is a module $(\Qpp \otimes_{\Qp} E) \cdot \e$ satisfying $N=0$,
$$\varphi(\e) = p^{k-i-1} (1 \otimes a^{-1}) \e ,$$ and, for $g
\in \Gal(F_2/\Qpp)$, $$ g(\e) = (1 \otimes \omt^m(g))(\e) .$$

Finally, the character $\epsilon^i \om^j \lambda_a$ of $G_{\Qp}$
also becomes semistable over $F_2$, and
$D_{\mathrm{st},k}^{F_2}(\epsilon^i \om^j \lambda_a)$ is a module
$(\Qpp \otimes_{\Qp} E) \cdot \e$ satisfying $N=0$, $$\varphi(\e) =
p^{k-i-1} (1 \otimes a^{-1}) \e,$$ and, for $g \in \Gal(F/\Qp)$,
$$ g(\e) = (1 \otimes \om^j(g))(\e).$$

\end{example}

Let $W_{F'}$ denote the Weil subgroup $W_{F'}$ of $G_{F'}$; if $\F'$
is the residue field of $F'$, recall that
there is a map $\alpha : W_{F'} \rightarrow \Z \subset
\Gal(\overline{\F'}/\F')$
sending arithmetic Frobenius to $1$.  Now, if $D$ is a filtered
$(\varphi,N,F/F',E)$-module, observe
that $W_{F'}$
acts in two different ways on $D$: by
restriction to $\Gal(F/F')$, but also by letting $g \in W_{F'}$ act as
$\varphi^{\alpha(g)}$.  These two actions satisfy the compatibilities
necessary to attach a Weil-Deligne representation of $W_{F'}$ to $D$.

\begin{defn} \label{galtype} Suppose that $E$ contains $F_0$.  If $V$ is
an $E$-representation of $G_{\Qp}$ which is semistable when
restricted to $G_F$, then the Weil-Deligne representation $\WD(V)$
attached to $V$ is $\WD(\DstF(V))$.  The \textit{Galois type} (or
\textit{type}) $\tau(V)$ of $V$ is defined to be $\WD(V) \,|_{I_p}$.
\end{defn}

From \cite[App. B.2]{CDT} we recall several properties of $\WD(V)$:
\begin{itemize}
\item $\WD(V)$ does not depend on the choice of $F$, \item
$\WD(V_1 \otimes V_2) \cong \WD(V_1) \otimes \WD(V_2)$, and \item
$\WD(\epsilon_F)$ is unramified, where $\epsilon_F$ denotes the
cyclotomic character of $G_F$.
\end{itemize}
In particular, these facts imply that we could equivalently have defined
$\tau(V)$
using $\WD(\DstkF(V))$ instead of
$\WD(\DstF(V))$

\begin{example} \label{rednotation} Let $\omega$ and $\omega_2$ denote the mod~$p$
reductions of $\om$ and $\omt$.  By abuse of notation, we
often refer to the restrictions $\omega \,|_{I_p}, \omega_2
\,|_{I_p},\om \,|_{I_p},\omt \,|_{I_p}$ simply as $\omega,
\omega_2, \om, \omt$. Suppose that $V$ is a $2$-dimensional
potentially semistable $E$-representation of $G_{\Qp}$ and suppose that
$\tau(V)$ is tamely ramified. Then either $\tau(V) = \om^i \oplus
\om^j$ for integers $i$ and $j$ (the ``principal series'' case) or
else $\tau(V) = \omt^m \oplus \omt^{pm}$ for an integer $m$ not
divisible by $p+1$ (the ``supercuspidal'' case).
\end{example}

\begin{prop}  \label{prin} Suppose that $V$ is an
indecomposable $2$-dimensional potentially semistable
$E$-representation of $G_{\Qp}$ with
\begin{itemize}
\item Hodge-Tate weights $(0,1)$, and\item type $\tau(V) = \om^i
\oplus \om^j$ for $i \not\equiv j
\pmod{p-1}$.
\end{itemize}
Let $\pi$ be a choice of $(-p)^{1/(p-1)}$ and set $F_1 = \Qp(\pi)$.
Then $V$ is crystalline over $F_1$ and $D_{\mathrm{st},2}^{F_1}(V)$
is of the form
$$ D = E\cdot \e_1 \oplus E \cdot \e_2 ,$$
$$ \varphi(\e_1) = x_1 \e_1 \ , \quad \varphi(\e_2) = x_2 \e_2 \ , \quad \ N = 0 ,$$
$$ \Fil^1(F_1 \otimes_{\Qp} D) = (F_1 \otimes_{\Qp} E)( \pi^{j-i} \e_1 + \e_2 )
,$$
$$ g \cdot \e_1 = \om(g)^i \e_1 \ , \quad g \cdot \e_2 = \om(g)^j \e_2 \
\text{for} \ g \in \Gal(F_1/\Qp) ,$$ with $x_1,x_2 \in \OO_E$
and $\val_p(x_1 x_2) = 1$.
\end{prop}

\begin{proof}
Since $\tau(V)$ is nonscalar, $V$ is potentially crystalline by
\cite[Lem. 2.2.2.2]{BreuilMezard}; moreover, $\tau(V)
\,|_{I_{F_1}}$ is trivial and so $V$ becomes crystalline over
$F_1$. Hence $N=0$.

From the construction of $\WD(V)$, it is easy to see that $D$ must
have a basis $\e_1,\e_2$ on which $\Gal(F_1/\Qp)$ acts via $g
\cdot \e_1 = \om(g)^i \e_1$ and $g \cdot \e_2 = \om(g)^j \e_2$.
Since $\varphi$ and $g$ commute, and again using the fact that
$\tau(V)$ is nonscalar, it follows that $\varphi(\e_1)=x_1 \e_1$
and $\varphi(\e_2)=x_2 \e_2$ for some $x_1$ and $x_2$.

Using the fact that $\Gal(F_1/\Qp)$ preserves the filtration, we
find $\Fil^1(F_1 \otimes_{\Qp} D)$ to be of the form $(F_1
\otimes_{\Qp} E)(\pi^{j-i} a \e_1 + b \e_2)$ for $a,b \in E$. Both
$a$ and $b$ must be nonzero: otherwise, the resulting
$(\varphi,N,F_1/\Qp,E)$-module would be a direct sum of two
$1$-dimensional $(\varphi,N,F_1/\Qp,E)$-modules, contradicting the
indecomposability of $V$. Replacing $\e_1$ by $a \e_1$ and $\e_2$
by $b \e_2$ in our basis for $D$, we see that $\Fil^1$ may be
taken to have the desired form. Finally, the weak admissibility of
$D$ implies that $x_1,x_2 \in \OO_E$ and that $\val_p(x_1 x_2)=1$.
\end{proof}

We denote the filtered modules of the preceding Proposition
by $D_{x_1,x_2}$.

\begin{prop}  \label{super} Suppose that $E$ contains $\Qpp$, and suppose
that $V$ is a
$2$-dimensional
potentially semistable $E$-representation of $G_{\Qp}$ with
\begin{itemize}
\item Hodge-Tate weights $(0,1)$, and
\item type $\tau(V) = \omt^m \oplus \omt^{pm}$ with $p+1 \nmid m$.
\end{itemize}
Write $m = i + (p+1)j$ with $i \in \{1,\ldots,p\}$ and $j \in
\Z/(p-1)\Z$.  Let $\varpi$ be a choice of $(-p)^{1/(p^2-1)}$, set
$F_2 = \Qpp(\varpi)$, and let $g_{\varphi}$ denote the element of
$\Gal(F_2/\Qp)$ which fixes $\varpi$ and is nontrivial on $\Qpp$.
Then $V$ is crystalline over $F_2$ and $D_{\mathrm{st},2}^{F_2}(V)$
is of the form
$$ D = (\Qpp \otimes E)\cdot \e_1 \oplus (\Qpp \otimes E) \cdot \e_2 ,$$
$$ \varphi(\e_1) = \e_2 \ , \quad \varphi(\e_2) = (1 \otimes x) \e_1 \ , \quad N = 0 ,$$
$$ \Fil^1(F_2 \otimes_{\Qpp} D) = (F_2 \otimes_{\Qp} E)( (\varpi^{(p-1)i}\otimes
a) \e_1 + (1 \otimes b ) \e_2 ) ,$$
$$ g \cdot \e_1 = (\omt(g)^m \otimes 1) \e_1 \ , \quad g \cdot \e_2 =
(\omt(g)^{pm} \otimes 1) \e_2 \ \text{for} \ g \in \Gal(F_2/\Qpp)
,$$
$$ g_{\varphi} \cdot \e_1 = \e_1 \ , \quad g_{\varphi} \cdot \e_2 = \e_2 ,$$
with $(a,b) \in E^2 \setminus (0,0)$, $x \in \OO_E$ and $\val_p(x)
= 1$.
\end{prop}

\begin{proof}
Exactly as in Proposition \ref{prin}, $V$ becomes crystalline over
$F_2$ and $N=0$. Let $\sigma_1$ and $\sigma_2$ denote the two
embeddings of $\Qpp$ into $E$.  For each $\mu=1,2$, the construction
of $\WD(V)$ implies that $D_{\sigma_\mu}$ has an $E$-basis $v_{\mu
1},v_{\mu 2}$ on which $g \in \Gal(F_2/\Qpp)$ acts as
$$g \cdot v_{\mu 1} = \sigma_\mu(\omt^{m}(g)) v_{\mu 1} \ \text{and} \ g \cdot
v_{\mu 2} = \sigma_\mu (\omt^{pm}(g)) v_{\mu 2} .$$ Since
$g_{\varphi}$ is an $E$-linear map on $D$ which swaps the two
subspaces $D_{\sigma_\mu}$ and satisfies the relation $g_{\varphi} g
g_{\varphi} = g^p$ for all $g \in \Gal(F_2/\Qpp)$, it follows
(possibly after multiplying some $v_{\mu \nu}$ by constants in $E$)
that $g_{\varphi} \cdot v_{\mu \nu} = v_{(3-\mu)\nu}$. Similarly,
since $\varphi$ swaps the $D_{\sigma_\mu}$ and commutes with the
action of $\Gal(F_2/\Qpp)$, there exist $c,d \in E$ such that
$\varphi(v_{11})=c v_{22}$, $\varphi(v_{21})=c v_{12}$,
$\varphi(v_{12})=d v_{21}$, and $\varphi(v_{22}) = d v_{11}$.

Taking $\e_1 = v_{11} + v_{21}$, $\e_2 = c(v_{12} + v_{22})$, and
$x=cd$, and using the fact that $\Fil^1$ must be preserved by
$\Gal(F_2/\Qp)$, we see without difficulty that in this basis $D$
has the desired form.
\end{proof}

We denote the filtered modules of Proposition \ref{super} by
$D_{m,[a:b]}$.

\begin{remark}  It is not difficult to see that these filtered modules
match those of ``type IV'' in \cite[\S 11]{FontaineMazur} (which deals
only with the case $E=\Qp$).
\end{remark}

By a similar argument, we also find the following.

\begin{prop}  \label{prin2} Suppose that $E$ contains $\Qpp$, and suppose that $V$ is
an indecomposable $2$-dimensional potentially semistable
$E$-representation of $G_{\Qp}$ with
\begin{itemize}
\item Hodge-Tate weights $(0,1)$, and
\item type $\tau(V) = \om^i \oplus \om^j$ for $i \not\equiv j \pmod{p-1}$.
\end{itemize}
Let $\varpi$, $F_2$, and the elements of $\Gal(F_2/\Qp)$ be as in
Proposition \ref{super}, and set $\pi = \varpi^{p+1}$. Then $V$ is
crystalline over $F_2$ and $D_{\mathrm{st},2}^{F_2}(V)$ is of the
form
$$ D = (\Qpp \otimes E)\cdot \e_1 \oplus (\Qpp \otimes E) \cdot \e_2 ,$$
$$ \varphi(\e_1) = (1 \otimes x_1) \e_1 \ , \quad \varphi(\e_2) = (1 \otimes x_2) \e_2 \ , \quad N = 0 ,$$
$$ \Fil^1(F_2 \otimes_{\Qpp} D) = (F_2 \otimes_{\Qp} E)(( \pi^{j-i} \otimes 1) \e_1 + \e_2 ) ,$$
$$ g \cdot \e_1 = (\om(g)^i \otimes 1) \e_1 \ , \quad g \cdot \e_2 =
(\om(g)^{j} \otimes 1) \e_2 \ \text{for} \ g \in \Gal(F_2/\Qpp)
,$$
$$ g_{\varphi} \cdot \e_1 = \e_1 \ , \quad g_{\varphi} \cdot \e_2 = \e_2 ,$$
with $x_1, x_2 \in \OO_E$ and $\val_p(x_1x_2) = 1$.
\end{prop}

Denote this filtered module by $D'_{x_1,x_2}$, and note that
$(D'_{x_1,x_2})^{\Gal(F_2/F_1)} = D_{x_1,x_2}$.  When $0 <
\val_p(x_1),\val_p(x_2) < 1$, we actually need to work with
$D'_{x_1,x_2}$ instead of $D_{x_1,x_2}$.

Finally, we conclude this section by checking the following.

\begin{prop} \label{weak-acc} The filtered modules listed in Propositions \ref{prin},
\ref{super}, \ref{prin2} are weakly admissible.
\end{prop}

\begin{proof} To verify the weak admissibility of one of these filtered
modules $D$, by \cite[Prop. 3.1.1.5]{BreuilMezard} one needs
to check that any $(F_0 \otimes E)$-submodule of $D$ that is preserved
by $\varphi$ satisfies $t_H \le t_N$ (in the notation of Fontaine
\cite{AST223}).  For the filtered modules of Proposition \ref{prin},
this is evident.  We give the details for the filtered modules of
Proposition \ref{super}; the details for Proposition \ref{prin2} are
similar.

Suppose that $\Qpp = \Qp(\tau)$ with $\tau^2 \in \Qp$.  Put $u =
\frac{1}{2}(1\otimes 1 - \tau \otimes \tau^{-1}) \in \Qpp \otimes E$
and set $v = 1\otimes 1 - u$, so that $u,v$ are idempotents
satisfying $\varphi(u) = v$, $\varphi(v)=u$.  It then easy to see
that any nonzero proper $\Qpp \otimes E$-submodule $D_0$ of $D$
which is preserved by $\varphi$ must be generated by a pair of
elements $u(a \e_1 + b\e_2)$, $v(bx \e_1 + a \e_2)$ with $a,b \in
E$: if $D_0$ contains two $E$-linearly independent elements of $uD$,
then it contains all of $uD$; hence it is equal to $D$.  One checks
that $t_H(D_0) = 0$ for such $D_0$.
\end{proof}

\section{Strongly divisible modules with tame descent data}
\label{strtame}

 In \cite{Br}, C. Breuil constructed a category of ``strongly
divisible modules'' over a local field $F'$, and he proved that it is
(anti-)equivalent to the category of Galois lattices inside
Barsotti-Tate representations of $G_{F'}$.  Following the strategy
of \cite[Sec. 5.4]{BCDT}, we formulate descent data on
strongly divisible modules, thereby extending Breuil's
antiequivalence to Galois lattices inside potentially
Barsotti-Tate representations. This is essentially formal, but for
simplicity we  work exclusively with descent data for tame
extensions $F/F'$.

\subsection{Galois lattices and p-divisible groups}

In this section, we review the relation between Galois lattices
and $p$-divisible groups.  Let $\rho : G_{F'} \rightarrow \GL(V)$
be a $p$-adic representation, and let $T$ be a Galois
$\Zp$-lattice inside $V$.  To $T$ we may associate a $p$-divisible
group over $F'$, as follows: each $T/p^n$ is a finite
representation of $G_{F'}$, hence corresponds to a finite flat
group scheme $\Gamma(n)$ over $F'$.  Then the $p$-divisible group
associated to $T$ is $\Gamma = \cup \Gamma(n)$.  Conversely, given
a $p$-divisible group $\Gamma$ over $F'$ we may recover the Galois
lattice
$$\underset{n}{\invlim} \, \Gamma(n)(\Qpbar),$$
and these two operations are readily seen to be inverse to one
another.

Breuil shows (\cite[Th. 5.3.2]{Br}) that $\Gamma$ extends
to a $p$-divisible group $\G$ over the integers $\OO_{F'}$ if and
only if the representation $\rho$ is crystalline with Hodge-Tate
weights in $\{0,1\}$.  More precisely, Breuil shows that if $\rho$
is crystalline with Hodge-Tate weights in $\{0,1\}$, then there
exists some lattice inside $V$ for which the associated
$p$-divisible group over $F'$ extends; but then by a
scheme-theoretic closure argument (see
\cite[Secs. 2.2, 2.3]{Raynaud}), for any lattice $T \subset V$ the $p$-divisible
group over $F'$ associated to $T$ extends to a $p$-divisible
group over $\OO_{F'}$.  Tate's full faithfulness theorem guarantees
that this extension is unique up to isomorphism.

Suppose now that $\rho$ is merely \textit{potentially} crystalline
with Hodge-Tate weights in $\{0,1\}$, and, more precisely, that $\rho$
becomes crystalline over $F$.  Let $T \subset V$ be a Galois
lattice.  Then $T$ regarded as a $G_F$-lattice does correspond to a
$p$-divisible group $\Gamma$ over $F$ which, as above, extends to a
$p$-divisible group $\G$ over $\OO_F$. However, the restriction from
$F'$ to $F$ also induces descent data on $\Gamma$.  Indeed, recall
that
$$\Gamma(n) = \Spec({\rm Maps}_{G_{F}}(T/p^n, \Qpbar)) .$$
The algebra on the right-hand side carries an action of
$\Gal(F/F')$: if $g \in \Gal(F/F')$, let $\tilde{g}$ be any
extension of $g$ to $\Gal(\Qpbar/F')$, and for $f \in {\rm
Maps}_{G_{F'}}(T/p^n, \Qpbar))$, we set $g \cdot f = \tilde{g} \circ f
\circ \tilde{g}^{-1}$.  This is easily seen to be well-defined and
compatible among different values of $n$, so that we obtain a
$g$-semilinear map $\langle g \rangle : \Gamma \rightarrow
\Gamma$.  It is convenient to factor $\langle g \rangle$ as
$$
\begin{CD}
\Gamma @>[g]>> {}^g \Gamma @>>> \Gamma\\
@VVV @VVV @VVV \\
\Spec(F) @>>> \Spec(F) @>g>> \Spec(F)
\end{CD}
$$
where the right-hand square is Cartesian, so that the $[g]$ are
maps of $p$-divisible groups over $F$ satisfying the compatibility
$[gh] = ({}^g[h]) \circ [g]$.  (Here and henceforth, the
superscript ${}^g$ denotes base change by $g$.)  Finally, by
Tate's full faithfulness theorem, each $[g]$ extends to a map $\G
\rightarrow {}^g\G$.  We again denote this by $[g]$, and we
note that the compatibility relation is automatically still
satisfied.

\begin{defn}  If $\G$ is a $p$-divisible group over $\OO_F$, then
\textit{descent data relative to $F'$} is a collection of maps
$[g] : \G \rightarrow {}^g\G$ for each $g \in \Gal(F'/F)$
satisfying $[gh] = ({}^g[h]) \circ [g]$.
\end{defn}

In the reverse direction, if $\G = \cup \G(n)$ is a $p$-divisible
group over $\OO_F$ with descent data relative to $F'$, we can
construct a $G_{F'}$-lattice.  Writing $\G(n)=\Spec(R_n)$, the
descent data comes from a compatible collection of
$\Gal(F/F')$-actions on the $R_n$. If $\sigma \in G_{F'}$, and $f
\in \G(n)(\OO_{\Qpbar}) = \Hom(R_n,\OO_{\Qpbar})$, we set $\sigma
\cdot f = \sigma \circ f \circ (\sigma ^{-1} \,|_{F})$.  Since the
descent data is actually descent data on the whole $p$-divisible
group, these $G_{F'}$-actions are compatible for varying $n$ and
yield a $G_{F'}$-action on
$$\underset{n}{\invlim} \, \G(n)(\OO_{\Qpbar}) .$$
Unsurprisingly, this construction is inverse to the construction
from a $G_{F'}$-lattice of a $p$-divisible group over $\OO_{F}$
with descent data relative to $F'$.  As a consequence, we have
the following.

\begin{prop} \label{prop:lattice} The above constructions describe an equivalence
between the category of Galois lattices inside potentially
crystalline $G_{F'}$-representations that become crystalline over
$F$ and have Hodge-Tate weights inside $\{0,1\}$, and the category
of $p$-divisible groups over $\OO_F$ with descent data relative to
$F'$.
\end{prop}

\subsection{Big rings and categories of filtered modules without descent data}
\label{modules}

In this section we review the definitions of various big rings and
categories of filtered modules from \cite{Br} and
\cite{BreuilStrong}.

Let $R$ be a complete discrete valuation ring of characteristic
zero, absolute ramification index $e$, and perfect residue field $k$
of characteristic $p$. Fix a uniformizer $\pi$ of $R$. Let $S=S_R$
be the $p$-adic completion of $W(k)[u,\frac{u^{ie}}{i!}]_{i \in
\N}$, and let $S_n = S/p^nS$. The map $\phi : S \rightarrow S$ is
the unique Frobenius-semilinear map sending $\phi(u)=u^p$ and
$\phi(u^{ie}/i!) = u^{iep}/i!$; we also use $\phi$ to denote
the map $S_n \rightarrow S_n$ induced by $\phi$. Let $N$ denote the
unique $W(k)$-linear derivation such that $N(u)=-u$ and
$N(u^{ie}/i!) = -ieu^{ie}/i!$, so that $N\phi = p\phi N$. Let $E(u)
\in S$ denote the minimal polynomial of $\pi$ over $W(k)$, and if $k
\ge 1$, let $\Fil^{k-1} S$ be the $p$-adic completion of the ideal of
$S$ generated by $E(u)^i/i!$ for $i \ge k-1$. (Note that we are
using $k$ to denote both the residue field and the weight; it should
not be possible to confuse the two uses with one another.) Then
$\phi(\Fil^{k-1} S) \subset p^{k-1}S$ for $k \le p$, and so for $k
\le p$, we let $\phi_{k-1}$ denote $\phi/p^{k-1}$ on $\Fil^{k-1} S$.
Finally, let $c$ denote $\p(E(u))$.

We now repeat (essentially verbatim) some notation and definitions
of \cite{BCDT} and \cite{Br}; we refer the reader to \cite[Secs. 5.3, 5.4]{BCDT} for details.

Let $\syn$ be the small $p$-adic formal syntomic site over $R$,
and let $\AbR$ denote the category of abelian sheaves on $\syn$.

If $\X \in {\rm Spf}(R)_{{\rm syn}}$, set $\X_n = \X \times_R
R/p^n$. The sheaf $\Ocris$ is the sheaf of $S_n$-modules on ${\rm
Spf}(R)_{{\rm syn}}$ associated to the presheaf
$$ \X \mapsto \left( W_n(k)[u] \otimes_{\phi^n,W_n(k)}
W_n(\Gamma(\X_1,\OOO_{\X_1}))\right)^{DP}, $$ where $\phi$ is
Frobenius on $W_n(k)$, where ``DP'' means that we take the divided
power envelope with respect to the kernel of the map
$$ W_n(k)[u] \otimes_{\phi^n,W_n(k)} W_n(\Gamma(\X_1,\OOO_{\X_1}))
\rightarrow \Gamma(\X_n, \OOO_{\X_n})$$
$$ s(u) \otimes (w_0,\ldots,w_{n-1}) \mapsto s(\pi)(\hat{w}_0^{p^n} + \cdots + p^{n-1} \hat{w}_{n-1}^p)$$
and relative to the usual divided power structure on the maximal
ideal of $W_n(k)$, and where $\hat{w}_i$ is a local lifting of
$w_i$. If $\OOO_n \in \AbR$ is the sheaf
$\OOO_n(\X)=\Gamma(\X_n,\OOO_{\X_n})$, then the above map induces a
morphism $\Ocris \rightarrow \OOO_n$, and we denote its kernel by
$\Jcris$.  The map $\phi : \Ocris \rightarrow \Ocris$ induced by
crystalline Frobenius satisfies $\phi(\Jcris) \subset p\Ocris$, and
there exists $\p : \Jcris \rightarrow \Ocris$ which may be thought of
as $\phi/p$. Define $\Ocrisi = \underset{n}{\varinjlim} \Ocris$ and
$\Jcrisi = \underset{n}{\varinjlim} \Jcris$; these limits are taken
over the multiplication-by-$p$ inclusions $\Ocris \rightarrow
\OOO^{\mathrm{cris}}_{n+1,\pi}$. See \cite[Sec 2.3]{Br} for
further details regarding these sheaves.

If $R=\OO_F$ is the ring of integers in a finite extension $F$ over
$\Qp$, recall from \cite[Sec 5.3]{Br} that we define
$A_{\mathrm{cris}} = \invlim W_n(\OO_{\Qpbar}/p\OO_{\Qpbar})^{DP}$.
Fix a system of roots $(\pi_n)_{n \ge 0}$ in $\Qpbar$ such that
$\pi_0=\pi$ and $\pi_{n+1}^p = \pi_n$, from which we construct an
element $\upi \in \Acris$ (see
\cite[Sec. 2.2.2]{BreuilStrong}).  Then $B_{\mathrm{cris}}^{+}=A_{\mathrm{cris}}
\otimes_{W(k)} F_0$, where $F_0$ is the fraction field of $W(k)$, and
$\Acrish = \invlim \Ocris(\OO_{\Qpbar})$ is isomorphic to the
$p$-adic completion of
$A_{\mathrm{cris}}[\frac{(u-\underline{\pi})^i}{i!}]_{i \in \N}$.

We refer the reader to \cite[Sec. 2.2.2]{BreuilStrong} for the
construction of the ring $\Asth$.  The ring $\Asth$ has a filtration
$\Fil^{\bullet} \Asth$, a Frobenius $\phi$, and a monodromy operator
$N$ which is the unique $A_{\mathrm{cris}}$-linear derivation such
that $N(X)=1+X$.  If $k \le p$, the Frobenius satisfies
$\phi(\Fil^{k-1}\Asth)\subset p^{k-1}\Asth$, so we let $\phi_{k-1}$
be $\phi/p^{k-1}$ on $\Fil^{k-1} \Asth.$ The choice of $\upi$ fixes
an $S$-module structure on $\Asth$ and an embedding $\Acrish
\rightarrow \Asth$ by sending $u \mapsto \underline{\pi}(1+X)^{-1}$,
and this embedding induces a filtration, Frobenius, and monodromy
operator on $\Acrish$ and a filtration on $\Acris$.  Set $\Acrisi =
\Acris \otimes_{W(k)} F_0/W(k)$, $\Acrishi = \Acrish \otimes_{W(k)}
F_0/W(k)$,  and $\Asthi = \Asth \otimes_{W(k)} F_0/W(k)$ with the
induced Frobenius, filtration, and monodromy operators (e.g., 
$\Fil^{k-1} \Asthi = (\Fil^{k-1} \Asth) \otimes F_0/W(k)$).

Let $k \in \{1,\ldots,p-1\}$.  Recalling
\cite[Sec. 2.2.1]{BreuilStrong}, we let $'\uMod^{k-1}$ denote the category of
quadruples consisting of
\begin{itemize}
\item an $S$-module $\M$,

\item an $S$-submodule $\Fil^{k-1} \M$ of $\M$ containing
$(\Fil^{k-1} S)\M$,

\item a $\phi$-semilinear map $\phi_{k-1} : \Fil^{k-1} \M
\rightarrow \M$ such that for all $s \in \Fil^{k-1} S$ and $x \in
\M$ we have $\phi_{k-1}(sx) = \phi_{k-1}(s)\phi(x)$ with
$\phi(x)=\frac{1}{c^{k-1}} \phi_{k-1}(E(u)^{k-1} x)$, and

\item a $W(k)$-linear map $N : \M \rightarrow \M$ satisfying:
$N(sx)=N(s)x + sN(x)$ for $s \in S, x \in \M$; and
$E(u)N(\Fil^{k-1} \M) \subset \Fil^{k-1} \M$ and
$\phi_{k-1}(E(u)N(x)) = cN(\phi_{k-1}(x))$ for $x \in \Fil^{k-1}
\M$.
\end{itemize}

Morphisms in $'\uMod^{k-1}$ are the $S$-linear maps preserving
$\Fil^{k-1}$ and commuting with $\phi_{k-1}$ and $N$.  We define
six additional categories as follows: $'\uMod^{k-1}_0$ is the
category obtained by omitting $N$ in the definition of
$'\uMod^{k-1}$, while $\uMod^{k-1}$ and $\uMod^{k-1}_0$ are the
full subcategories of $'\uMod^{k-1}$ and $'\uMod^{k-1}_0$ with the
following extra conditions:
\begin{itemize}
\item $\M$ is of the form $\oplus_i S_{n_i}$ for some finite list
of positive integers $n_i$, and

\item $\phi_{k-1}(\Fil^{k-1}\M)$ generates $\M$ over $S$.
\end{itemize}

Next, $\Mod^{k-1}$ and $\Mod^{k-1}_0$ are the full subcategories
of $'\uMod^{k-1}$ and $'\uMod^{k-1}_0$ with the following extra
conditions:
\begin{itemize}
\item $\M$ is a free $S$-module and $\Fil^{k-1} \M \cap p\M =
p\Fil^{k-1}\M$, and

\item $\phi_{k-1}(\Fil^{k-1}\M)$ generates $\M$ over $S$.
\end{itemize}

Finally, let $\Mod^{k-1}_{\mathrm{cris}}$ be the full subcategory of
objects of $\Mod^{k-1}$ with the property that $N(\M) \subset I\M$,
where $I$ is the ideal $\sum_{i \ge 1} \frac{u^i}{\lfloor i/e
\rfloor!} S$ in $S$.

The category $\Mod^{k-1}_0$ is called the category of strongly
divisible modules (of weight $k$).  Let $R=\OO_F$ be the integers
in a finite extension $F$ of $\Qp$, and $\M$ be a strongly
divisible module of weight $2$ for $R$. By
\cite[Prop. 5.1.3(1)]{Br}, there exists a unique $W(k)$-linear
endomorphism $N$ of $\M$ such that:
\begin{itemize}
\item $N(sx) = N(s)x + sN(x)$ for $s \in S$ and $x \in \M$,

\item $N\p = \phi N$, and

\item $N(\M) \subset I \M$,  where $I$ is the ideal $\sum_{i \ge 1}
\frac{u^i}{\lfloor i/e \rfloor !}S$ in $S$.
\end{itemize}
Thus for $R=\OO_F$, the categories $\Mod^{1}_0$ and
$\Mod^{1}_{\mathrm{cris}}$ are equivalent. Before proceeding to the
next section, we note the following examples:
\begin{itemize}
\item $S$ is an object of $\Mod^{k-1}_{\mathrm{cris}}$,

\item each $S_n$ is an object of $\uMod^{k-1}$,

\item $\Acrish$, $\Acrishi$, $\Asth$, and $\Asthi$ are objects of
$'\uMod^{k-1}$,

\item $\Ocris(\OO_{\Qpbar})$ is an object of $'\uMod^{1}_0$, and

\item regarding $\Acris$ as an $S$-module via $u \cdot x = \upi x$
we can make $\Acris$ and $\Acrisi$ into objects of
$'\uMod^{k-1}_0$; then the maps $\Asth \rightarrow \Acris$ and
$\Asthi \rightarrow \Acrisi$ sending $X \mapsto 0$ are morphisms
in $'\uMod^{k-1}_0$.
\end{itemize}

\subsection{$p$-divisible groups and strongly divisible modules with tame descent data}

Let $\G = \cup \G(n)$ be a $p$-divisible group over $R$.  Then for
each $n$ we may regard $\G(n)$ as a sheaf on $\syn$ and we define
$$\M_{\pi}(\G(n)) = \Hom_{\AbR}(\G(n),\Ocrisi) =
\Hom_{\AbR}(\G(n),\Ocris) ,$$
$$ \Fil^1 \M_{\pi}(\G(n)) = \Hom_{\AbR}(\G(n),\Jcrisi) =
\Hom_{\AbR}(\G(n),\Jcris) ,$$ and $$\p : \Fil^1 \M_{\pi}(\G(n))
\rightarrow \M_{\pi}(\G(n))$$ induced by $\p : \Jcrisi \rightarrow
\Ocrisi$.  Next, define
$$(\M_{\pi}(\G),\Fil^1\M_{\pi}(\G),\p) = (\underset{n}{\invlim} \M_{\pi}(\G(n)),\underset{n}{\invlim} \Fil^1\M_{\pi}(\G(n)),\underset{n}{\invlim} \p).$$
We often denote the triples
$(\M_{\pi}(\cdot),\Fil^1\M_{\pi}(\cdot),\p)$ by
 $\M_{\pi}(\cdot)$ or, suppressing the fixed
uniformizer $\pi$, by $\M(\cdot)$.  By
\cite[Cor. 4.2.2.7, Lem. 4.2.2.8]{Br}, the $\M(\G(n))$ are objects in $\uMod^1_0$,
while $\M(\G)$ is an object in the category $\Mod^1_0$ of strongly
divisible modules.

Now suppose that $g : R \rightarrow R$ is a continuous
automorphism of $R$. For simplicity of notation, we assume
that $g(\pi) = h_g \pi$ with $h_g \in W(k)$. This assumption is
not strictly necessary until Corollary \ref{prop:dd}, but we do
not need the extra generality for our applications.  Define $\ghat
: W(k)[[u]] \rightarrow W(k)[[u]]$ by $\ghat(\sum w_i u^i) = \sum
g(w_i) h_g^i u^i$, and similarly, let $\ghat : S \rightarrow S$ be
the unique ring isomorphism such that $\ghat\left(w_i
\frac{u^{i}}{\lfloor i/e \rfloor!}\right) = g(w_i)
\frac{u^{i}}{\lfloor i/e \rfloor!} h_g^{i} $.  We also let $\ghat$
denote the isomorphism induced on $S_n$.

If $\X \in \syn$, let ${}^g \X = {\rm Spf}(R) \times_{g^{*},{\rm
Spf(R)}} \X$, and as in \cite[Sec. 5.4]{BCDT}, we define
$$ \Ocrisg(\X) = \Ocris({}^g\X) , \quad \Jcrisg(\X) = \Jcris({}^g\X),$$
so that $\Ocrisg \in \AbR$ is the sheaf associated the the
presheaf
\begin{eqnarray*}
\X & \mapsto & \left( W_n(k)[u] \otimes_{\phi^n,W_n(k)}
W_n(\Gamma({}^g \X_1, \OOO_{{}^g \X_1}))\right)^{DP} \\
 & = &\left( W_n(k)[u] \otimes_{\phi^n,W_n(k)}
W_n(R \otimes_{g,R} \Gamma(\X_1, \OOO_{\X_1}))\right)^{DP} .
\end{eqnarray*}
Then there is a canonical isomorphism
$$ \Ocris \otimes_{S_n,\ghat} S_n \xrightarrow{\sim} \Ocrisg $$
coming from the $\ghat$-semilinear map from  $W_n(k)[u] \otimes_{\phi^n,W_n(k)}
W_n(\Gamma(\X_1,\OOO_{\X_1}))$ to  $W_n(k)[u]
\otimes_{\phi^n,W_n(k)} W_n(R \otimes_{g,R} \Gamma(\X_1,
\OOO_{\X_1}))$ sending
$$ s \otimes (w_0,\ldots,w_{n-1}) \mapsto \ghat(s) \otimes
(1\otimes w_0,\ldots,1\otimes w_{n-1}),$$ and this induces
$$ \Jcris \otimes_{S_n,\ghat} S_n \xrightarrow{\sim} \Jcrisg .$$
Lemma 5.4.4 of \cite{BCDT} tells us that the diagram $$
\begin{CD}
\Jcris \otimes_{S_n,\ghat} S_n @>\sim>> \Jcrisg \\
@V{\p\otimes\phi}VV  @V{\p}VV \\
\Ocris \otimes_{S_n,\ghat} S_n @>\sim>> \Ocrisg
\end{CD}$$
is commutative.  Moreover, looking at the presheaves, it is evident
that the above diagrams for $n$ and $n+1$ are compatible under the
multiplication-by-$p$ inclusion (i.e., we have a commutative
cube, where the front and back faces are the above diagrams for $n$
and $n+1$, and the four front-to-back maps are induced by $\Ocris
\rightarrow \OOO^{\mathrm{cris}}_{n+1,\pi}$).  We have the following
proposition, which is an analogue of \cite[Cor. 5.4.5]{BCDT}
and is proved in essentially the same manner.

\begin{prop}\label{prop:pdiv}
Let $g : R \rightarrow R$ be a continuous automorphism such that
$g \pi = h_g \pi$ with $h_g \in W(k)$.
\begin{enumerate}
\item Let $\G$ be a $p$-divisible group over $R$.  Then there are
canonical isomorphisms in $\Mod^1_0$:
$$ (\M_{\pi}(\G) \otimes_{\ghat} S, \Fil^1\M_{\pi}(\G)
\otimes_{\ghat} S, \p \otimes \phi) \xrightarrow{\sim}
(\M_{\pi}({}^g \G),\Fil^1\M_{\pi}({}^g\G),\p) .$$

\item If $f: \G \rightarrow \G'$ is a morphism of $p$-divisible
groups over $R$, then there is a commutative diagram in
$\Mod^1_0$:
$$\begin{CD} \M_{\pi}(\G') \otimes_{\ghat} S
@>\M_{\pi}(f)>>
\M_{\pi}(\G) \otimes_{\ghat} S \\
@VVV @VVV \\
\M_{\pi}({}^g\G') @>\M_{\pi}({}^gf)>> \M_{\pi}({}^g\G)
\end{CD} .
$$

\item If $g_1,g_2$ are two continuous automorphisms such that $g_i
\pi = h_{g_i} \pi$ with $h_{g_i} \in W(k)$ for $i=1,2$, then on
$$(\M_{\pi}(\G) \otimes_{\ghat_1} S) \otimes_{\ghat_2} S
\cong \M_{\pi}(\G) \otimes_{\widehat{g_2g_1}} S$$ one has $(\p
\otimes \phi) \otimes \phi = \p \otimes \phi$.

\end{enumerate}
\end{prop}

\begin{proof} First, prove the $p^n$-torsion analogue of this
proposition; more precisely, prove the same statements for $\G(n)$ in
$\uMod^1_0$, and note that they are compatible under the
inclusions $\G(n)\rightarrow\G(n+1)$.  Then pass to the inverse
limit.
\end{proof}

We then have the following analogue of
\cite[Cor. 5.4.6]{BCDT}.

\begin{cor}\label{cor:onemap}  Let $\G$ be a $p$-divisible group over $R$.  To
give a morphism $\langle g \rangle : \G \rightarrow \G$ such that
the diagram
$$\begin{CD}
\G @>\langle g\rangle>> \G \\
@VVV @VVV \\
\Spec(R) @>\Spec(g)>> \Spec(R)
\end{CD}$$
is commutative and the induced morphism $[g] : \G \rightarrow
{}^g\G$ is a morphism of $p$-divisible groups over $R$ is
equivalent to giving an additive map $\ghat: \M_{\pi}(\G)
\rightarrow \M_{\pi}(\G)$ such that
\begin{itemize}
\item for all $s \in S$ and $x \in \M_{\pi}(\G)$,
$\ghat(sx)=\ghat(s)\ghat(x)$;

\item $\ghat(\Fil^1\M(\G)) \subset \Fil^1\M(\G)$ and $\p \circ
\ghat = \ghat \circ \p$.
\end{itemize}
\end{cor}

\begin{proof} As in the proof of \cite[Cor. 5.4.6]{BCDT},
the map $\ghat$ is the composition
$$ \M_{\pi}(\G) \rightarrow \M_{\pi}(\G)\otimes_{\ghat}S
\rightarrow \M_{\pi}({}^g\G) \rightarrow \M_{\pi}(\G) $$ where the
leftmost map is $x \mapsto x \otimes 1$ and the rightmost map is
$\M_{\pi}([g])$.
\end{proof}

We now specify hypotheses that we frquently need to assume.

\begin{hyp}\label{hyp} Suppose $R=\OO_F$, and suppose that
$F/F'$ is a tamely ramified Galois extension with ramification
index $e(F/F')$. Fix a uniformizer $\pi \in F$ such that
$\pi^{e(F/F')} \in F'$, and write $g(\pi)=h_g \pi$ for each $g \in
\Gal(F/F')$.  Let $\ghat: S \rightarrow S$ be defined as before.
\end{hyp}

Then Corollary \ref{cor:onemap} and parts (2) and (3) of Proposition
\ref{prop:pdiv} together imply the following.

\begin{cor}\label{prop:dd}  Under Hypotheses \ref{hyp},
let $\G$ be a $p$-divisible group over $\OO_F$. Giving descent
data on $\G$ relative to $F'$ is equivalent to giving additive
bijections $\ghat : \M_{\pi}(\G) \rightarrow \M_{\pi}(\G)$ for all
$g \in \Gal(F/F')$ such that
\begin{itemize}
\item $\ghat(sx) = \ghat(s)\ghat(x)$ for $s \in S$, $x \in
\M_{\pi}(\G)$, and $g \in \Gal(F/F')$,

\item $\ghat(\Fil^1\M(\G))\subset\Fil^1\M(\G)$ and $\p \circ \ghat
= \ghat \circ \p$ for all $g \in \Gal(F/F')$, and

\item $\ghat_1 \circ \ghat_2 = \widehat{g_1 \circ g_2}$ for
$g_1,g_2 \in \Gal(F/F')$.
\end{itemize}
\end{cor}

This motivates the following definition.

\begin{defn} Assume Hypotheses \ref{hyp}.  If $\Mod$ is any one of
the categories of Section \ref{modules} (${}' \Mod^{k-1}$, etc.),
then the category $\Mod_{\mathrm{dd}}$ consists of
objects $\M$ of $\Mod$ together with additive bijections
$\ghat : \M \rightarrow \M$ for each $g$ in $\Gal(F/F')$
and satisfying the following compatibilities:
\begin{itemize}

\item $\ghat(sx) = \ghat(s)\ghat(x)$ for $s \in S$, $x \in \M$, $g
\in \Gal(F/F')$,

\item $\ghat(\Fil^{k-1} \M) \subset \Fil^{k-1}(\M)$ and
$\phi_{k-1} \circ \ghat = \ghat \circ \phi_{k-1}$ for each $g \in
\Gal(F/F')$,

\item $\ghat_1 \circ \ghat_2 = \widehat{g_1 \circ g_2}$ for $g_1,
g_2 \in \Gal(F/F')$,

\item $N \circ \ghat = \ghat \circ N$ if the category $\Mod$ is
equipped with an $N$.
\end{itemize}

Morphisms in $\Mod_{\mathrm{dd}}$ are those of $\Mod$ which commute
with $\ghat$ for all $g \in \Gal(F/F')$.
\end{defn}

So we may rephrase Corollary \ref{prop:dd} as follows: under
Hypotheses \ref{hyp}, the category of p-divisible groups over
$\OO_F$ with descent data relative to $F'$ is equivalent to the
category $\Mod^1_{0,\mathrm{dd}}$.

\begin{defn}  Under Hypotheses \ref{hyp}, we
refer to $\Mod^{1}_{0,\mathrm{dd}}$ as the category of
\textit{strongly divisible modules with tame descent data} (of
weight $2$).
\end{defn}

\begin{remark} \label{remark:DP} Retain Hypotheses \ref{hyp}, and
let $\G$ be a $p$-divisible group over $\OO_F$ with descent data
relative to $F'$.  Then $\G(1)$ is a finite flat group scheme
killed by $p$ with descent data relative to $F'$, and the filtered
$\p$-module $\M_{\pi}(\G(1))$ obtains descent data in the sense of
\cite[Th. 5.6.1]{BCDT}.  By construction, this descent data
is exactly the collection of maps induced on $$\M(\G)/p\M(\G)
\otimes_{S_1} k[u]/u^{ep}$$ by the descent data on $\M(\G)$.

\end{remark}

\subsection{Galois lattices and strongly divisible modules with descent data}

In the two preceding sections, we have seen how to pass between
Galois lattices inside potentially crystalline Galois
representations of $G_{F'}$ with Hodge-Tate weights in $\{0,1\}$
and $p$-divisible groups over $\OO_F$ with descent data relative
to $F'$, and between these and strongly divisible modules with
descent data. We now describe how to pass directly to Galois
lattices from strongly divisible modules with descent data.

 We may extend the natural
action of $G_F$ on $\Acrish$ to an action of $G_{F'}$. In fact,
more generally if $A$ is a syntomic $\OO_F$-algebra with a
$\Gal(F/F')$-semilinear action of $G_{F'}$, then $G_{F'}$ acts on
$\Ocris(A)$ as follows: if $g \in G_{F'}$, then $x \in \Ocris(A)$
maps to $(g\cdot x) \otimes 1$ under the composition
$$ \Ocris(A) \rightarrow \Ocris(A \otimes_{g^{-1}} \OO_F)
\xrightarrow{\sim} \Ocris(A) \otimes_{g^{-1}} S_n ,$$ where the
first map is induced by the $\OO_F$-algebra map $A \rightarrow A
\otimes_{g^{-1}} \OO_F$, $a \mapsto g(a) \otimes 1$.  Under
Hypotheses \ref{hyp}, it is not difficult to see (by checking on
presheaves) that $g\cdot u = h_g u$ for the element $u \in \Acrish$.
Therefore $g$ preserves the filtration and commutes with the $\phi$
induced on $\Acrish$ from the ring $\Asth$ as defined in \cite[Sec. 2.2.2]{BreuilStrong}, and furthermore the $\Acris$-linear map $f_{\pi}: \Acrish
\rightarrow B_{\mathrm{dR}}^{+}$ sending $u$ to $\pi$ is actually a
$G_{F'}$-morphism.  Thus we may regard $\Ocris(\OO_{\Qpbar})$ and
$\Acrish$ as objects of $\uMod^1_{0,\mathrm{dd}}$ and ${}' \uMod^1$
respectively.

Let $\M$ be a strongly divisible module with tame descent data,
and let $\G = \cup \G(n)$ be the $p$-divisible group over $\OO_F$ with
descent data relative to $F'$ such that $\M \cong \M_{\pi}(\G)$.
Forgetting the descent data momentarily, by \cite[Th. 4.2.2.9]{Br} and
the construction in
\cite[Sec. 4.2.1]{Br} we know that
$$ \G(n)(\OO_{\Qpbar}) = \Hom_{{}' \uMod^1_0}(\M/p^n\M,
\Ocris(\OO_{\Qpbar}))$$ is an isomorphism of $G_F$-modules.  The
crucial point is that it is actually an isomorphism of
$G_{F'}$-modules, where $\tilde{g} \in G_{F'}$ acts on the
right-hand side via $f \mapsto \tilde{g} \cdot (f \circ
\ghat^{-1})$. (To simplify notation, in this section we use
$\tilde{g}$ to denote an element of $G_{F'}$, and $g$ to denote
its restriction $\tilde{g} \,|_F$.)   More generally, if
$\G(n)=\Spec(R_n)$ and $A$ is a syntomic $\OO_F$-algebra with a
$\Gal(F/F')$-semilinear action of $G_{F'}$, we show that the
canonical bijection
$$ \G(n)(A) = \Hom_{\OO_F}(R_n,A) \xrightarrow{\sim}
\Hom_{(\p,\Fil^1)}(\M(\G(n)),\Ocris(A))$$ is a $G_{F'}$-module
isomorphism. Here, the subscript $(\p,\Fil^1)$ denotes morphisms
that commute with $\p$ and preserve $\Fil^1$, with
$\Fil^1(\Ocris(A))=\Jcris(A)$.  Indeed, consider the diagram
$${\footnotesize
\xymatrix{ \Hom_{\OO_F}(R_n,A) \ar[d] \ar[r] &
\Hom_{(\p,\Fil^1)}(\M(\G(n)),
\Ocris(A)) \ar[d] \\
\Hom_{\OO_F}(R_n \otimes_{g^{-1}} \OO_{F},A) \ar[d] \ar[r] &
\Hom_{(\p,\Fil^1)}(\M({}^{g^{-1}}\G(n)),
\Ocris(A)) \ar[d] \\
\Hom_{\OO_F}(R_n \otimes_{g^{-1}} \OO_{F},A \otimes_{g^{-1}}
\OO_{F}) \ar[dd] \ar[r] & \Hom_{(\p,\Fil^1)}(\M({}^{g^{-1}}\G(n)),
\Ocris(A \otimes_{g^{-1}} \OO_{F})) \ar[d] \\
& \Hom_{(\p,\Fil^1)}(\M(\G(n)) \otimes_{\ghat^{-1}} S_n,
\Ocris(A) \otimes_{\ghat^{-1}} S_n) \ar[d] \\
\Hom_{\OO_F}(R_n,A) \ar[r] & \Hom_{(\p,\Fil^1)}(\M(\G(n)),
\Ocris(A)) } }$$ in which
\begin{itemize}

\item the top square is functorial, induced by $[g^{-1}] : \G(n)
\rightarrow {}^{g^{-1}}\G(n)$, and hence commutes;

\item the middle square is functorial, induced by $A \rightarrow A
\otimes_{g^{-1}} \OO_F$, $a \mapsto \tilde{g}(a) \otimes 1$, and hence
commutes;

\item the left-hand vertical map in the bottom square is
``untwisting'', that is, takes a map sending $r\otimes 1 \mapsto
a\otimes 1$ to a map sending $r \mapsto a$; the first right-hand
vertical map in the bottom square is induced by the isomorphism
$\Ocris \otimes_{\ghat^{-1}} S_n \xrightarrow{\sim} \Ocrisg$; and
the second right-hand vertical map is again untwisting.
\end{itemize}

The actions of $\tilde{g}$ on $\Hom_{\OO_F}(R_n,A)$ and
$\Hom_{(\p,\Fil^1)}(\M(\G(n)),\Ocris(A))$ are the composites of
the left-hand and right-hand vertical maps in the above diagram,
respectively; hence it suffices to verify that the bottom square
commutes. Indeed, if ${\mathscr A},{\mathscr B} \in ({\rm
Ab}/\OO_{F})$ are any two abelian sheaves, then one checks locally
on sections that the composition
\begin{eqnarray*} \Hom_{({\rm
Ab}/\OO_{F})}({\mathscr A},{\mathscr B}) & \rightarrow &
\Hom_{(\p,\Fil^1)}(\Ocrisgi({\mathscr B}),\Ocrisgi({\mathscr A}))
\\
 & \xrightarrow{\sim} & \Hom_{(\p,\Fil^1)}(\Ocris({\mathscr
B})\otimes_{\ghat^{-1}} S_n,\Ocris({\mathscr
A})\otimes_{\ghat^{-1}} S_n) \\
 & \xrightarrow{\sim} &
\Hom_{(\p,\Fil^1)}(\Ocris({\mathscr B}),\Ocris({\mathscr A}))
\end{eqnarray*}
is just the natural map $\Hom_{({\rm Ab}/\OO_{F})}({\mathscr
A},{\mathscr B}) \rightarrow \Hom_{(\p,\Fil^1)}(\Ocris({\mathscr
B}),\Ocris({\mathscr A}))$.  This yields the conclusion.  Passing to
the inverse limit over $n$, we obtain the following.
\begin{thm} \label{thm:dd}
Assume Hypotheses \ref{hyp}. Suppose that $\G$ is a $p$-divisible
group over $\OO_F$ with descent data relative to $F'$, and let
$\M$ be the corresponding strongly divisible module with descent
data.  Then there is an isomorphism of $G_{F'}$-lattices
$$ T_{p}(\G) = \invlim \G(n)(\OO_{\Qpbar}) \cong \Hom_{{}' \uMod^1_0}(\M,\Acrish) .$$
\end{thm}

\subsection{From $\Acrish$ to $\Asth$} \label{sec35}

Assume Hypotheses \ref{hyp}, and let $\M$ be a strongly divisible
module with tame descent data. Recall that $\phi : \M \rightarrow
\M$ is defined as $\phi(x) = \frac{1}{c} \p(E(u)x)$.

From the equivalence of categories between $\Mod^1_0$ and
$\Mod^1_{\mathrm{cris}}$, our strongly divisible module $\M$ obtains
a monodromy operator $N$.  Because our fixed uniformizer satisfies
$\pi^{e(F/F')} \in F'$, it follows that $\ghat(E(u))=E(u)$. Since
$\ghat$ commutes with $\p$, it also commutes with $\phi$, and
so with $N$ as well: indeed, for the latter, note that $\ghat^{-1} N
\ghat$ satisfies the three properties of
\cite[Prop 5.1.3(1)]{Br}, and then invoke the uniqueness of $N$.

Moreover, any $S$-linear map from $\M$ to $\Acrish$, or to another
strongly divisible module, that preserves $\Fil^1$ and commutes with
$\p$, automatically commutes with $N$.  (If $f: \M \rightarrow
\M'$ is such a map, one sees iteratively that the $S$-linear map
$\Delta = f \circ N - N \circ f$ has $\Delta(\M) \subset
\phi^{m}(I)\M'$ for all $m$.)  Thus the equivalence of categories
between $\Mod^1_0$ and $\Mod^1_{\mathrm{cris}}$ extends to an
equivalence of $\Mod^1_{0,\mathrm{dd}}$ and
$\Mod^1_{\mathrm{cris},\mathrm{dd}}$, and
$$\Hom_{{}' \uMod^1_0}(\M,\Acrish) =
\Hom_{{}' \uMod^1}(\M,\Acrish).$$  Henceforth, when we refer to a
strongly divisible module with tame descent data (of weight $2$), we
are typically referring to an object of
$\Mod^1_{\mathrm{cris},\mathrm{dd}}$ (i.e., the corresponding object
of $\Mod^1_{0,\mathrm{dd}}$ endowed with its canonical $N$).

It is not difficult to see that our action of $G_{F'}$ on $\Acrish$
extends uniquely to $\Asth$.  We then have the following.

\begin{prop} \label{prop:ast} The embedding $\Acrish \rightarrow \Asth$ induces an
isomorphism of $G_{F'}$-lattices
$$\Hom_{{}' \uMod^1}(\M,\Acrish) = \Hom_{{}' \uMod^1}(\M,\Asth).$$
\end{prop}

\begin{proof}  The induced map is evidently injective, so we need
to prove surjectivity.  If $\gamma \in \Hom_{{}'
\uMod^1}(\M,\Asth)$ it is not difficult to see that $\gamma(\M)
\subset \Acrish[\frac{1}{p}]$, for example using \cite[Prop. 5.1.3]{Br} and the fact that $\{x \in \Asth \ | \ Nx = 0\} =
\Acris$. Therefore given $\gamma \in \Hom_{{}'
\uMod^1}(\M,\Acrish)$ such that $\gamma(\M) \subset p\Asth$, we
need only show that $\gamma(\M)\subset p\Acrish$.

We remark first that $p$ divides $\p\left(\frac{(u -
\upi)^i}{i!}\right)$ in $\Acrish$ if $i \ge 2$.  Indeed,
$$ \phi\left(\frac{(u - \upi)^i}{i!}\right) = \frac{1}{i!} \left(p! \sum_{j=1}^{p}
\frac{(u-\upi)^j}{j!} \frac{\upi^{p-j}}{(p-j)!}\right)^i $$ and
$\frac{p!^i}{i!p}$ is divisible by $p$ for $i \ge 2$ since $p$ is
odd. Now suppose $x \in \Fil^1\M$, and write
$$\gamma(x) = \sum_{i \ge 0} a_i \frac{(u - \upi)^i}{i!} \in \Fil^1
\Acrish ,$$ with each $a_i$ in $\Acris$.   (In particular, $a_0
\in \Acris \cap \Fil^1 \Acrish$.)  Then
$$ \p(\gamma(x))
\in \p(a_0) + \phi(a_1)\left(\frac{u^p - \upi^p}{p}\right) +
p\Acrish .$$ Since $\p(\gamma(x))=\gamma(\p(x)) \in \gamma(\M)
\subset p\Asth$ it follows that $p$ divides $\p(a_0)$ in $\Acris$
(use that the map $\Asth \rightarrow \Acris, X \mapsto 0$ sends
$u$ to $\upi$), and therefore
$$\p(\gamma(x)) \in \phi(a_1)\frac{u^p - \upi^p}{p} + p\Acrish \subset
A$$ where $A \subset \Acrish$ is the subset of elements of the
form
$$ pb + \sum_{i=1}^{p} b_i \frac{\upi^{p-i}}{(p-i)!}
\frac{(u-\upi)^i}{i!}$$ with $b \in \Acrish$ and each $b_i \in
\Acris$.
 Since $\M$ is generated over $S$ by $\phi_1(\Fil^1\M)$, it
follows that $\gamma(\M)$ is contained in the subset of $\Acrish$
generated over $S$ by $A$, and we deduce that every element of
$\gamma(\M)$ is of the form $\sum_{i \ge 0} b_i
\frac{(u-\upi)^i}{i!}$ with $b_1 = \upi^{p-1} c_1 + p c_2$ with
$c_1,c_2 \in \Acrish$. In particular this applies to $a_1$, so
that $\phi(a_1)$ is divisible by $p$ and $\gamma(\M) \subset
p\Acrish$.
\end{proof}

\begin{remark} Note that on the level of rings, it is not true
that $\Acrish \cap p\Asth = p\Acrish$.
\end{remark}

We note for future reference that this $G_{F'}$-lattice may, by
the proof of \cite[Prop. 2.3.2.4]{BreuilStrong}, be
written as
$$ \invlim_n  \Hom_{{}' \uMod^1}(\M/p^n\M,\Asthi) .$$
If $\M$ is a strongly divisible module with descent data, we
define $G_{F'}$-modules
$$ \VstFtwo(\M/p^n\M) = \Hom_{{}' \uMod^1}(\M/p^n\M,\Asthi)$$
and
$$ \TstFtwo(\M) = \invlim_n \VstFtwo(\M/p^n\M)\widehat{\ }(1) ,$$
where $\widehat{\ }$ denotes the ($\Qp/\Zp$)-dual and where the
$(1)$ is a twist by the cyclotomic character.

If $D$ is a filtered $(\varphi,N,F/F',\Zp)$-module, we say $\M$ is
contained in $S[1/p] \otimes_{F_0} D$ if $\M \otimes_{W(k)} F_0
\cong S[1/p] \otimes_{F_0} D$, the isomorphism respecting $N$,
$\phi$, the filtration, and descent data (which acts on $S[1/p]
\otimes_{F_0} D$ in the obvious manner).  We recall (see, e.g., 
\cite[Sec. 3.2.3]{BreuilMezard}) that $$\Fil^1 (S
\otimes_{F_0} D) = \left\{ \sum s_i(u) \otimes d_i \ | \ \sum
s_i(\pi) d_i \in \Fil^1 D_F \right\}.$$ (Recall that $D_F$
denotes $F \otimes_{F_0} D$.)

\begin{lemma} \label{lemma:cont} If $\M$ is a strongly divisible
module with tame descent data, then there exists a filtered
$(\varphi,N,F/F',\Qp)$-module such that
\begin{itemize}
\item $\M$ is contained in $S[1/p] \otimes_{F_0} D$,

\item $N=0$ on $D$, and

\item $\Fil^i D_F = D_F$ if $i \le 0$, and $\Fil^i D_F = 0$ if $i \ge
2$.
\end{itemize}
\end{lemma}

\begin{proof}  Forgetting descent data momentarily, by
\cite[Prop. 5.1.3(2)]{Br} we obtain a filtered
$(\varphi,N,\Qp)$-module $D$ satisfying the above conditions on
$N$ and $\Fil^i D_F$ and such that $\M \otimes_{W(k)} F_0 \cong
S[1/p] \otimes_{F_0} D$, the isomorphism respecting $N$, $\phi$,
and the filtration. However, since this isomorphism identifies $D$
with $\ker(N)$ on $\M \otimes_{W(k)} F_0$, and since each $\ghat$
commutes with $N$, it follows that each $\ghat$ acts on $D$.  Thus
$D$ is actually a filtered $(\varphi,N,F/F',\Zp)$-module.
\end{proof}

Finally, we have the following.

\begin{thm}  Retain the hypotheses of Corollary \ref{prop:dd}.
Suppose that $\rho : G_{F'} \rightarrow \GL(V)$ becomes
crystalline over $F$ and has Hodge-Tate weights in $\{0,1\}$. The
functor $\TstFtwo$ is an equivalence between the category of
strongly divisible modules with tame descent data contained in
$S[1/p] \otimes_{F_0} \DstFtwo(V)$ and the category of
$G_{F'}$-lattices in $\rho$.
\end{thm}

\begin{proof} By Lemma \ref{lemma:cont}, Theorem \ref{thm:dd}, Corollary \ref{prop:dd}, and
Propositions \ref{prop:lattice}  and
\ref{prop:ast}, it suffices to prove that if $\M$ is contained in
$S[1/p] \otimes_{F_0} \DstFtwo(V)$, then $\TstFtwo(\M)$ is a
$G_{F'}$-lattice in $\rho$. This follows by the same proof as
\cite[Lem. 3.2.3.1]{BreuilMezard}: one simply notes that each
map in that proof is now a $G_{F'}$-map, and not just a $G_F$-map.
\end{proof}

\section{Coefficients}
\label{coef}

Throughout this section, we assume Hypotheses \ref{hyp}.

We now wish to add coefficients to our theory of strongly divisible
modules. Specifically, let $E$ be a finite extension of $\Qp$, let
$\OO_E$ be its ring of integers, and let $R$ be a complete local
noetherian flat $\OO_E$-algebra with maximal ideal $\mm_R$ and
residue field a finite extension of the residue field $\kk_E$ of
$\OO_E$. Let $\kk_F$ be the residue field of $F$. We construct a category
$R-\Mod^{k-1}_{\mathrm{cris},\mathrm{dd}}$, the category of strongly
divisible $R$-modules with tame descent data, having roughly the
following properties:
\begin{itemize}
\item there is a functor $\Tst$ from $R-\Mod^{k-1}_{\mathrm{cris},\mathrm{dd}}$ to
$R$-representations of $G_{F'}$ for each $R$, compatible with base
change $R \rightarrow R'$, and

\item when $k=2$ and $R=\OO_E$, the functor $T_{\mathrm{st},2}$ is an equivalence of
categories between $R-\Mod^{1}_{\mathrm{cris},\mathrm{dd}}$ and the
category of $\OO_E$-lattices inside representations of $G_{F'}$ with
Hodge-Tate weights $\{0,1\}$ and becoming crystalline over~$F$,
coinciding with $\TstFtwo$ when $E=\Qp$.
\end{itemize}
Our exposition follows that of
\cite[Sec. 3.2]{BreuilMezard} as closely as possible (verbatim in many
places), but some changes are forced by the lack of any
restrictions on $e = e(F)$.

Set $S_{F,R}$ to be the ring
$$ \left\{ \sum_{j=0}^{\infty} r_j \frac{u^j}{\lfloor j/e \rfloor
!}, \ \text{where} \ r_j \in W(\kk_F) \otimes_{\Zp} R, \ r_j
\rightarrow 0 \ \mm_R\text{-adically} \ \text{as} \ j \rightarrow
\infty\right\}.$$ Extend the definitions of $\Fil, \phi, \phi_k,
N, \ghat$ to $S_{F,R}$ in the evident ($R$-linear) manner; for
example, $\Fil^{k-1} S_{F,R}$ is the $\mm_R$-adic completion of
the ideal generated by the $E(u)^j / j!$ for $j \ge k-1$.

We remark that if $I$ is any ideal of $R$, then $$I S_{F,R} \cap
\Fil^{k-1} S_{F,R} = I \Fil^{k-1} S_{F,R}.$$  Indeed, every element
of $S_{F,R}$ may be written uniquely in the form $$\sum_{j \ge 0}
r_j(u) (E(u)^j / j!)$$ with $r_j(u)$ a polynomial of degree
less than $e(F)$ over $W(\kk_F) \otimes R$. For an element of $I
S_{F,R} \cap \Fil^{k-1} S_{F,R}$, it follows (by uniqueness) that
$r_j(u) = 0$ for $j < k-1$ and the coefficients of $r_j(u)$ lie in
$W(\kk_F) \otimes I$ for $j \ge k-1$.  Since $R$ is noetherian, such
an element is actually in $I \Fil^{k-1} S_{F,R}$.

Note that if $R$ is the ring of integers in a local field, then we
actually have $S_{F,R} = R \otimes_{\Zp} S_F$.  We often
abbreviate $S_{F,R}$ by $S_R$.

\begin{defn} A strongly divisible $R$-module with tame descent
data is a finitely generated free $S_R$-module $\M$, together with
a sub-$S_R$-module $\Fil^{k-1}\M$, maps $\phi,N : \M
\rightarrow \M$, and additive bijections $\ghat : \M \rightarrow
\M$ for each $g \in \Gal(F/F')$, satisfying the following
conditions:
\begin{enumerate}
\item $\Fil^{k-1}\M$ contains $(\Fil^{k-1}S_R)\M$,

\item $\Fil^{k-1}\M \cap I\M = I \Fil^{k-1}\M$ for all ideals $I$
in $R$,

\item $\phi(sx) = \phi(s)\phi(x)$ for $s \in S_R$ and $x \in \M$,

\item $\phi(\Fil^{k-1}\M)$ is contained in $p^{k-1}\M$ and
generates it over $S_R$,

\item $N(sx) = N(s)x + sN(x)$ for $s \in S_R$ and $x \in \M$,

\item $N\phi = p\phi N$,

\item $E(u)N(\Fil^{k-1}\M) \subset \Fil^{k-1}\M$,

\item $N(\M) \subset J\M$ where $J$ is the ideal
$\sum_{j \ge 1} \frac{u^j}{\lfloor j/e \rfloor!} S_R$ in $S_R$,

\item $\ghat(sx) = \ghat(s)\ghat(x)$ for all $s \in S_R$, $x \in
\M$, $g \in \Gal(F/F')$,

\item $\ghat_1 \circ \ghat_2 = \widehat{g_1 \circ g_2}$ for all
$g_1,g_2 \in \Gal(F/F')$,

\item $\ghat(\Fil^{k-1}\M) \subset \Fil^{k-1}\M$ for all $g \in
\Gal(F/F')$, and

\item $\phi$ and $N$ commute with $\ghat$ for all $g \in
\Gal(F/F')$.

\end{enumerate}
The category $\Rmod$ consists of strongly divisible $R$-modules with
tame descent data, along with $S_R$-linear morphisms that preserve
$\Fil^{k-1}$ and commute with $\phi$, $N$, and descent data.
\end{defn}

\begin{example}  If $R = \OO_E = \Zp$ and $k=2$, then $\Rmod$ is the
category $\Mod^1_{\mathrm{cris},\mathrm{dd}}$.
\end{example}

\begin{example} If $F=F'=\Qp$, then our strongly divisible
$R$-modules are precisely those strongly divisible $R$-modules of
\cite[D\'ef. 3.2.1.1]{BreuilMezard} which satisfy the extra
condition $N(\M) \subset J\M$; that is, our strongly divisible
$R$-modules are all ``crystalline'', whereas those of
\cite{BreuilMezard} may be ``semistable''.
\end{example}

\begin{defn} Let $I \subset R$ be an ideal containing $\mm_R^n$
for $n$ sufficiently large.  An \textit{object of}
$\uMod^{k-1}_{\mathrm{dd}}$ \textit{with an action of} $R/I$ is an
object $\NN$ of $\uMod^{k-1}_{\mathrm{dd}}$ together with an algebra
map $R/I \rightarrow \textrm{End}_{\uMod^{k-1}_{\mathrm{dd}}}(\NN)$.
Such an $\NN$ is an $S_R/IS_R$-module.
\end{defn}

\begin{example} If $\M$ is a strongly divisible $R$-module and $I$ is an arbitrary
ideal of $R$, let $\Fil^{k-1}(\M/I\M)$ be the image of
$\Fil^{k-1}\M/I\Fil^{k-1}\M \hookrightarrow \M/I\M$.  If $R/I$ is
flat, then $\M/I\M$ together with $\Fil^{k-1}(\M/I\M)$ and the
reductions modulo $I$ of $\phi$, $N$, $\ghat$, is a strongly
divisible $R/I$-module.  If $R/I$ is Artinian, then $\M/I\M$
together with $\Fil^{k-1}(\M/I\M)$ and the reductions modulo $I$ of
$\phi$, $N$, $\ghat$, is an object of $\uMod^{k-1}_{\mathrm{dd}}$
with an action of $R/I$.
\end{example}

We have the following weaker version of
\cite[Lem 3.2.1.3]{BreuilMezard}, adapted for the fact that given a morphism
$f: \NN^r \rightarrow \NN$ in $\uMod^{k-1}$, the identity
$f(\Fil^{k-1}\NN^r) = \Fil^{k-1}\NN \cap f(\NN^r)$ may not hold
when the ramification index $e$ is large.

\begin{lemma} \label{lemma46} Suppose $I$ is an ideal of $R$ containing $\mm_R^n$
for $n$ sufficiently large, $R'$ is a local Artinian $\OO_E$-algebra
with residue field a finite extension of $\kk_E$, and $R/I \rightarrow R'$ is a
local $\OO_E$-algebra morphism.  Suppose that $\NN$ is an object of
$\uMod^{k-1}_{\mathrm{dd}}$ with $R/I$-action, and that either
\begin{enumerate}
\item $\NN = \M/I\M$ for some strongly divisible $R$-module $\M$ with
tame descent data, or

\item $R'$ is isomorphic to $(R/(p^r,I))^n$ as an $R/I$-module.
\end{enumerate}
Then $\NN \otimes_{R/I} R'$ is an object of
$\uMod^{k-1}_{\mathrm{dd}}$ with $R'$-action, and $\NN \rightarrow
\NN \otimes_{R/I} R'$ is a morphism in $\uMod^{k-1}_{\mathrm{dd}}$.
\end{lemma}

\begin{proof}  The result is clear if $R'$ is a free $R/I$-module,
so we may assume that $R \rightarrow R'$ is surjective.  In case
(2), the result follows as in \cite{BreuilMezard} from the fact
that $p^r \NN \cap \Fil^{k-1} \NN$ does equal $p^r \Fil^{k-1}\NN$.
In case (1), suppose that $R' = R/I'$ with $I' \supset I$.  Then
$N \otimes_{R/I} R' = \M/I'\M$.
\end{proof}

\begin{cor} \label{cor47} Suppose that $R \rightarrow R'$ is a finite local map of
complete local noetherian flat $\OO_E$-algebras, and suppose that either
\begin{enumerate}
\item this map is surjective, or
\item every non-zero ideal $I$ of $R'$ has $I^m = (p^r)$ for some positive integers
$m$ and $r$.
\end{enumerate}
If $\M$ is a strongly divisible $R$-module with descent data, then
$\M \otimes_R R'$, equipped with $\phi \otimes 1$, $N \otimes 1$,
$\ghat \otimes 1$, and the image of $\Fil^{k-1} \M \otimes R'$, is
a strongly divisible $R'$-module with descent data.
\end{cor}

\begin{proof} In the first case, use the proof of Lemma
\ref{lemma46}
and the fact that every ideal of $R'$ is the image of an ideal of
$R$, and pass to the appropriate inverse limit.  Write $\M' = \M
\otimes_R R'$.  In the second case, we use the fact that we know
$p^r \M' \cap \Fil^{k-1}\M'$ does equal $p^r \Fil^{k-1} \M'$. 
If $I\Fil^{k-1}\M' \subsetneq I\M' \cap \Fil^{k-1}\M'$, then
inductively we would also have $I^m \Fil^{k-1}\M' \subsetneq I^m
\M' \cap \Fil^{k-1}\M'$, a contradiction.
\end{proof}

\begin{remark} \label{finext} In particular, the conclusions of Corollary \ref{cor47}
hold for any local $\OO_E$-algebra map of the form $R
\rightarrow \OO_{E'}$ with $E'$ a finite extension of $E$.
\end{remark}

If $\NN$ is an object of $\uMod^{k-1}_{\mathrm{dd}}$, we may define,
as in Section \ref{sec35}, $G_{F'}$-modules
$$ V_{\mathrm{st,k}}^{F'}(\NN) = \Hom_{{}' \uMod^{k-1}}(\NN,\Asthi)$$
and
$$ T_{\mathrm{st},k}^{F'}(\NN) = V_{\mathrm{st,k}}^{F'}(\NN)\widehat{\ }(k-1) .$$
When $\NN$ has an action of $R/I$, so does
$T_{\mathrm{st},k}^{F'}(\NN)$, as in \cite[Sec. 3.2.2]{BreuilMezard}.

\begin{lemma}  \label{red-modp} Let $I \subset I'$ be ideals of $R$ containing $\mm_R^n$ for
$n$ sufficiently large.  Let $\M$ be a strongly divisible $R$-module
of rank $d$ with tame descent data.  Then we have the following.
\begin{enumerate}
\item The map $T_{\mathrm{st},k}^{F'}(\M/I\M) \rightarrow
T_{\mathrm{st},k}^{F'}(\M/I'\M)$ is surjective.

\item The $R/I$-module $T_{\mathrm{st},k}^{F'}(\M/I\M)$ is free of rank
$d$.
\end{enumerate}
\end{lemma}

\begin{proof}  The proof is exactly the same as that of
\cite[Lem. 3.2.2.1]{BreuilMezard}, replacing
\cite[Prop. 3.2.3.1, Cor. 3.2.3.2]{BreuilENS} with \cite[Lems. 2.3.1.1, 2.3.1.2, 2.3.1.3]{BreuilStrong}.
\end{proof}

\begin{lemma}  \label{bc-modp} Let $I$ be an ideal of $R$ containing $\mm_R^n$ for
$n$ sufficiently large, let $R'$ be an Artinian local $\OO_E$-algebra
with residue field a finite extension of $\kk_E$, and let $R/I
\rightarrow R'$ be a local morphism of $\OO_E$-algebras.  If $\M$ is
a strongly divisible $R$-module with tame descent data, then
$$ T^{F'}_{\mathrm{st},k}(\M/I\M) \otimes_{R/I} R' \cong T^{F'}_{\mathrm{st},k}(\M/I\M
\otimes_{R/I} R').$$
\end{lemma}

\begin{proof} The proof is exactly the same as that of
\cite[Lem. 3.2.2.2]{BreuilMezard}, substituting Lemma \ref{red-modp}
for \cite[Lem. 3.2.2.1]{BreuilMezard}.
\end{proof}

\begin{defn} If $\M$ is a strongly divisible $R$-module, set
$$T^{F'}_{\mathrm{st},k}(\M) = \underset{n}{\invlim} \ T^{F'}_{\mathrm{st},k}(\M/\mm_R^n\M) .$$
This is naturally an $R[G_{F'}]$-module.
\end{defn}

Finally, using Lemmas \ref{red-modp} and \ref{bc-modp} and passing
to the limit, we have the following.

\begin{cor} \label{cor:reduction} Let $\M$ be a strongly divisible $R$-module with
descent data.
\begin{enumerate}
\item $T^{F'}_{\mathrm{st},k}(\M)$ is a free $R$-module of rank $d$
with a continuous action of $G_{F'}$, and
$$T^{F'}_{\mathrm{st},k}(\M)/\mm_R^n \overset{\sim}{\rightarrow}
T^{F'}_{\mathrm{st},k}(\M/\mm_R^n).$$
\item If $R'$ is another complete local noetherian flat
$\OO_E$-algebra with residue field a finite extension of $\kk_E$,
and if $R \rightarrow R'$ is a local map such that $\M \otimes_R R'$ is a strongly
divisible $R'$-module with descent data, then
$$  T^{F'}_{\mathrm{st},k}(\M)
\otimes_R R' \overset{\sim}{\rightarrow}  T^{F'}_{\mathrm{st},k}(\M
\otimes_R R') .$$
\end{enumerate}
\end{cor}

Suppose $k=2$.  It remains to verify in this case that when
$R=\OO_E$, the category of strongly divisible $R$-modules with
descent data corresponds to the category of lattices in
potentially Barsotti-Tate $E$-representations of $G_{F'}$.  Let
$\M$ be a strongly divisible $\OO_E$-module.  Regarding $\M$ as a
strongly divisible $\Zp$-module, from Lemma \ref{lemma:cont} we
obtain a filtered $(\varphi,N,F/F',\Qp)$-module $D$ such that
$$ \M \otimes_{W(k)} F_0 \cong S_{\Zp}[1/p] \otimes_{F_0} D ,$$
and such that $D = \{ x \in \M \otimes_{W(k)} F_0 \ | \ Nx = 0\}$.
Since $Nx = 0$ implies $N(\alpha x) = 0$ for any $\alpha \in \OO_E$, it
follows that the action of $\OO_E$ on $\M$ preserves $D$; in this
manner $D$ is a filtered $(\varphi,N,F/F',E)$-module, and
$$ \M \otimes_{W(k)} F_0 \cong S_{\OO_E}[1/p] \otimes_{F_0 \otimes
E} D .$$  Suppose that the filtered $(\varphi,N,F/F',E)$-module
$D$ is $D_{\mathrm{st},2}^{F}(\rho)$ for the potentially
Barsotti-Tate representation $\rho : G_{F'} \rightarrow \GL_d(E)$
becoming Barsotti-Tate over $F$.  By the proof of
\cite[Lem. 3.2.3.1]{BreuilMezard} (and noting that each map is now a $G_{F'}$-map,
and not just a $G_F$-map) we conclude that
$T_{\mathrm{st},2}^{F'}(\M)$ is a $G_{F'}$-stable $\OO_E$-lattice in
$\rho$.

We now check the following.

\begin{prop} Each $G_{F'}$-stable $\OO_E$-lattice $T$ in $\rho$ is
isomorphic to $T_{\mathrm{st},2}^{F'}(\M)$ for some strongly
divisible $\OO_E$-module with descent data $\M$.
\end{prop}

\begin{proof} We know $T$ is $\Zp$-isomorphic to $T_{\mathrm{st},k}^{F'}(\M)$
for a strongly divisible $\Zp$-module with descent data $\M$.  We
know from Corollary \ref{prop:dd} that $T$ gives rise to $\M$ via
a $p$-divisible group $\Gamma$ with descent data; since $T$ is an
$\OO_E$-module, the $p$-divisible group $\Gamma$ has an action of
$\OO_E$ by Tate's full faithfulness theorem in \cite{tate}, and so we
obtain a map
$$\OO_E \rightarrow \End_{\uMod^1_{0,\mathrm{dd}}}(\M) .$$

We must check that this makes $\M$ into a strongly divisible
$\OO_E$-module.  To do this, we first note that $\M/(\Fil^1
{S_{\Zp}}) \M$ is a torsion-free $\OO_E$-module: indeed, if $m
\not\in (\Fil^1 {S_{\Zp}})\M$ but $am \in (\Fil^1 {S_{\Zp}}) \M$,
then $p^r m \in (\Fil^1 {S_{\Zp}}) \M$ for sufficiently large $r$.
Hence there would exist $m' \not\in (\Fil^1 {S_{\Zp}}) \M$ such
that $pm' \in (\Fil^1 {S_{\Zp}}) \M$, which is not the case.  The
proof now proceeds just as the proof of
\cite[Prop. 3.2.3.2]{BreuilMezard}.

Finally, recalling that the isomorphism
$T_{\mathrm{st},2}^{F'}(\M)[1/p] \cong \rho$ is compatible with
$\OO_E$-structures, we conclude that the strongly divisible
$\OO_E$-module $\M$ gives rise to the $\OO_E$-lattice $T$.
\end{proof}

\subsection*{Objects killed by $p$}

As before, let $\kk_F$ be the residue field of $F$, and let $e$ be the absolute
ramification index of $F$.  Let $\BrMod$ be the category of
\textit{Breuil modules}, that is, the category of triples $(\M',
\Fil^1\M', \phi_1)$ such that:
\begin{itemize}
\item $\M'$ is a finite rank free $\kk_F[u]/u^{ep}$-module,
\item $\Fil^1\M'$ is a submodule of $\M'$ containing $u^e \M'$, and
\item $\phi_1 : \Fil^1\M' \rightarrow \M'$ is an additive map such
that $\phi_1(hv) = h^p \phi_1(v)$ for any $h \in \kk_F[u]/u^{ep}$
and $v \in \Fil^1 \M'$, and $\phi_1(\Fil^1 \M')$ generates $\M'$.
\end{itemize}

If $\M$ is an object of $\uMod^1$ which is killed by $p$, then $\M
\otimes_{S_1} \kk_F[u]/u^{ep}$ is an object of $\BrMod$; and in fact
(see \cite[Prop. 2.1.2.2]{Br}) this induces an equivalence of
categories $T_0$ between the subcategory of $\uMod^1$ of objects
that are killed by $p$ and $\BrMod$, with quasi-inverse $T'_0$ given
by $\M' \mapsto \M' \otimes_{\kk_F[u]/u^{ep}} S_1$. Moreover (see
Rem. \ref{remark:DP}) this extends to an equivalence between the
subcategory of $\uMod^1_{\mathrm{dd}}$ of objects that are killed by
$p$, and the Breuil modules with descent data $\BrMod_{\mathrm{dd}}$
(described in \cite[Th. 5.6.1]{BCDT} and
\cite[Sec. 3.3]{SavittCompositio}).

Here we note that if $\M$ is an object of $\uMod^1_{\mathrm{dd}}$
with an action of $R/I$, then the corresponding Breuil module
$T_0(\M)$ also has an action of $R/I$.  We thus obtain an
equivalence of categories between the subcategory of
$\uMod^1_{\mathrm{dd}}$ of objects that are killed by $p$ and have
an action of $R/I$, and Breuil modules with descent data and an
action of $R/I$.  (See the proof of
\cite[Prop. 2.2.2.1]{BreuilENS} to verify that the isomorphisms $T'_0(T_0(\M))\cong
\M$ and $T_0(T'_0(\M')) \cong \M'$ are compatible with the actions
of $R/I$.)  Hence when we study objects of $\uMod^1_{\mathrm{dd}}$
with an action of $R/I$ which are killed by $p$, it suffices to
consider the corresponding Breuil modules. By abuse of notation, if
$\M'$ is a Breuil module with descent data, we write
$T_{\mathrm{st},2}^{F'}(\M')$ for
$T_{\mathrm{st},2}^{F'}(T'_0(\M'))$.

It is worth remarking that while it is certainly not the case
that every Breuil module with an action of $R/I$ is free as a
$(\kk_F \otimes R/I)[u]/u^{ep}$-module, this is true of
Breuil modules arising as the reductions of strongly divisible
modules.

We make the following observation.

\begin{lemma} \label{maximal} Suppose that $\M'$ is a Breuil module with descent
data satisfying $\Fil^1 \M' = u^{e}\M'$.  If $\M''$ is another
Breuil module with descent data such that
$T_{\mathrm{st},2}^{F'}(\M')=T_{\mathrm{st},2}^{F'}(\M'')$, then
there is a nontrivial map $\M'' \rightarrow \M'$.  In the
terminology of \cite[Def. 8.1]{SavittCompositio}, $\M'$ is
the maximal Breuil module of $\M''$.
\end{lemma}

\begin{proof}  By the compatibility between Breuil modules and
Dieudonn\'e modules (see \cite[Th. 5.1.3(3)]{BCDT}), we see
that the group scheme $\G'$ corresponding to $\M'$ under the
contravariant functor $\G_{\pi}$ (of \cite[Th. 5.1.3(1)]{BCDT}) is
\'etale. Let $(\G')^{+}$ be the maximal prolongation of the generic
fibre of $\G'$ (see \cite{Raynaud}); by the universal property of the
connected-\'etale sequence (see, e.g., \cite[3.7(I)]{Tate2}), we
find that $\G' = (\G')^{+}$.  If $\G'' = \G_{\pi}(\M'')$, we
conclude that there is a map $\G' \rightarrow \G''$ which induces an
isomorphism on generic fibres, and therefore also a map $\M''
\rightarrow \M'$.
\end{proof}

\begin{remark} \label{minimal-breuil}  Similarly, if $\Fil^1 \M' =
\M'$ we have the minimal Breuil module.
\end{remark}

\section{Strongly divisible modules for characters}
\label{sec:chars}

In this section, we compute the strongly divisible $\OO_E$-modules
corresponding to lattices in the characters of Examples
\ref{characters} and \ref{char2}, in the case $k=2$.  The purpose is
to list Breuil modules with descent data and an action of $\kk_E =
\OO_E/\mm_E$ to which $T^{F'}_{\mathrm{st},2}$ associates the
reduction mod $\mm_E$ of these characters.

These particularly simple strongly divisible $\OO_E$-modules are
given by the following propositions.

\begin{prop}  \label{strchar1} Let $F_1 = \Qp(\zeta_p)$, fix $\pi = (-p)^{1/(p-1)}$
as the choice of uniformizer in $\OO_{F_1}$,  and consider the
character $\epsilon \om^j \lambda_a$ of $G_{\Qp}$ as in Example
\ref{characters}. Then a strongly divisible $\OO_E$-module in
$S_{\OO_E} \otimes D_{\mathrm{st},2}^{F_1}(\epsilon \om^j
\lambda_a)$ is given by
$$ \M = S_{\OO_E} \cdot \e , \quad  \Fil^1 \M  = \Fil^1 S_{\OO_E} \cdot \e ,$$
$$ \phi(\e) = a^{-1} \e , \quad \ N \e = 0 ,$$
$$ \ghat(\e) = \om^j(g) \e \ \text{for} \ g \in \Gal(F_1/\Qp).$$
\end{prop}

\begin{proof} The proof is clear.  For example, since $\Fil^1 D_{\mathrm{st},2}^{F_1}(\epsilon
\om^j \lambda_a) = 0$, it follows that $\Fil^1 \M = \{ s(u) \e \ |
\ s(\pi) = 0 \}$; that is, $\Fil^1 \M = \Fil^1 S_{\OO_E} \cdot \e$.
\end{proof}

Similarly we have the following.

\begin{prop}  Let $F_2 = \Qpp(\varpi)$, as in Example \ref{char2}, and fix $\varpi$
as the choice of uniformizer in $\OO_{F_2}$.

\begin{enumerate}

\item A strongly divisible $\OO_E$-module in
$S_{\OO_E} \otimes D_{\mathrm{st},2}^{F_2}(\epsilon \om^j
\lambda_a)$ is given by
$$ \M = S_{\OO_E} \cdot \e , \quad \Fil^1 \M  = \Fil^1 S_{\OO_E} \cdot \e ,$$
$$ \phi(\e) = (1 \otimes a^{-1}) \e , \quad N \e = 0 ,$$
$$ \ghat(\e) = (1 \otimes \om^j(g)) \e \ \text{for} \ g \in \Gal(F_2/\Qp).$$

\item Suppose that $E$ is a finite extension of $\Qpp$.  A strongly divisible $\OO_E$-module in
$S_{\OO_E} \otimes D_{\mathrm{st},2}^{F_2}(\omt^m (\epsilon
\lambda_a)|_{G_{\Qpp}})$ is given by:
$$ \M = S_{\OO_E} \cdot \e , \quad \Fil^1 \M  = \Fil^1 S_{\OO_E} \cdot \e ,$$
$$ \phi(\e) = (1 \otimes a^{-1}) \e , \quad N \e = 0 ,$$
$$ \ghat(\e) = (1 \otimes \omt^m(g)) \e \ \text{for} \ g \in \Gal(F_2/\Qpp).$$
\end{enumerate}
\end{prop}

Next, we have the following.

\begin{prop} \label{prop:chars1} Let $F_1=\Qp(\zeta_p)$, fix $\pi = (-p)^{1/(p-1)}$
as the choice of uniformizer in $\OO_{F_1}$, set $e_1 = p-1$, and
let $\kk_E$ be the residue field of $E$. Let $\M'$ be the Breuil
module with descent data and action of $\kk_E$ given by
$$ \M' = (\kk_E[u]/u^{e_1 p}) \e , \quad  \Fil^1 \M' = u^{e_1}
\M'$$
$$ \phi_1(u^{e_1} \e) = \overline{a}^{-1} \e , \quad \ghat(\e) = \omega^j(g) \e \ \text{for} \ g \in \Gal(F_1/\Qp).$$
Here $\overline{a}$ is the reduction of $a$ modulo $\mm_E$. Then
$T_{\mathrm{st},2}^{\Qp}(\M') = \lambda_{\overline{a}} \cdot
\omega^{j+1}$.
\end{prop}

\begin{proof}  Note that $\kk_E[u]/u^{e_1 p} = \kk_E \otimes
\F_p[u]/u^{e_1 p}$.   The proposition follows directly from (1) of
Corollary \ref{cor:reduction}, once one checks that $\M'$ is the
Breuil module corresponding to the reduction modulo $\mm_E$ of the
strongly divisible $\OO_E$-module $\M$ in
Proposition~\ref{strchar1}. This is easy: for example, $(u^{e_1} +
p) \e \in \Fil^1\M$ implies $u^{e_1} \e \in \Fil^1\M'$; and the
equality $\phi_1 ((u^{e_1} + p)\e) = (\frac{u^{e_1 p}}{p} +
1)a^{-1} \e$ in $\M$ implies $\phi_1(u^{e_1} \e) =
\overline{a}^{-1}\e$ in $\M'$.
\end{proof}

We denote the above Breuil modules by
$\M_E(F_1/\Qp,e_1,\overline{a}^{-1},j)$. Similarly we have the
following.

\begin{prop} \label{prop:chars2} Let $F_2=\Qpp(\varpi)$, fix $\varpi$ as the choice of
uniformizer in $\OO_{F_2}$, set $e_2 = p^2 -1 $, suppose that $E$
contains $\Qpp$, and let $\kk_E$ be the residue field of $E$.

\begin{enumerate}

\item Let $\M'$ be the Breuil module with descent data and action
of $\kk_E$ given by:
$$ \M' = (\Fpp \otimes \kk_E)[u]/u^{e_2 p} \e , \quad \Fil^1 \M' = u^{e_2} \M'$$
$$ \phi_1(u^{e_2} \e) = (1\otimes \overline{a}^{-1}) \e , \quad
\ghat(\e) = (1 \otimes \omega^j(g)) \e \ \text{for} \ g \in
\Gal(F_2/\Qp).$$ Then $T_{\mathrm{st},2}^{\Qp}(\M') =
\lambda_{\overline{a}} \cdot \omega^{j+1}.$

\item Let $\M'$ be the Breuil module with descent data and action
of $\kk_E$ given by:
$$ \M' = (\Fpp \otimes \kk_E)[u]/u^{e_2 p} \e , \quad \Fil^1 \M' = u^{e_2} \M'$$
$$ \phi_1(u^{e_2} \e) = (1\otimes \overline{a}^{-1}) \e , \quad
\ghat(\e) = (1 \otimes \omega_2^m(g)) \e \ \text{for} \ g \in
\Gal(F_2/\Qpp).$$ Then $T_{\mathrm{st},2}^{\Qpp}(\M') =
(\lambda_{\overline{a}}) \,|_{G_{\Qpp}} \cdot \omega_2^{m+p+1}.$
\end{enumerate}
\end{prop}

\begin{proof} The proof is the same as that of
Proposition \ref{prop:chars1}.
\end{proof}

The Breuil modules in parts (1) and (2) of the above Proposition
are denoted by $\M_E(F_2/\Qp,e_2,\overline{a}^{-1},j)$ and
$\M_E(F_2/\Qpp,e_2,\overline{a}^{-1},m)$, respectively.

\begin{remark} When comparing Propositions \ref{prop:chars1} and
\ref{prop:chars2} with \cite[Th. 6.3]{SavittCompositio}, one
should remember that $T_{\mathrm{st},2}$ is a Tate twist of the dual
of $V_{\mathrm{st},2}$.  For example, when $E=\Qp$, the Breuil
modules in (1) of Proposition \ref{prop:chars2} are identified in
\cite[Thm 6.3]{SavittCompositio} with the character
$\lambda_{\overline{a}^{-1}} \cdot \omega^{-j}$.
\end{remark}

\begin{remark} By Lemma \ref{maximal}, the Breuil modules of
Propositions \ref{prop:chars1} and \ref{prop:chars2} are maximal.
\end{remark}

\section{Some strongly divisible modules}
\label{calcs}

In this section, we list strongly divisible modules inside the
weakly admissible filtered modules $D_{x_1,x_2}$, $D'_{x_1,x_2}$,
and $D_{m,[a:b]}$ of Propositions \ref{prin}, \ref{prin2}, and
\ref{super}, and we use them to prove the main results of our paper.

\subsection{Elements of $S$} \label{special-elements}

We begin by constructing certain elements of the rings
$S_{F_1,\OO_{E}}$ and $S_{F_2,\OO_{E}}$.  Recall the notation of
Propositions \ref{prin} and \ref{prin2}, and define $w \in
\OO_E^{\times}$ via $x_1 x_2 = pw$. Set $e_1=e(F_1/\Qp)=p-1$ and
$e_2 = e(F_2/\Qp) = p^2-1$.

\begin{lemma} \label{elements:one} Let $x \in \OO_E$.  If $j = 1$, suppose further
that $x^2 \not\equiv w \pmod{\mE}$.  Then there exists a unique
element $V_x \in S_{F_1,\OO_E}$ satisfying
\begin{equation}\label{eqC}
V_x = 1 + \frac{x^2}{w} u^{p(p-1)(j-1)}\left(\frac{u^{e_1p}}{p} +
1\right) \phi(V_x) .
\end{equation}
\end{lemma}

\begin{proof} Suppose that $ V_x = \sum_n v_n u^n $ solves \ref{eqC}.
Then for $n > 0$, $v_n$ satisfies
\begin{equation}\label{eqB}
v_n = \frac{x^2}{w}\left( v_k + \frac{v_{k-e_1}}{p} \right)
\end{equation}
where $$ kp + p(p-1)(j-1) = n$$ and $v_k$ is taken to be zero if $k$
is not a nonnegative integer. Since $n > 0$, both $k$ and $k-e_1$
are strictly smaller than $n$, and so the existence and uniqueness
of $V$ (as a formal power series) follow inductively as soon as we
know that the constant term in \eqref{eqC} can be satisfied.

If $j > 1$, the condition on $v_0$ is simply $v_0=1$.  For $j=1$,
the constant term in \eqref{eqC} is
$$ v_0 = 1 + \frac{x^2}{w} v_0 . $$
This has a solution $v_0 \in \OO_E$ exactly as long as $x^2
\not\equiv w \pmod{\mE}$.

It remains to check that $V_x$ is actually an element of
$S_{F_1,\OO_E}$. Indeed, it follows inductively from \eqref{eqB}
that if the denominator of $v_n$ has $p$-adic valuation at least
$N$, then $n \ge e_1 p (p^N - 1)/(p-1)$.  In particular
$$ v_n \in \frac{1}{p^{\lfloor n/e_1 p \rfloor}} \OO_E .$$  Since
$\frac{u^n}{p^{\lfloor n/e_1 p \rfloor}} \rightarrow 0$ in $S_E$
as $n \rightarrow \infty$ , the desired conclusion follows.
\end{proof}

Similarly, we define $U_x \in S_{F_1,\OO_E}$ satisfying
$$ U_x = 1 + \frac{x^2}{w} u^{p(p-1)(p-2-j)}\left(\frac{u^{e_1 p}}{p} +
1\right) \phi(U_x) ,$$ which exists provided that $x^2
\not\equiv w \pmod{\mE}$ in the case $j=p-2$, and is then unique.

We define analogous elements $V'_x$ and $U'_x$ in $S_{F_2,\OO_E}$
by replacing $u$ everywhere by $u^{p+1}$ (e.g., replacing $u^{e_1}$
by $u^{e_2}$). For example, $V'_x$ satisfies
$$V'_x = 1 + (1
\otimes x^2 w^{-1}) u^{pe_2 (j-1)}\left(\frac{u^{e_2 p}}{p} +
1\right) \phi(V'_x) .$$

We remark that each coefficient of $u$ in $V'_x$ and $U'_x$ is a
power series in $x$.  As a result, putting variables $X_1, X_2$ for
$x$ in $V'_x$ and $U'_x$ respectively, we obtain elements $V_{X_1},
U_{X_2} \in S_{F_2,\OO_E[[X_1,X_2]]/(X_1 X_2 - wp)}$ which
specialize to $V'_{x_1}$ and $U'_{x_2}$ under the map
$\OO_E[[X_1,X_2]]/(X_1 X_2 - wp) \rightarrow \OO_E$ sending $X_1,X_2
\mapsto x_1,x_2$ when $0 < \val_p(x_1),\val_p(x_2) < 1$. Similarly,
if $\val_p(x)=0$ put $x = \tilde{x} ( 1 + y)$ with $\tilde{x}$ the
Teichm\"uller lift of the image of $x$ in $\kk_E$. Putting
$\tilde{x}(1+Y)$ for $x$ in $V_x$ and $U_x$ respectively, we obtain
elements $V_Y, U_Y \in S_{F_1,\OO_E[[Y]]}$ which specialize to $V_x,
U_x$ under the map $\OO_E[[Y]] \rightarrow \OO_E$ sending $Y \mapsto
y$.

Next, recall the notation of Proposition \ref{super}, and define $w
\in \OO_E^{\times}$ via $x = pw$.  Write $m = i + (p+1)j$ with $i
\in \{0,\ldots,p\}$ and $j \in \Z/(p-1)\Z$. It is easy to see that
$D_{m,[a:b]}\cong D_{pm,[bw:-a]}$, so without loss of generality we
may assume that $a=1$ and $\val_p(b) \ge 0$.  Then we have the
following.

\begin{lemma} \label{elements:two} If $i < p$, there is a unique $W \in S_{F_2,\OO_E}$
satisfying
\begin{equation}
\label{eqW}
 W = -(1 \otimes w) + \left(1 +
\frac{u^{pe_2}}{p}\right) (1 \otimes b^2) W \phi(W)
u^{pe_2(p-i)}.
\end{equation}
\end{lemma}

\begin{proof} This follows inductively in the same manner as Lemma
\ref{elements:one}.  For the base case, note that since $i < p$
the constant term $w_0$ is just $-(1 \otimes w)$.
\end{proof}

When $i=p$, we must solve the identity \eqref{eqW} somewhat more
carefully. The constant term solves
\begin{equation} \label{constant-term}
w_0 = -(1 \otimes w) + (1 \otimes b^2) w_0^2 .
\end{equation}
Therefore,  as long as $\OO_E$ contains a root of the
quadratic $b^2 z^2 - z - w$ --- that is, as long as $1 + 4wb^2$ is a
square in $E$ --- the recursion can get started with $w_0 = 1 \otimes
z$.  If $\val_p(b)
> 0$, by Hensel's lemma this is always possible; taking the square
root of $1 + 4wb^2$ which is $1 \pmod{\mm_E}$, the corresponding
root $z = (1 - \sqrt{1 + 4wb^2})/2b^2 \in \OO_E$ can be expressed as
a power series in $b$.  If $\val_p(b) = 0$ and $1 + 4wb^2 \not\equiv
0 \pmod{\mm_E}$, write $b = \tilde{b} (1 + \beta)$ with
$\val_p(\beta)>0$ and $\tilde{b}$ the Teichm\"uller lift of the
image of $b$ in $\kk_E$.  Then either root $z$ of the quadratic $b^2 z^2 - z- w$
may be chosen and
expressed as a power series in $\beta$; in this case we must assume
that $1 + 4w\tilde{b}^2$ is a square in $E$.  Finally, if $1 + 4wb^2
\equiv 0 \pmod{\mm_E}$, we must assume that $1 + 4wb^2$ is a square
in $E$; in this case our root of $b^2 z^2 - z - w$ may not be
expressed as a power series in terms of $b$, but we shall see later
that this does not matter.  We obtain the following.

\begin{lemma}
\label{extend}
\begin{enumerate}
\item  If $i = p$, and if $1 + 4wb^2$ is a square in $E$
when $\val_p(b)=0$, then there is $W \in S_{F_2,\OO_E}$ satisfying
$$ W = -(1 \otimes w) + \left(1 + \frac{u^{pe_2}}{p}\right) (1
\otimes b^2) W \phi(W) u^{pe_2(p-i)}.$$

\item  If $i = 1$, $\val_p(b) > 0$, and $w$ is a
square in $E$, then there is $X \in S_{F_2,\OO_E}^{\times}$ satisfying
$$ X(1 \otimes wb) = 1\otimes w - \left(1 + \frac{u^{pe_2}}{p}\right)X\phi(X).$$
\end{enumerate}
\end{lemma}

\begin{proof} 
(1) The  paragraph before the Lemma solves for the constant term
of $W$.  The recursion for the coefficient $w_n$ of $u^n$ is
$$ w_n = (1 \otimes b^2) w_n w_0 + \textrm{lower terms}.$$
Since $w_0 = 1 \otimes z$ and $b^2 z \not\equiv 1 \pmod{\mm_E}$,
the recursion can be solved to obtain $W \in S_{F_2,\OO_E}$.

(2) The constant term of $X$ may be taken to be $1 \otimes
  x_0$ where $x_0$ is either root of $x_0^2 + wbx_0 - w$ in
  $\OO_E^{\times}$. 
 The recursion for the coefficient $x_n$ of $u^n$ is $x_n(x_0 + wb) =
  \mathrm{lower \ terms}$, and so the recursion can be solved to
  obtain $X \in S_{F_2,\OO_E}^{\times}$.  
\end{proof}

Moreover, if $\val_p(b) > 0$, then in all cases by putting the
variable $B$ for $b$ we obtain an element $W_B$ of
$S_{F_2,\OO_E[[B]]}$ which specializes to $W$ under the map
$\OO_E[[B]] \rightarrow \OO_E$ sending $B \mapsto b$.  If $\val_p(b)
= 0$ and we are away from the situation $i=p$ and $1 + 4w\tilde{b}^2
\equiv 0 \pmod{\mm_E}$, assume that $1 + 4w\tilde{b}^2$ is a square
in $E$; then by putting $\tilde{b} ( 1 + B)$ for $b$ we obtain an
element $W'_B$ of $S_{F_2,\OO_E[[B]]}$ which specializes to $W$
under the map $\OO_E[[B]] \rightarrow \OO_E$ sending $B \mapsto
\beta$. (In fact, when $\val_p(b)=0$ and $i=p$, there are two such
$W'_B$: one for each root of $b^2 z^2 - z - w = 0$.)

Similarly, if $\val_p(b) > 0$ then by putting the variable $B$ for $b$
we obtain an element $X_B$ of $S_{F_2,\OO_E[[B]]}$ which specializes
to $X$ under the map $\OO_E[[B]] \rightarrow \OO_E$ sending $B \mapsto b$.
Note that the image of $X$ in $(\Fpp \otimes \kk_E)[u]/u^{e_2 p}$ is
$1\otimes c$ with $c$ a square root of $\overline{w}$.  

\subsection{Strongly divisible modules} \label{subsec62}

With the special elements $U, V, W$ in hand, we now present the
strongly divisible modules that are contained inside the filtered
modules of Propositions \ref{prin} and \ref{super}.

First, suppose we are in the situation of Proposition \ref{prin} or
\ref{prin2}. Without loss of generality (twisting by an appropriate
character) it suffices to consider the case $i=0$.  We begin by
noting the following lemma.

\begin{lemma} \label{reduciblecases} In the two cases
\begin{itemize}
\item $\val_p(x_1)=0$, $j=1$, and $x_1^2 \equiv w \pmod{\mE}$;

\item $\val_p(x_2)=0$, $j=p-2$, and $x_2^2 \equiv w \pmod{\mE}$;
\end{itemize}
the mod $p$ reduction of the representation corresponding to
$D_{x_1,x_2}$ does not have trivial centralizer.
\end{lemma}

\begin{proof} In the first case Example
\ref{characters} tells us that the representation corresponding to
$D_{x_1,x_2}$ is an extension of $\epsilon \lambda_{x_1^{-1}}$ by
$\om \lambda_{x_1 w^{-1}}$, and the condition that $x_1^2 \equiv w
\pmod{\mE}$ forces $x_1^{-1} \equiv x_1 w^{-1} \pmod{\mE}$.
Therefore the two characters $\epsilon \lambda_{x_1^{-1}}$ and
$\om \lambda_{x_1 w^{-1}}$ have the same reduction modulo $p$. The
second case is similar.
\end{proof}

In the remainder of this section, we therefore assume that we
are not in either of the two cases of Lemma \ref{reduciblecases}.
Set $\D_{x_1,x_2} = S_{F_1,\OO_E} \otimes D_{x_1,x_2}$ if
$\val_p(x_1),\val_p(x_2)$ are integers and $\D_{x_1,x_2} =
S_{F_2,\OO_E} \otimes D'_{x_1,x_2}$ if $0 < \val_p(x_1),\val_p(x_2)
< 1$.   Then we have the following.

\begin{prop}  \label{str-prin} Put $F=F_1$ if $\val_p(x_1),\val_p(x_2)$ are integers
and $F=F_2$ if $0 < \val_p(x_1),\val_p(x_2) < 1$.  There exists a
strongly divisible $\OO_E$-module with descent data
$$\M_{x_1,x_2} = S_{F,\OO_E} \cdot g_1 + S_{F,\OO_E} \cdot g_2$$
inside $\D_{x_1,x_2}$, where:
\begin{enumerate}

\item if $\val_p(x_1)=0$ and $\val_p(x_2)=1$, then
\begin{eqnarray*}
g_1 & = &-x_1 \e_1 \\
g_2 & = & \e_2 + \frac{x_1^2}{w} \frac{u^{pj-e_1}}{p} (u^{e_1} +
p) V_{x_1} \e_1 \ ;
\end{eqnarray*}

\item if $\val_p(x_1)=1$ and $\val_p(x_2)=0$, then
\begin{eqnarray*}
g_1 & = & -x_1 \e_1 + x_2 \frac{u^{p(e_1-j)-e_1}}{p}(u^{e_1} + p) U_{x_2} \e_2 \\
g_2 & = & \e_2 \ ;
\end{eqnarray*}

\item if $0 < \val_p(x_1),\val_p(x_2) < 1$, then if $k=(p+1)j$,
\begin{eqnarray*}
g_1 & = & -x_1 \e_1 + x_2 \frac{u^{p(e_2-k)-e_2}}{p}(u^{e_2} +
p)U'_{x_2}
\e_2 \\
g_2 & = & \e_2 + \frac{x_1^2}{w}\frac{u^{pk-e_2}}{p}(u^{e_2}+
p)V'_{x_1} \e_1 \ .
\end{eqnarray*}
\end{enumerate}
\end{prop}

\begin{proof} Abbreviate $\M = \M_{x_1,x_2}.$  In each case, the only nontrivial steps are to
compute $\Fil^1 \M$, to verify that it satisfies $\Fil^1 \M \cap
I\M = I \Fil^1 \M$, and to check that $\phi(\Fil^1 \M)$ lies
inside $p\M$ and generates it over $S_{F,\OO_E}$ or, equivalently,
that $\p(\Fil^1\M)$ lies inside $\M$ and generates it over
$S_{F,\OO_E}$. Note that in each case, $g_1$ and $g_2$ are both
eigenvectors for the action of $\Gal(F_1/\Qp)$ (resp.,
$\Gal(F_2/\Qp)$).

We begin with case (1), in which $\val_p(x_1)=0$.  It is easy to
check that
$$\Fil^1 \M = S_{F_1,\OO_E} \cdot (-u^j g_1 + x_1 g_2) + (\Fil^1
S_{F_1,\OO_E})\M,$$ that $\phi(g_1) = x_1 g_1$, and using the
defining equation for $V_{x_1}$ from Lemma \ref{elements:one},
that
$$ \phi(g_2) = x_2 g_2 + u^{pj-e_1} (u^{e_1} + p V_{x_1}) g_1
.$$  From this it follows that $$\phi(-u^j g_1 + x_1 g_2) =
p(wg_2 + x_1 u^{pj-e_1} V_{x_1} g_1),$$ and we see easily from this
that $\phi(\Fil^1 \M) \subset p\M$ and generates it.  The fact
that $\Fil^1 \M \cap I\M = I \Fil^1 \M$ follows without difficulty
from the analogous fact for $S_{F_1,\OO_E}$.

Similarly, in case (2), in which $\val_p(x_2)=0$, we have
$$\Fil^1 \M = S_{F_1,\OO_E} \cdot(x_2 g_1 + wu^{e_1-j} g_2) +
(\Fil^1 S_{F_1,\OO_E})\M .$$  We see that $\phi(g_2) = x_2 g_2$
and, by the defining equation for $U_{x_2}$, that
$$\phi(g_1) = x_1 g_1 - w u^{p(e_1 - j) - e_1} (u^{e_1} + p
U_{x_2}) g_2.$$  It follows that $$\phi(x_2 g_1 + wu^{e_1-j} g_2) =
p(wg_1 - wx_2 u^{p(e_1-j)-e_1} U_{x_2} g_2),$$ and the other
properties of $\M$ follow as above.

Finally, we turn to case (3), where $0 < \val_p(x_1),\val_p(x_2) <
p$.  We note that if polynomials $s(u),t(u)$ over $W(k) \otimes
\OO_E$ are such that $(1 \otimes x_1) s + u^k t$ is divisible by
$u^{e_2} + p$, then $(s,t)$ is a linear combination of $(-u^k,1
\otimes x_1)$ and $(1 \otimes x_2,(1\otimes w)u^{e_2-k})$. It
follows that $\Fil^1 \M$ is the submodule of $\M$ generated by
$-u^k g_1 + (1 \otimes x_1) g_2$, $(1 \otimes x_2) g_1 + (1
\otimes w)u^{e_2-k} g_2$, and $(\Fil^1 S_{F_2,\OO_E})\M$.
Moreover, if $(s,t) = \alpha(-u^k,1 \otimes x_1) + \beta(1 \otimes
x_2,(1 \otimes w)u^{e_2-k})$ and the coefficients of $s,t$ are in
$I$, then so are the coefficients of $\alpha,\beta$, and so
$I\Fil^1 \M = I\M \cap \Fil^1 \M$.  It remains to compute
$\phi(g_1)$ and $\phi(g_2)$, and to verify that $\phi(-u^k g_1 +
(1 \otimes x_1) g_2)$ and $\phi((1 \otimes x_2) g_1 + (1 \otimes
w)u^{e_2-k} g_2)$ lie in $p\M$.

Set
$$ D = \left(1 + U'_{x_2} V'_{x_1} \left(\frac{u^{e_2 p}}{p} +
2u^{(p-1)e_2} + pu^{(p-2)e_2}\right)\right) ,$$ an invertible
element of $S_{F_2,\OO_E}$.   Inverting the matrix that yields
$g_1$ and $g_2$ in terms of $(1 \otimes x_1) \e_1$ and $\e_2$
gives
\begin{eqnarray*}
(1 \otimes x_1) \e_1 & = & D^{-1} \left( -g_1 + (1 \otimes x_2)
\frac{u^{p(e_2-k)-e_2}}{p}(u^{e_2}+p)U'_{x_2} g_2 \right) \\
\e_2 & = & D^{-1} \left((1 \otimes x_1 w^{-1})
\frac{u^{pk-e_2}}{p} (u^{e_2} + p) V'_{x_1} g_1 + g_2\right)
\end{eqnarray*}
Substituting into the expressions
\begin{eqnarray*}
\phi(g_1) & = & (1 \otimes x_1) g_1 - (1 \otimes w)
u^{p(e_2-k)-e_2}(u^{e_2} + pU'_{x_2})
\e_2 ,\\
\phi(g_2) & = & (1 \otimes x_2) g_2 - u^{pk - e_2}(u^{e_2} +
pV'_{x_1}) ((1 \otimes x_1) \e_1)
\end{eqnarray*}
(which are obtained using the defining equations for $V'_{x_1},
U'_{x_2}$) and simplifying yields
\begin{eqnarray*}
\phi(g_1) & = & (1 \otimes x_1) D^{-1} \left(1 +
\left(\frac{u^{e_2 p}}{p} +
u^{(p-1)e_2}\right)V'_{x_1} (U'_{x_2}-1)\right)g_1 \\
 & & - (1 \otimes w)D^{-1} u^{p(e_2-k)-e_2}(u^{e_2} + pU'_{x_2}) g_2 ,\\
\phi(g_2) & = & D^{-1} u^{pk-e_2} (u^{e_2} + pV'_{x_1})g_1 \\
 & & + (1 \otimes x_2) D^{-1} \left(1 + \left(\frac{u^{e_2 p}}{p} +
u^{(p-1)e_2}\right)U'_{x_2} (V'_{x_1}-1)\right)g_2 .
\end{eqnarray*}
This confirms that $\phi(g_1),\phi(g_2) \in \M$.  We then compute
$\phi(-u^k g_1 + (1 \otimes x_1) g_2)$ to be $pD^{-1}$ times
\begin{eqnarray*}
 & u^{pk-e_2} V'_{x_1} \left(
(1 \otimes x_1) - (1 \otimes x_2) (u^{e_2}+p) \left(\frac{u^{e_2
p}}{p} +
1\right)\frac{u^{pe_2 (p-1-j)}}{p} \phi(U'_{x_2})\right) g_1\\
 + & (1\otimes w) \left( 1 + \left(\frac{u^{e_2 p}}{p} +
 u^{(p-1)e_2}\right) U'_{x_2} V'_{x_1} + \frac{u^{e_2 p}}{p} \left(1
 - U'_{x_2}\right)\right) g_2
\end{eqnarray*}
and $\phi(x_2 g_1 + w u^{e_1 - j} g_2)$ to be $p(1\otimes
w)D^{-1}$ times
\begin{eqnarray*}
 & \left( 1 + \left(\frac{u^{e_2 p}}{p} +
 u^{(p-1)e_2}\right)U'_{x_2} V'_{x_1} + \frac{u^{e_2 p}}{p}(1 -
 V'_{x_1})\right) g_1 \\
+ & u^{p(e_2-k)-e_2} U'_{x_2} \left( -(1\otimes x_2) + (1 \otimes
x_1) (u^{e_2}+p) \left(\frac{u^{e_2 p}}{p} + 1\right) \frac{u^{p
e_2 j}}{p} \phi(V'_{x_1})\right) g_2.
\end{eqnarray*}
In each case, the image lies inside $p\M$.  Moreover, one checks
without difficulty (by working modulo $u$) that $\p(-u^k g_1 + x_1
g_2)$ and $\p(x_2 g_1 + w u^{e_2 - k} g_2)$ generate $\M$ over
$S_{F_2,\OO_E}$.  This completes the proof.
\end{proof}

We turn next to the strongly divisible modules in the situation of
Proposition \ref{super}.  (Happily, this is actually simpler
than the previous situation.)  Extend $E$ if necessary (i.e., when
required by Lem. \ref{extend}) to assume that $1 + 4wb^2$ or $w$ is a
square in $E$, and write $\D_{m,[1:b]} = S_{F_2,\OO_E} \otimes
D_{m,[1:b]}$.  (Recall that we have without loss of generality
assumed $a=1$ and $\val_p(b) \ge 0$.) Set $k=(p-1)i$.  We then have
the following.

\begin{prop} \label{str-super} There exists a strongly divisible $\OO_E$-module with
descent data $$ \M_{m,[1:b]} = S_{F_2,\OO_E} \cdot g_1 +
S_{F_2,\OO_E} \cdot g_2 $$ inside $\D_{m,[1:b]}$, where if $i > 1$
or $\val_p(b)=0$ then 
\begin{eqnarray*}
g_1 & = & \e_1, \\
g_2 & = & \frac{\e_2}{p}+ (1\otimes b) W \frac{u^{p(e_2-k)}}{p}
\e_1,
\end{eqnarray*}
while if $i=1$ and $\val_p(b)>0$ then we define instead
\begin{eqnarray*}
g_1 & =& \mathbf{e}_1 + \frac{X}{pw} u^{p(p-1)} \mathbf{e}_2 \\
g_2 & = & \mathbf{e}_2.
\end{eqnarray*}
\end{prop}

We remark that the first set of formulas for $g_1,g_2$ will still define a
strongly divisible module when $i=1$ and $\val_p(b)>0$; however, it is not the
strongly divisible module that we wish to consider later on.

\begin{proof}  Put $\M = \M_{m,[1:b]}$.  Suppose first that $i > 1$ or
  $\val_p(b)=0$.  We begin by noting that
$$ \left(-\e_1 + \frac{1 \otimes b}{p} u^{e_2-k} \e_2\right)
+ \left(\frac{u^{e_2(p-i)}}{p}(1\otimes
b^2)W\right)\left(u^{e_2}+p\right) \e_1$$ is equal to
$$ (1 \otimes b)u^{e_2-k} g_2 + (u^{e_2(p-i)} (1\otimes b^2)W -
1)g_1,$$ and so this element of $\M$ lies in $\Fil^1\M$.  We remark
that $u^{e_2(p-i)} (1\otimes b^2)W - 1$ is a unit in
$S_{F_2,\OO_E}$: this is clear when $i < p$; when $i=p$ use
\eqref{constant-term} to see that $b^2 w_0 - 1 \not\equiv 0
\pmod{\mm_E}$.  Noting that $g_2$ is not an element of $\Fil^1 \M$
(when $i=p$, this again uses the fact that $b^2 w_0 - 1 \not\equiv
0 \pmod{\mm_E}$) we find that
$$\Fil^1 \M = S_{F_2,\OO_E} \cdot ((1 \otimes b)u^{e_2-k} g_2 + (u^{e_2(p-i)} (1\otimes b^2)W -
1)g_1) + (\Fil^1 S_{F_2,\OO_E})\M.$$ From this, it is easy to
check that $I\M \cap \Fil^1\M = I\Fil^1 \M$.  It remains to
compute $\phi(g_1)$ and $\phi(g_2)$, and to verify that $\phi((1
\otimes b)u^{e_2-k} g_2 + (u^{e_2(p-i)} (1\otimes b^2)W - 1)g_1)$
lies in $p\M$.  Indeed
$$ \phi(g_1) = \e_2 = p g_2 - b W u^{p(e_2-k)} g_1 $$
and
$$ \phi(g_2) = \left(w - (1\otimes b^2) W\phi(W)
\frac{u^{pe_2(p+1-i)}}{p}\right) g_1 + b\phi(W)
u^{p^2(e_2-k)}g_2.$$ Then, after significant cancellation and
using the defining equation for $W$ from Lemma \ref{elements:two}
(when $i < p$) or Lemma \ref{extend} (when $i=p$), we find
$$ \phi((1
\otimes b)u^{e_2-k} g_2 + (u^{e_2(p-i)} (1\otimes b^2)W - 1)g_1) =
pw W^{-1} g_2 .$$ For future reference, we record that
$\phi((u^{e_2} + p)g_2)$ is equal to $p\left(\frac{u^{e_2 p}}{p} +
1\right)$ times
$$
\left(\left((1\otimes w) - (1\otimes b^2) W\phi(W)
\frac{u^{pe_2(p+1-i)}}{p}\right)g_1 + (1 \otimes b) \phi(W)
u^{p^2(e_2-k)} g_2\right).
$$
In particular, the coefficient of $g_1$ in this expression is a
unit in $S_{F_2,\OO_E}$, so $\phi_1(\Fil^1\M)$ does generate $\M$
over $S_{F_2,\OO_E}$.

Now suppose instead that $i=1$ and $\val_p(b) > 0$.  Observe that $h := u^{p-1} g_1 +
  \left(\frac{X}{w} + (1\otimes b)\right)g_2$ lies in $\Fil^1 \M$.
  Since $\frac{X}{w} + (1 \otimes b)$ is a unit in $S_{F_2,\OO_E}$ and
  $g_1$ does not lie in $\Fil^1 \M$, we deduce that $\Fil^1 \M =
  S_{F_2,\OO_E} \cdot h + (\Fil^1 S_{F_2,\OO_E})\M$.  From this it is
  easy to check that $ I\M \cap \Fil^1 \M = I \Fil^1 \M$.  Finally, we
  compute that
\begin{eqnarray*}
\phi(g_1) & = & \phi(X) u^{p^2(p-1)} g_1 + \left(1-X\phi(X)\frac{u^{pe_2}}{pw}\right) g_2 \\
\phi(g_2) & = & pwg_1 - X u^{p(p-1)} g_2
\end{eqnarray*}
both lie in $\M$; using the defining relation for $X$ we find
$\phi_1(h) = (1 \otimes w) X^{-1} g_1 \in \M$ and conclude that $\M$
is a strongly divisible module.
\end{proof}

\subsection{Reduction mod $\mm_E$} \label{subsec63} For each of the strongly
divisible modules $\M$ of Section \ref{subsec62}, corresponding to a
lattice in a Galois representation, we compute the reduction modulo
$\mm_E$ of that lattice; that is, we compute
$T_{\mathrm{st},2}^{\Qp}(\M/\mm_E)$.

Suppose first that we are in the situation of Propositions
\ref{prin} and \ref{prin2}, excluding the cases of Lemma
\ref{reduciblecases}. We have the following.

\begin{thm}  \label{red-prin} Let $\M=\M_{x_1,x_2}$ be one of the strongly divisible
modules of Proposition \ref{str-prin}.  Then we have the following.
\begin{enumerate}
\item If $\val_p(x_1) = 0$, then $T_{\mathrm{st},2}^{\Qp}(\M/\mm_E)$
depends only on the reduction $\ox_1$ of $x_1 \pmod{\mm_E}$, and
in fact,
$$
T_{\mathrm{st},2}^{\Qp}(\M/\mm_E) \cong
\begin{pmatrix}
\lambda_{\ox_1^{-1}} \omega & * \\
0 & \lambda_{\ox_1\ow^{-1} } \omega^{j}
\end{pmatrix}
$$
with $* \neq 0$.

\item If $\val_p(x_2) = 0$, then $T_{\mathrm{st},2}^{\Qp}(\M/\mm_E)$ depends only on
the reduction $\ox_2$ of $x_2 \pmod{\mm_E}$, and in fact,
$$
T_{\mathrm{st},2}^{\Qp}(\M/\mm_E) \cong
\begin{pmatrix}
\lambda_{\ox_2^{-1}} \omega^{1+j} & * \\
0 & \lambda_{\ox_2\ow^{-1} }
\end{pmatrix}
$$
with $* \neq 0$.

\item If $0 < \val_p(x_1),\val_p(x_2) < 1$, then
$T_{\mathrm{st},2}^{\Qp}(\M/\mm_E)$ is independent of $x_1$ and
$x_2$ and satisfies
$$ T_{\mathrm{st},2}^{\Qp}(\M/\mm_E) \,|_{I_p} \otimes_{\kk_E} \Fpbar \cong \omega_2^{1+j} \oplus
\omega_2^{p(1+j)} .$$
\end{enumerate}
\end{thm}

\begin{proof}
(1) By inspection, the Breuil module $\M' = T_0(\M/\mm_E)$ is
generated by $g_1$ and $g_2$ over $\kk_E[u]/u^{e_1 p}$ with
$$\Fil^1 \M' = \kk_E[u]/u^{e_1 p} \cdot (-u^j g_1 + \ox_1 g_2) +
\kk_E[u]/u^{e_1 p} \cdot (u^{e_1} g_1) ,$$
$$ \p(-u^j g_1 + \ox_1 g_2) = \ow g_2 + \ox_1 u^{pj-e_1}
\overline{V}_{x_1} g_1,$$ and $\p(u^{e_1} g_1) = \ox_1 g_1$.
Also, $\ghat(g_1) = g_1$ and $\ghat(g_2) = \om^j(g)(g_2)$.

Let $\M_1 = \M_E(F_1/\Qp,e_1,\ox_1,0)$. It follows from Proposition
\ref{prop:chars1} that $T_{\mathrm{st},2}^{\Qp}(\M_1) =
\lambda_{\ox_1^{-1}} \omega$.  Let $\M_2 = \M_E(F_1/\Qp,e_1,\ow
\ox_1^{-1},j-1)$.   By Proposition \ref{prop:chars2}, we have
$$T_{\mathrm{st},2}^{\Qp}(\M_2) = \lambda_{\ox_1 \ow^{-1}} \omega^j.$$

But it is clear that $\M'$ has a submodule that is isomorphic to
$\M_1$.  Moreover, there is a map from $\M' \rightarrow \M_2$
sending $g_1 \mapsto 0$ and $g_2 \mapsto u^{p} \e$.  It follows that
$T_{\mathrm{st},2}^{\Qp}(\M')$ has the desired form, and to see that
$* \neq 0$, it suffices by Lemma \ref{maximal} to check that there
is no nontrivial map $\M' \rightarrow \M_1$.  This is a standard
calculation (that uses crucially the assumption that $\ow \neq
\ox_1^2$ when $j=1$).

(2) This is similar to (1).

(3) Extend $E$ so that it contains $\Qpp$ and so that $\kk_E$
contains a square root of $\ow$.  (We see the reason for the
latter assumption towards the end of the argument.)   Note that
$\overline{U}'_{x_2} = \overline{V}'_{x_1} = 1$ in $\Fpp \otimes
\kk_E[u]/u^{e_2 p}$, so that $\overline{D} = 1 + 2 u^{(p-1)e_2}$.
Therefore the Breuil module $T_0(\M/\mm_E)$ is $$\M' = (\Fpp
\otimes \kk_E)[u]/u^{e_2 p} \cdot g_1 \oplus (\Fpp \otimes
\kk_E)[u]/u^{e_2 p} \cdot g_2 $$ with $\Fil^1 \M'$ generated over
$(\Fpp \otimes \kk_E)[u]/u^{e_2 p}$ by $-u^k g_1$ and $(1 \otimes
\ow) u^{e_2-k}g_2$ with
 $$\p(-u^k g_1) = (1 \otimes \ow) \overline{D}^{-1} (1 + u^{(p-1)e_2}) g_2 ,$$
 $$ \p((1\otimes \ow)u^{e_2-k}g_2) = (1 \otimes \ow) \overline{D}^{-1}
 (1 +
 u^{(p-1)e_2})g_1 ,$$
$$ \ghat(g_1) = g_1 , \ \ \ghat(g_2) = (\omt(g)^k \otimes 1)
g_2.$$ Replacing $g_1$ by $\overline{D}^{-1} (1 +
u^{(p-1)e_2})g_1$ and $g_2$ by $-\overline{D}^{-1} (1 +
u^{(p-1)e_2})g_2$ simplifies the form of the filtration and
Frobenius to:
$$\Fil^1 \M' = \Fpp \otimes \kk_E[u]/u^{e_2 p} \cdot (u^k g_1) +
\Fpp \otimes \kk_E[u]/u^{e_2 p} \cdot (u^{e_2-k}g_2) $$ with
$$\p(u^k g_1) = (1 \otimes \ow) g_2 , \ \ \p(u^{e_2-k} g_2) = -g_1.$$

Restrict the descent data on $\M'$ to $\Gal(F_2/\Qpp)$, which
amounts to restricting the representation
$T_{\mathrm{st},2}^{\Qp}(\M')$ to $G_{\Qpp}$.  Denote this new
Breuil module by $\M'_2$. Let $\M'' = \M_E(F_2/\Qpp, e_2, c, n)$.
One checks that there is a nontrivial map from $\M'_2 \rightarrow
\M''$ given by $$ g_1 \mapsto u^{p(p-j)} \alpha \e, $$
$$ g_2 \mapsto u^{p(1+j)} \beta \e,$$
provided that
\begin{itemize}
\item $\phi(\beta) (1\otimes c) = -\alpha$,

\item $\phi(\alpha) (1\otimes c) = (1 \otimes \ow) \beta$,

\item $(\omt^{p(j-p)}\otimes 1) \alpha = (1 \otimes \omt^{n})
\alpha$,

\item $(\omt^{j-p} \otimes 1 ) \beta = (1 \otimes \omt^{n})\beta$.
\end{itemize}
Then it is possible to satisfy the above conditions with $c =
\sqrt{-\ow}$ and either $n=p(j-p)$ or $n=j-p$: in the former case,
take $\alpha \in \Fpp \otimes \kk_E$ which is annihilated by
$(\omt^{p(j-p)} \otimes 1) - (1 \otimes \omt^{p(j-p)})$, and in the
latter case, take $\alpha$ which is annihilated by $(\omt^{p(j-p)}
\otimes 1) - (1 \otimes \omt^{j-p}).$   By Proposition \ref{prop:chars2}, it
follows that $$T_{\mathrm{st},2}^{\Qp}(\M') \,|_{G_{\Qpp}} \cong
\lambda_{\sqrt{-\ow}^{-1}} \,|_{G_{\Qpp}} \otimes ( \om_2^{j+1}
\oplus \om_2^{p(1+j)}) .$$  The result follows.
\end{proof}

In the situation of Proposition \ref{super}, we have the following.

\begin{thm}  \label{red-super} Let $\M = \M_{m,[1:b]}$ be one of the strongly divisible
modules of Proposition \ref{str-super}.  Then we have the following.
\begin{enumerate}
\item  If $\val_p(b)=0$ and $1 < i < p$, then
$T_{\mathrm{st},2}^{\Qp}(\M/\mm_E)$ depends only on the reduction
$\ob$ of $b \pmod{\mm_E}$, and
$$
T_{\mathrm{st},2}^{\Qp}(\M/\mm_E) \cong
\begin{pmatrix}
\lambda_{\overline{bw}^{-1}} \omega^{i+j} & * \\
0 & \lambda_{-\ob } \omega^{1+j}
\end{pmatrix}
$$
with $* \neq 0$ and peu ramifi\'e if $i=2$.

\item  If $\val_p(b)=0$ and $i=1$, then
$T_{\mathrm{st},2}^{\Qp}(\M/\mm_E)$ depends only on the reduction
$\ob$ of $b \pmod{\mm_E}$, and
$$
T_{\mathrm{st},2}^{\Qp}(\M/\mm_E) \otimes_{\kk_E} \Fpbar \cong
\begin{pmatrix}
\lambda_{r_{+}} \omega^{1+j} & * \\
0 & \lambda_{r_{-}} \omega^{1+j}
\end{pmatrix}
$$
where $r_{\pm} =  -\frac{1}{2}(\ob \pm \sqrt{\ob^2 + 4\ow^{-1}})$
and $* = 0$ if $r_+ \neq r_-$.  In any case
$T_{\mathrm{st},2}^{\Qp}(\M/\mm_E)$ does not have trivial
endomorphisms.

\item  If $\val_p(b)=0$ and $i=p$, then
$T_{\mathrm{st},2}^{\Qp}(\M/\mm_E)$ depends only on the reduction
$\ob$ of $b \pmod{\mm_E}$, and
$$
T_{\mathrm{st},2}^{\Qp}(\M/\mm_E) \cong
\begin{pmatrix}
\lambda_{\overline{bw_{-}/w}} \omega^{1+j} & * \\
0 & \lambda_{\overline{bw_{+}/w}} \omega^{1+j}
\end{pmatrix}
$$
where $w_{+}$ is the root of $b^2 z^2 - z - w = 0$ such that the
constant term of $W$ is $-(1\otimes w_+)$, and $w_-$ is the other
root. If $w_+ \not\equiv w_- \pmod{\mm_E}$ (i.e., if $1 + 4b^2 w
\not \equiv 0 \pmod{\mm_E}$), then $* \neq 0$; the two choices for
$W$ give lattices with different reductions. If $1 + 4b^2 w \equiv
0 \pmod{\mm_E}$, then $* =0 $.

\item If $i>1$ and $\val_p(b) > 0$, then $T_{\mathrm{st},2}^{\Qp}(\M/\mm_E)$ is
independent of $b$ and
$$T_{\mathrm{st},2}^{\Qp}(\M/\mm_E) \,|_{I_p} \otimes_{\kk_E} \Fpbar \cong \omega_2^{m+p}
\oplus \omega_2^{pm+1} .$$

\item If $i=1$ and $\val_p(b) > 0$, then
  $T_{\mathrm{st},2}^{\Qp}(\M/\mm_E)$ is independent of $b$ and 
$$
T_{\mathrm{st},2}^{\Qp}(\M/\mm_E) \cong
\begin{pmatrix}
\lambda_{-c^{-1}} \omega^{1+j} & * \\
0 & \lambda_{c^{-1} } \omega^{1+j}
\end{pmatrix}
$$ with $* \neq 0$.
\end{enumerate}
\end{thm}

\begin{proof}  The Breuil module $\M' = T_0(\M/\mm_E)$ in all cases
satisfies
$$ \M' = (\Fpp \otimes \kk_E)[u]/u^{e_2 p} \cdot g_1 \oplus (\Fpp \otimes \kk_E)[u]/u^{e_2
p} \cdot g_2$$ with $\Fil^1 \M'$ generated by $(1\otimes
\ob)u^{e_2-k} g_2 + (u^{e_2(p-i)}(1 \otimes \ob^2)\overline{W} -
1) g_1$ and $u^{e_2} g_2$, with
$$ \p((1\otimes \ob)u^{e_2-k} g_2
+ (u^{e_2(p-i)}(1 \otimes \ob^2)\overline{W} - 1) g_1) = (1
\otimes \ow) \overline{W}^{-1} g_2$$ and $$\p(u^{e_2} g_2) = (1
\otimes \ow)g_1 + (1\otimes
\ob)\phi(\overline{W})u^{p^2(e_2-k)}g_2$$ and $\ghat(g_1) =
(\omt^m(g) \otimes 1)g_1$, $\ghat(g_2) = (\omt^{pm}(g) \otimes
1)g_2$.

Suppose first that $\val_p(b)=0$.  If $i<p$, then $\overline{W} = -
1 \otimes \ow$ and $p^2(e_2 - k) > pe_2$.  Set $X = 1 +
u^{e_2(p-i)}(1 \otimes \ob^2) \ow$, and observe that $\phi(X) =
1$, so that $\p$ simplifies to:
$$ \p(g_1 - (1\otimes \ob) X^{-1} u^{e_2-k} g_2) = g_2$$
$$ \p(u^{e_2} g_2) = (1\otimes \ow) g_1.$$
Write $g'_1 = g_1 + C u^{kp} g_2$.  Observing that
$$ u^k g'_1 = u^k(g_1 - (1\otimes \ob)X^{-1} u^{e_2 -k} g_2)
+ (C u^{e_2 (i-1)} + (1\otimes \ob)X^{-1}) u^{e_2} g_2, $$ we
obtain $\p(u^k g'_1) = (1 \otimes \ow)(\phi(C)u^{e_2 p (i-1)} +
(1\otimes \ob)) g'_1$, provided that
$$ (1 \otimes \ow)(\phi(C)u^{e_2 p (i-1)} + (1\otimes
\ob)) C = 1 .$$ If $1 < i < p$, this is satisfied with $C = (1
\otimes \overline{bw})^{-1}$.  If $i = 1$, this is satisfied with
$C$ equal to either root of $c^2 + \ob c - \ow^{-1} = 0$,
extending $E$ if necessary to ensure that this equation has roots
in $\kk_E$.

If $1 < i <p$, this shows that $\M'$ has a sub-Breuil module
$\M''$ generated by $g_1'$ with $\Fil^1 \M'' = u^k \M''$
satisfying $\p(u^k g_1') = \overline{bw} g_1'$ and $\ghat(g_1') =
(\omt^m(g) \otimes 1) g_1'$.  Since there is a map
$$\M'' \rightarrow \M_E(F_2/\Qp,e_2,\overline{bw},i+j-1)$$
obtained by sending $g_1' \mapsto u^{p(p+1-i)} \e$, we see that
$T_{\mathrm{st},2}^{\Qp}(\M')$ has a subcharacter equal to
$\lambda_{\overline{bw}^{-1}} \omega^{i+j}$.  By considering the
determinant, the quotient character must be $\lambda_{-\overline{b}}
\omega^{1+j}$; alternately, one may check that there is a nontrivial
map from $\M' \rightarrow \M_E(F_2/\Qp,e_2,-\ob^{-1},j)$ (sending
$g_2 \mapsto u^{pi}\e$ and $g_1 \mapsto -\overline{bw}^{-1}
u^{p^2i}\e$). Finally, to see that $* \neq 0$, by Lemma
\ref{maximal} one checks that there is no map $\M' \rightarrow
\M_E(F_2/\Qp,e_2,\overline{bw},i+j-1)$ (assume such a map exists,
and use the commutativity with $\p$ and $\ghat$ to see that this
implies $i=p$).  The peu ramifi\'e claim follows from Lemma
\ref{peu-ram}.

On the other hand, if $i=1$ and the roots of $c^2 + \ob c - \ow^{-1}
=0 $ are distinct, this shows that $\M'$ has \textit{two} sub-Breuil
modules, hence two distinct subcharacters equal to
$\lambda_{r_{\pm}} \omega^{1+j}$, where $r_{\pm} = -\frac{1}{2}(\ob
\pm \sqrt{\ob^2 + 4\ow^{-1}})$.  It follows that the representation
is split.  If the roots of $c^2 + \ob c - \ow^{-1}=0$ are equal (i.e., if $4 + \ow \ob^2 = 0$),
then we only obtain one subcharacter,
equal to $\lambda_{-\ob/2} \omega^{1+j}$. But then, by considering
the determinant, we see that the quotient character is the same as
the subcharacter (since $(-\ob/2)^2 = -\ow^{-1}$) and so
$T_{\mathrm{st},2}^{\Qp}(\M')$ does not have trivial endomorphisms.

Now consider $\val_p(b)=0$ and $i=p$.  Here $\overline{W}=1\otimes
\overline{w}_+$, where $w_+$ is a chosen root of $b^2 w_+^2 - w_+ -
w = 0$. Let $w_-$ be the other root.  Since $b^2 w_+ - 1 = w/w_+$,
setting $\beta = 1 \otimes bw_+/w \pmod{\mm_E}$ we have $\p(g_1 +
\beta u^{p-1} g_2) = g_2$. Set $g_1' = g_1 + \beta u^{p^2 (p-1)}
g_2$. Since $p^2 (p-1) \ge 2e_2$, we see that $\Fil^1 \M'$ is
generated by $u^{e_2} g_2$ and $g'_1 + \beta u^{p-1} g_2$ with
$$\p(u^{e_2} g_2) = (1 \otimes \ow) g'_1 , \ \ \p(g'_1 + \beta
u^{p-1} g_2) = g_2.$$ Setting $g''_1 =  g'_1 - (1 \otimes
\overline{bw}_0^{-1}) u^{p^2 (p-1)} g_2$, one computes that
$\p(u^{p(p-1)} g''_1) = -(1 \otimes \overline{bw}_+) g''_1$.
Therefore $\M'$ has a sub-Breuil module generated by $g''_1$ with
$\Fil^1 \M'' = u^{p(p-1)} g''_1$.   There is a map $\M''
\rightarrow \M_E(F_2/\Qp,e_2,-\overline{bw}_+,j)$ sending $g''_1
\mapsto u^p \e$, so the subcharacter is
$\lambda_{\overline{bw_{-}/w}} \omega^{1+j}$.  Considering the
determinant, the quotient character is
$\lambda_{\overline{bw_{+}/w}} \omega^{1+j}$. Finally, one checks
when there exists a map $\M' \rightarrow
\M_E(F_2/\Qp,e_2,-\overline{bw}_+,j)$: one sees that such a map
must be of the form $g_2 \mapsto u^{p^2} \e$ and $g_1 \mapsto 0$.
This commutes with $\p$ on $g_1 + \beta u^{p-1} g_2$ if and only
if $-\beta^2 \ow = 1$, which occurs if and only if $1 + 4b^2 w
\equiv 0 \pmod{\mm_E}$.  In particular, $* \neq 0$ if $1 + 4b^2 w
\not\equiv 0 \pmod{\mm_E}$.  This settles part (3).

In part (4), the hypothesis that $\val_p(b) > 0$ simplifies $\M'=
T_0(\M/\mm_E)$ dramatically: namely, $\Fil^1\M'$ is generated by
$g_1$ and $u^{e_2} g_2$ with $\p(g_1) = g_2$ and $\p(u^{e_2} g_2) = (1
\otimes \ow) g_1$.  The identification of
$T_{\mathrm{st},2}^{\Qp}(\M')$ proceeds as in case (3) of Theorem
\ref{red-prin}.  In particular, let $\M'_2$ denote $\M'$ with the
descent data restricted to $\Gal(F_2/\Qpp)$.  Then a map $\M'_2
\rightarrow \M_E(F_2/\Qpp,e_2,c,n)$ must be of the form $g_1 \mapsto
\alpha u^{p^2} \e$ and $g_2 \mapsto \phi(\alpha) (1 \otimes c) u^{p}
\e$, and such a map exists if and only if $c$ is a square root of
$\ow$ and $\alpha$ is annihilated by $(\omt^{m-1} \otimes 1) - (1
\otimes \omt^{n})$. Extending $E$ if necessary so that $\ow$ has a
square root in $\kk_E$, such a map then exists for $n = m-1$ and for
$n = p(m-1)$. In the former case we get the character $(\lambda_{c^{-1}})
\,|_{G_{\Qpp}} \omega_2^{m+p}$, and in the latter case we get the
character $(\lambda_{c^{-1}}) \,|_{G_{\Qpp}} \omega_2^{pm+1}$. The result
follows.

For part (5), write $\M' = T_0(\M/\mm_E)$.  Then $\Fil^1 \M'$ is
  generated by  $u^{p-1} g_1 + c^{-1} g_2$ and $u^{e_2} g_1$, with
  $\phi_1(u^{p-1} g_1 + c^{-1} g_2) = cg_1$ and  $\phi_1(u^{e_2} g_1)
= u^{p^2(p-1)} cg_1 + g_2$.  Note that $\phi_1(u^{p(p-1)} g_2) = -cg_2$.  There is evidently a
 nontrivial map $\M' \rightarrow \M_E(F_2/\Q_p, e_2, c,j)$ sending
 $g_2 \mapsto 0$ and $g_1 \mapsto u^{p^2} \mathbf{e}$.  On the other
 hand if $f: \M' \rightarrow \M_E(F_2/\Qp,e_2,d,n)$ is a nontrivial map
 sending $g_1 \mapsto \alpha \mathbf{e}$ and $g_2 \mapsto \beta
 \mathbf{e}$, then $\alpha, \beta$ must both be polynomials in $u^p$
 since $g_1,g_2$ are in the image of $\phi_1$. On the other hand if
 $\beta\neq 0$ then the relation $f\circ \phi_1 = \phi_1 \circ f$ on $u^{p(p-1)}g_2$
implies that $\beta$ is a unit times $u^p$; but then $f(u^{p-1}g_1 +
c^{-1}g_2) \in \langle u^{e_2} \mathbf{e} \rangle$ implies that
$\alpha$ has a linear term, a contradiction.  Therefore $\beta=0$, and
then it is easy to check that $c=d$ and $j=n$.  It follows that $*
\neq 0$.
\end{proof}

\begin{lemma} \label{peu-ram}  Let $\kk$ be a finite
field of characteristic $p$, and suppose that $\rhobar : G_{\Qp}
\rightarrow \GL_2(\kk)$ is tr\`es ramifi\'e. If $\rhobar
\,|_{G_F}$ extends to a finite flat $\kk$-vector space scheme over
the ring of integers $\OO_F$, then $ p \,| e(F)$.
\end{lemma}

\begin{proof}  This lemma follows from the proof of
\cite[Lem. 8.2]{edix} and the discussion that follows it.  Namely, let $r =
[\kk:\Fp]$, so that $\rhobar$ corresponds to an element
$\sigma=(x_1,\ldots,x_r)$ in $(\Qp^{\times}/(\Qp^{\times})^p)^r$;
the assumption that $\rhobar$ is tr\`es ramifi\'e implies that
$\sigma$ does not lie in $(\Zp^{\times}/(\Zp^{\times})^p)^r$, that is,
that some $\val_p(x_i)\not\equiv 0 \pmod{p}$.  If the image
$\sigma_F$ of $\sigma$ in $(F^{\times}/(F^{\times})^p)^r$ then lies
in $(\OO_F^{\times}/(\OO_F^{\times})^p)^r$, it is evident that $p
\,| e(F)$.
\end{proof}

\begin{remark} The behavior in the cases $i=1$ and $i=p$,
$\val_p(b)=0$, is the same as that observed in
\cite[Prop. 8.4]{SavittCompositio}.  We also note that this provides examples
of a Galois representation containing both a lattice whose
reduction is split and a lattice whose reduction is reducible and
nonsplit.
\end{remark}

\begin{cor} \label{found-all} Let $\rho : G_{\Qp} \rightarrow \GL_2(E)$ be a
potentially crystalline representation with Hodge-Tate weights
$\{0,1\}$, and $T$ a Galois-stable lattice inside $\rho$ such that
the reduction $T/\mm_E$ has trivial endomorphisms.
\begin{enumerate}
\item If $\tau(\rho) = \om^i \oplus \om^j$ with $i \not\equiv j
\pmod{p-1}$, then $(T/\mm_E) \,|_{I_p} \otimes_{\kk_E} \Fpbar$ has
one of the three forms
\begin{itemize}
\item
$\begin{pmatrix}
\omega^{1+i} & * \\
0 & \omega^j
\end{pmatrix}$,
\item $\begin{pmatrix}
\omega^{1+j} & * \\
0 & \omega^i
\end{pmatrix}$,
\item $\omega_2^{1 + \{j-i\} + (p+1)i} \oplus \omega_2^{p -
\{j-i\} + (p+1)j}$ where $\{a\}$ denotes the unique integer in
$\{0,\ldots,p-2\}$ which is congruent to $a \pmod{p-1}$.
\end{itemize}

\item If $\tau(\rho) = \omt^m \oplus \omt^{pm}$ with $p + 1 \nmid \ m$,
then $(T/\mm_E) \,|_{I_p} \otimes_{\kk_E} \Fpbar$ has one of the
four forms
\begin{itemize}
\item
$\begin{pmatrix}
\omega^{i+j} & * \\
0 & \omega^{1+j}
\end{pmatrix}$ with $*$ peu ramifi\'e when $i=2$,
\item $\begin{pmatrix}
\omega^{1+j} & * \\
0 & \omega^{i+j}
\end{pmatrix}$ with $*$ peu ramifi\'e when $i=p-1$,
\item $\omega_2^{p+m} \oplus \omega_2^{1+pm}$,
\item $\omega_2^{1+m} \oplus \omega_2^{p(1+m)}$.
\end{itemize}
\end{enumerate}
\end{cor}

\begin{proof}  Part (1) follows, twisting by $\om^i$, from the corresponding
result for type $1 \oplus \om^{j-i}$.  We know that
$D_{\mathrm{st},2}^{F_1}(\rho)$ is described by Proposition
\ref{prin}. If $\val_p(x_1)=0$ or $\val_p(x_2)=0$, then $\rho$ is
actually reducible and the only possible possible reduction of
$\rho$ with trivial endomorphisms is given by part (1) or (2) of
Theorem \ref{red-prin}.  If $0 < \val_p(x_1),\val_p(x_2) < 1$, then
the reduction is given by part (3) of Theorem \ref{red-prin}, and it is
irreducible (hence unique).

For part (2), recall the isomorphism $D_{m,[a:b]} \cong
D_{pm,[bw:-a]}$.  Since $p+1 \nmid \ m$, we know that
$D_{\mathrm{st},2}^{F_2}(\rho)$ is $D_{m,[a:b]}$ for some $[a:b]$.
Suppose first that $\val_p(a)=\val_p(b)$.  If $i \neq 1,p$, then
applying part (1) of Theorem \ref{red-super} to $D_{m,[a:b]}$ yields
a lattice with a reduction of the first kind in the list, and
applying the same result to $D_{pm,[bw:-a]}$ yields a lattice with a
reduction of the second kind.  These are distinct, and so are the
two nontrivial reductions of $\rho$ with trivial endomorphisms (see, 
e.g., Lem. 9.1.1
of Breuil's Barcelona notes \cite{BreuilBarcelona} for the proof
that there are at most two).  If $i = 1,p$, then part (3) of Theorem
\ref{red-super} gives two distinct reductions (since we have assumed
that $T/\mm_E$ has trivial endomorphisms). 

Suppose next that $\val_p(b) >
\val_p(a)$.  If $i > 1$ then part (4) of Theorem \ref{red-super} yields a
reduction of the third kind on the above list, necessarily unique
since it is irreducible; if $i=1$, then part (5) of Theorem
\ref{red-super} yields two reductions (one for each choice of $x_0$ in
(2) of Lemma \ref{extend}) that are both of the first and
second kind.

Finally, if $\val_p(b) < \val_p(a)$, then the previous paragraph
applied to $D_{pm,[bw:-a]}$ gives a unique reduction of the fourth kind when
$i < p$, and two reductions of the first/second kind when $i=p$.
\end{proof}

\begin{remark} \label{really-all} Observe that the reductions in Corollary \ref{found-all}
are the same as those in
\cite[Conjs. 1.2.2, 1.2.3]{CDT}. Note also that it follows from the proof of Corollary
\ref{found-all} that, up to isomorphism, we have actually listed
in Propositions \ref{str-prin} and \ref{str-super} \textit{all}
lattices (in such $\rho$) whose reductions have trivial
endomorphisms.
\end{remark}

\begin{remark} \label{added-note}
We elaborate on the need for the trivial endomorphisms hypothesis in
Corollary \ref{found-all}.  Let $\rho: G_{\Qp} \rightarrow \GL_2(E)$ be a
potentially crystalline representation with Hodge-Tate weights $\{0,1\}$
and $\tau(\rho)$ tame and nonscalar.  If $\rho$ is decomposable, then its reduction is
easy to compute, so we assume that $\rho$ is indecomposable.  If $\rho$ is
not one of the reducible representations considered in Lemma
\ref{reduciblecases}, then in either Theorem \ref{red-prin} or Theorem
\ref{red-super} we have computed the reduction of at least one lattice
contained in $\rho$.  Therefore, in all cases we know the
semisimplification of the reduction of $\rho$.  However, when the
semisimplification is split, we do not claim
to have found all of the lattices $T$ contained in $\rho$ such that
$T/\mE$ has nontrivial endomorphisms.
\end{remark}

\subsection{Application to modular forms}  \label{modforms}
We now apply the results of Section \ref{subsec63} to give a new
computation of the reduction mod $p$ of the local (at $p$)
respresentation attached to a modular form of weight $2$ for
$\Gamma_1(pN)$ (a result due variously to Deligne, Serre, Fontaine,
Gross \cite{GrossTameness}, Edixhoven \cite{edix},...).

\begin{prop} Let $N$ be a positive integer relatively prime to
$p$, let $\chi_p$ be the Teichm\"uller character modulo $p$, and let $\chi_N$
be a Dirichlet character modulo $N$.  Suppose that $f \in
S_2(\Gamma_1(pN),\chi_p^j \chi_N)$ is a normalized cuspidal
newform with $j \in \{1,\ldots,p-2\}$.  Let $\rho_{f,p}$ be the
restriction to $G_{\Qp}$ of the mod $p$ Galois representation
attached to $f$. Then we have the following.
\begin{itemize}
\item If $f$ has slope $0$, then
$ \rhobar_{f,p} \cong
\begin{pmatrix}
\lambda_{\chi_N(p)/a_p} \omega^{j+1} & * \\
0 & \lambda_{a_p}
\end{pmatrix}$.
\item If $f$ has slope $1$, then
$ \rhobar_{f,p} \cong
\begin{pmatrix}
\lambda_{a_p/p} \omega & * \\
0 & \lambda_{\chi_N(p)(p/a_p)} \omega^j
\end{pmatrix}$.
\item If $f$ has slope in the interval $(0,1)$, then
$\rhobar_{f,p} \,|_{I_p} \cong \omega_2^{1+j} \oplus
\omega_2^{p(1+j)}$.
\end{itemize}
\end{prop}

\begin{proof} Let $\rho_{f,p}$ be the restriction to $G_{\Qp}$ of the $p$-adic
Galois representation attached to $f$, so that $\rhobar_{f,p}$ is a
reduction of $\rho_{f,p}$ mod $p$.  We briefly summarize the
(well-known) computation of $D_{\mathrm{st},2}^{F_1}(\rho_{f,p})$
(see, e.g., \cite[Sec. 3.4]{BreuilBarcelona} for more details (of
a dual version)). Faltings
\cite{FaltingsHodgeTate,FaltingsCrystalline} shows that $\rho_{f,p}$
is potentially crystalline, becoming crystalline over $F_1$ with
Hodge-Tate weights $(0,1)$.  By theorems of Saito
\cite{SaitoModularForms} and Deligne, Langlands, and Carayol
\cite{Carayol}, we find that $\tau(\rho_{f,p}) = 1 \oplus \om^j$ and,
if $\rho_{f,p}$ is indecomposable,
$D_{\mathrm{st},2}^{F_1}(\rho_{f,p}) = D_{pa_p^{-1},a_p
\chi_N(p)^{-1}}$. The result now follows from Theorem
\ref{red-prin}. (Note that we do not really need the strong input of
Theorem \ref{red-prin} in the case where $f$ has integer slope,
because $\rho_{f,p}$ is reducible, but we do require it when the
slope is is not an integer.)
\end{proof}

\begin{remark} Techniques of Coleman and Iovita \cite{CI2} may be
used to prove that the representation $\rho_{f,p}$ attached to a
weight $2$ newform for $\Gamma_0(p^2N)$ with $(p,N)=1$ becomes
crystalline over $F_2$.  Thus Theorem \ref{red-super} reduces the
problem of computing $\rhobar_{f,p}$ to the problem of computing
$D_{\mathrm{st}}(\rho_{f,p})$ for such forms.
\end{remark}

\subsection{Families of Galois lattices}  We now describe
explicitly how to arrange our Galois lattices into families. Recall
the elements $V_Y, U_Y \in S_{F_1,\OO_E[[Y]]}$, $V_{X_1}, U_{X_2}
\in S_{F_2,\OO_E[[X_1,X_2]]/(X_1 X_2 - wp)}$, and $W_B, W'_B \in
S_{F_2,\OO_E[[B]]}$ which we described in section
\ref{special-elements}.

\begin{remark} For brevity, we omit the description of $N$ in the
strongly divisible modules below.  In each case, the desired
description is clear from the corresponding strongly divisible
$\OO_E$-modules we have already constructed (and well-defined
using, e.g., the fact that $x_1^2 g_1, x_2 g_2 \in \M_{x_1,x_2}$
in the case when $0 < \val_p(x_1),\val_p(x_2) < 1$, and that $p$
divides $N(W)$).
\end{remark}

\begin{prop} \label{families} There exist strongly divisible modules with descent
data as follows.
\begin{enumerate}
\item  Denoting $x_1 = \tilde{x}_1 (1+Y)$, $x_2 = pw
\tilde{x}^{-1}_1 (1+Y)^{-1}$ and assuming $\tilde{x}_1^2
\not\equiv w \pmod{\mm_E}$ when $j=1$,
$$ \M_{Y_1} =S_{F_1,\OO_E[[Y]]} \cdot g_1 + S_{F_1,\OO_E[[Y]]} \cdot g_2,$$
$$ \Fil^1 \M_{Y_1} = S_{F_1,\OO_E[[Y]]} \cdot (-u^j g_1 + x_1 g_2 ) + ( \Fil^1 S_{F_1,\OO_E[[Y]]}) \M,$$
$$ \phi(g_1) = x_1 g_1, \ \ \phi(g_2) = x_2 g_2 +
u^{pj-e_1}(u^{e_1} + pV_Y) g_1,$$
$$ \ghat(g_1) = g_1, \ \ \ghat(g_2) = \om^j(g) g_2 .$$

\item  Denoting $x_2 = \tilde{x}_2 (1+Y)$, $x_1 = pw
\tilde{x}^{-1}_2 (1+Y)^{-1}$ and assuming $\tilde{x}_2^2
\not\equiv w \pmod{\mm_E}$ when $j=p-2$,
$$ \M_{Y_2} =S_{F_1,\OO_E[[Y]]} \cdot g_1 + S_{F_1,\OO_E[[Y]]} \cdot g_2,$$
$$ \Fil^1 \M_{Y_2} = S_{F_1,\OO_E[[Y]]} \cdot (x_2 g_1 + w u^{e_1-j} g_2 ) + ( \Fil^1 S_{F_1,\OO_E[[Y]]}) \M,$$
$$ \phi(g_1) = x_1 g_1 - wu^{p(e_1-j) - e_1}(u^{e_1} + pU_Y)g_2, \ \ \phi(g_2) = x_2 g_2, $$
$$ \ghat(g_1) = g_1, \ \ \ghat(g_2) = \om^j(g) g_2 .$$

\item  Denoting
$$ D = \left(1 + U_{X_2} V_{X_1} \left(\frac{u^{e_2 p}}{p} +
2u^{(p-1)e_2} + pu^{(p-2)e_2}\right)\right) ,$$
$$ \M_{X_1,X_2} =S_{F_2,\OO_E[[X_1,X_2]]/(X_1 X_2 - pw)} \cdot g_1 + S_{F_2,\OO_E[[X_1, X_2]]/(X_1 X_2 - pw)} \cdot g_2,$$
\begin{eqnarray*}
\Fil^1 \M_{X_1,X_2} & = & S_{F_2,\OO_E[[X_1,X_2]]/(X_1 X_2 - pw)}
\cdot
(-u^k g_1 + (1 \otimes X_1) g_2) \\
 & & + S_{F_2,\OO_E[[X_1,X_2]]/(X_1
X_2 - pw)} \cdot ( (1 \otimes X_2) g_1 + (1 \otimes
w)u^{e_2-k} g_2) \\
 & & + (\Fil^1 S_{F_2,\OO_E[[X_1, X_2]]/(X_1 X_2 -
pw)}) \M,
\end{eqnarray*}
\begin{eqnarray*}
\phi(g_1) & = & (1 \otimes X_1) D^{-1} \left(1 +
\left(\frac{u^{e_2 p}}{p} +
u^{(p-1)e_2}\right)V_{X_1} (U_{X_2}-1)\right)g_1 \\
 & & - (1 \otimes w)D^{-1} u^{p(e_2-k)-e_2}(u^{e_2} + pU_{X_2}) g_2, \\
\phi(g_2) & = & D^{-1} u^{pk-e_2} (u^{e_2} + pV_{X_1})g_1 \\
 & & + (1 \otimes X_2) D^{-1} \left(1 + \left(\frac{u^{e_2 p}}{p} +
u^{(p-1)e_2}\right)U_{X_2} (V_{X_1}-1)\right)g_2,
\end{eqnarray*}
$$ \ghat(g_1) = g_1, \ \ \ghat(g_2) = (\om^j(g) \otimes 1) g_2 .$$

\item Denoting $b = \tilde{b}(1 + B)$, letting $2 \le i \le p$,
and if $i=p$ assuming that $1 + 4 w^2 \tilde{b} \not\equiv 0
\pmod{\mm_E}$ and is a square in $E$,
$$ \M'_B = (S_{F_2,\OO_E[[B]]}) \cdot g_1 + (S_{F_2,\OO_E[[B]]}) \cdot
g_2 ,$$
\begin{equation*} \begin{split}\Fil^1 \M'_B = S_{F_2,\OO_E[[B]]} \cdot ((1 \otimes b)u^{e_2-k} g_2 + (u^{e_2(p-i)} (1\otimes b^2)W'_B -
1)g_1)\\ + (\Fil^1 S_{F_2,\OO_E[[B]]})\M,\end{split}\end{equation*}
\begin{eqnarray*}
\phi(g_1) & = & p g_2 - b W'_B u^{p(e_2-k)} g_1 , \\
\phi(g_2) & = & \left(w - (1\otimes b^2) W'_B \phi(W'_B)
\frac{u^{pe_2(p+1-i)}}{p}\right) g_1 + b\phi(W'_B)
u^{p^2(e_2-k)}g_2,
\end{eqnarray*}
$$ \ghat(g_1) = (\omt^m(g)\otimes 1)g_1, \ \ \ghat(g_2) = (\omt^{pm}(g)\otimes 1) g_2 .$$

\item If $i > 1$,
$$ \M_B = (S_{F_2,\OO_E[[B]]}) \cdot g_1 + (S_{F_2,\OO_E[[B]]}) \cdot
g_2 ,$$
\begin{equation*}\begin{split}\Fil^1 \M_B = S_{F_2,\OO_E[[B]]} \cdot ((1 \otimes B)u^{e_2-k} g_2 + (u^{e_2(p-i)} (1\otimes B^2)W_B -
1)g_1) \\+ (\Fil^1 S_{F_2,\OO_E[[B]]})\M,\end{split}\end{equation*}
\begin{eqnarray*}
\phi(g_1) & = & p g_2 - B W_B u^{p(e_2-k)} g_1 ,\\
\phi(g_2) & = & \left(w - (1\otimes B^2) W_B \phi(W_B)
\frac{u^{pe_2(p+1-i)}}{p}\right) g_1 + B\phi(W_B)
u^{p^2(e_2-k)}g_2,
\end{eqnarray*}
$$ \ghat(g_1) = (\omt^m(g)\otimes 1)g_1, \ \ \ghat(g_2) =
(\omt^{pm}(g)\otimes 1) g_2 .$$

\item If $i=1$ and assuming that $w$ is a square in $E$,
$$ \M_X = (S_{F_2,\OO_E[[B]]}) \cdot g_1 \oplus (S_{F_2,\OO_E[[B]]})\cdot g_2
,$$
$$ \Fil^1 \M_X = S_{F_2,\OO_E[[B]]} \cdot (u^{p-1} g_1 + (w^{-1} X_B + (1
\otimes B) )g_2) + (\Fil^1 S_{F_2,\OO_E[[B]]}) \M_X,$$
\begin{eqnarray*}
\phi(g_1) & = & \phi(X_B) u^{p^2(p-1)} g_1 + \left(1-X_B \phi(X_B )\frac{u^{pe_2}}{pw}\right) g_2  , \\
\phi(g_2) & = & pwg_1 - X_B u^{p(p-1)} g_2,
\end{eqnarray*}
$$\widehat{g}(g_1) = (\tilde{\omega}_2^m \otimes 1)g_1, \qquad 
\widehat{g}(g_2) = (\tilde{\omega}_2^{pm} \otimes 1)g_2.$$
\end{enumerate}
\end{prop}

\begin{proof} In each case, the proof that these formulas define a strongly divisible module
is identical to the proof that the corresponding strongly
divisible $\OO_E$-modules with descent data of Proposition
\ref{str-prin} or \ref{str-super} are indeed strongly divisible
$\OO_E$-modules.
\end{proof}

We adopt the following notation.  For each strongly divisible
$R$-module $\M$ in Proposition \ref{families}, set $R(\M) = R$ (so
that, for example, $R(\M_{Y_1}) = \OO_E[[Y_1]]$).  Set $\tau(\M) = 1
\oplus \om^j$ in the first three cases, and $\tau(\M) = \omt^m
\oplus \omt^{pm}$ in the final two cases.  Finally, set $\rhobar(\M)
= T_{\mathrm{st},2}^{\Qp}(\M/\mm_{R(\M)})$.

\subsection{Deformation rings}

We now come to our main results.

\begin{thm} \label{main-prin} Conjecture 1.2.2 of \cite{CDT} holds;
that is, suppose that $\rhobar : G_{\Qp} \rightarrow \GL_2(\kk_E)$
has trivial endomorphisms.  Suppose that $\tau \cong \om^i \oplus
\om^j$ with $i \not\equiv j \pmod{p-1}$.  Then we have the
following:
\begin{enumerate}
\item $R(2,\tau,\rhobar)_{\OO_E} = 0$ if $\rhobar \,|_{I_p}
\otimes_{\kk_E} \Fpbar \not\in \left\{
\begin{pmatrix}
\omega^{1+i} & * \\
0 & \omega^j
\end{pmatrix},
\begin{pmatrix}
\omega^{1+j} & * \\
0 & \omega^i
\end{pmatrix},
\omega_2^k \oplus \omega_2^{pk}\right\}$ with $k =
1+\{j-i\}+(p+1)i$;

\item $R(2,\tau,\rhobar)_{\OO_E} = \OO_E[[Y]]$ if $\rhobar \,|_{I_p}
\otimes_{\kk_E} \Fpbar \in \left\{
\begin{pmatrix}
\omega^{1+i} & * \\
0 & \omega^j
\end{pmatrix},
\begin{pmatrix}
\omega^{1+j} & * \\
0 & \omega^i
\end{pmatrix}\right\}$;

\item
$R(2,\tau,\rhobar)_{\OO_E} = \OO_E[[X_1,X_2]]/(X_1 X_2 - pw)$ if
$\rhobar \,|_{I_p} \otimes_{\kk_E} \Fpbar \cong \omega_2^k \oplus
\omega_2^{pk}$ with $k = 1+\{j-i\}+(p+1)i$, assuming that $E$
contains $\Qpp$ and that $\kk_E$ contains a square root of
$\det(\rhobar(\Frob_p))$.
\end{enumerate}
\end{thm}

\begin{thm} \label{main-super} Conjecture 1.2.3 of \cite{CDT} holds; 
that is, suppose that $\rhobar : G_{\Qp} \rightarrow \GL_2(\kk_E)$
has trivial endomorphisms.  Suppose that $\tau \cong \omt^m \oplus
\omt^{pm}$ with $p+1 \nmid \, m$.
\begin{enumerate}
\item $R(2,\tau,\rhobar)_{\OO_E} = \OO_E[[B]]$ if $\rhobar
\,|_{I_p} \otimes_{\kk_E} \Fpbar \in \left\{
\begin{pmatrix}
\omega^{i+j} & * \\
0 & \omega^{1+j}
\end{pmatrix},
\begin{pmatrix}
\omega^{1+j} & * \\
0 & \omega^{i+j}
\end{pmatrix}\right\}$, the first
$*$ peu ramifi\'e when $i=2$ and the second when $i=p-1$;

\item $R(2,\tau,\rhobar)_{\OO_E} = \OO_E[[B]]$ if $\rhobar
|_{I_p} \otimes_{\kk_E} \Fpbar \in \left\{ \omega_2^{p+m} \oplus
\omega_2^{1+pm}, \omega_2^{1+m} \oplus \omega_2^{p(1+m)}
\right\}$;

\item $R(2,\tau,\rhobar)_{\OO_E} = 0$ otherwise.
\end{enumerate}
\end{thm}

\begin{thm} \label{th623} Conjecture 2.2.2.4 of \cite{BreuilMezard} (and so, in
particular, \cite[Conj. 1.2.1]{CDT}) holds for $k=2$ and
$\tau$ tame.
\end{thm}

\begin{proof}  We remark that it suffices to prove Theorem
\ref{main-super} and part (2) of Theorem \ref{main-prin} after
extending $E$ in a manner dependent only on $\rhobar$: indeed,
once this result (and the corresponding case of
\cite[Conj. 2.2.2.4]{BreuilMezard}) has been established,
\cite[Lems. 5.1.8, 2.2.2.5]{BreuilMezard} yield the result for our original
$E$.

Part (1) of Theorem \ref{main-prin} and part (3) of Theorem
\ref{main-super} follow immediately from Corollary \ref{found-all}.
In the cases concerning type $\om^i \oplus \om^j$, we may suppose
without loss of generality that $i=0$.   We claim that for each
strongly divisible module $\M$ of Proposition \ref{families}, the
$R(\M)$-representation $T_{\mathrm{st},2}^{\Qp}(\M)$ is actually the
universal deformation of $\rhobar$ to
$R(2,\tau(\M),\rhobar(\M))_{\OO_E}$.  As in the proof of \cite[Th. 5.3.1]{BreuilMezard}, after all of the work that we have
done (the fact that we have found every lattice in a deformation of
$\rhobar$ of type $(2,\tau(\M))$; cf. Prop.
\ref{weak-acc} and Rems. \ref{finext}, \ref{really-all}), it is
essentially formal that there is a canonical injection
$$R(2,\tau(\M),\rhobar(\M))_{\OO_E} \rightarrow R(\M).$$
Abbreviate $R=R(\M)$.  It remains to show that this map is a
surjection; once this is done, the rest of Theorems \ref{main-prin}, \ref{main-super},
and \ref{th623} follows
as in \cite[Sec. 5.3]{BreuilMezard}.

For this surjectivity, it suffices to see that
$T_{\mathrm{st},2}^{\Qp}(\M/(\mm_R^2,\mm_E))$ cannot be defined over
a $\kk_E$-subalgebra of $R/(\mm_R^2,\mm_E)$.  The method used in
\cite{BreuilMezard} is unavailable, as $T_{\mathrm{st},2}$ is not
fully faithful, so we must resort to another (somewhat more
unpleasant) method.  We outline the proof, after which we give the
proof in detail in the most daunting case (part (3) of Th.
\ref{main-prin}).

In most of our cases, $R/(\mm_R^2,\mm_E) = \kk_E[X]/(X^2)$ for a
variable $X$.  Consider the Breuil module $\M'_X =
T_0(\M/(\mm_R^2,\mm_E))$.  If the representation $\rhobar_X =
T_{\mathrm{st},2}^{\Qp}(\M'_X)$ is defined over a
$\kk_E$-subalgebra, that subalgebra can only be $\kk_E$, and in
particular $\rhobar_X$ (regarded simply as a representation over
$\kk_E$) has a subrepresentation $\rhobar'$ such that the
composition $\rhobar' \rightarrow \rhobar_X \rightarrow \rhobar(\M)$
is an isomorphism, where the rightmost map is reduction modulo $X$.
By a scheme-theoretic closure argument, $\M'_X$ has a sub-Breuil
module $\M'$ (with action of $\kk_E$) so that $\M' \rightarrow \M'_X
\rightarrow \M'_X/X\M'_X$ corresponds to a map on group schemes
which is an isomorphism on generic fibres.  (Recall that since $\M$
is a strongly divisible module, reduction modulo $X$ actually
corresponds to the map $\M'_X \rightarrow \M'_X/X\M'_X$ on Breuil
modules.)  In practice, it is too complicated to show directly that
such $\M'$ does not exist.  Fortunately, we know that in every case
(possibly restricting $\rhobar$ to $G_{\Qpp}$ or extending $E$ if
necessary), $\rhobar(\M)$ has a subcharacter $\chi$.  From Remark
\ref{minimal-breuil} and the results of section \ref{sec:chars}, we
can compute the minimal Breuil module $\M''$ corresponding to
$\chi$.  What one proves is that the image of every map $\M''
\rightarrow \M'_X$ falls inside $X\M'_X$, and so the map $\M''
\rightarrow \M'_X/X\M'_X$ is zero and the sought-for $\M'$ cannot
exist.

We demonstrate how this argument can be applied to part (3) of
Theorem \ref{main-prin}.  In this case $R/(\mm_R^2,\mm_E) =
\overline{R} = \kk_E[X_1,X_2]/(X_1^2,X_1 X_2, X_2^2)$, so let $L$
be a linear form in $X_1$ and $X_2$ and suppose that $\rhobar_X$
is defined over the subalgebra $\kk_E[L]$.  Let the corresponding
subrepresentation of $\rhobar_X$ be $\rhobar_L$.  Then the
representation $\rhobar_X/(L)$ is actually defined over $\kk_E$,
and we may apply the argument of the previous paragraph.

We now do this explicitly.  Suppose $E$ is sufficiently large that
$-\ow$ is a square in $\kk_E$.  Since $X_1^2 = X_2^2 = 0$, we see
that $V_{X_1} = U_{X_2} = 1$ in $\Fpp \otimes
\overline{R}[u]/u^{e_2 p}$.  We compute $\M'_{X_1,X_2} =
T_0(\M/(\mm_R^2,\mm_E))$ explicitly from Proposition
\ref{families} and the calculations in the proof of Proposition
\ref{str-prin} and obtain, after a simplifying change of basis,
that $\M'_{X_1,X_2}$ may be generated by $g_1, g_2$ in such a way
that $\Fil^1 \M'_{X_1,X_2}$ is generated by $h_1 = -u^k g_1 + (X_1
- X_2 u^{e_2(p-1-j)})g_2$ and $h_2 = (1\otimes \ow) u^{e_2-k} g_2
+ (X_2 - X_1 u^{e_2 j}) g_1$ satisfying
$$\p(-u^k g_1 + (X_1 - X_2 u^{e_2(p-1-j)})g_2) = (1 \otimes \ow) g_2,$$
$$\p((1\otimes \ow) u^{e_2-k} g_2 + (X_2 - X_1 u^{e_2 j}) g_1) =
(1 \otimes \ow) g_1.$$

The minimal Breuil module $\M''$ of the desired subrepresentation
$\chi$ of $\rhobar(\M)$ restricted to $G_{\Qpp}$ is such that
$\Fil^1\M'' = \M''$, and for some generator $\e$, we have $\p(\e)
= (1 \otimes c)\e$ with $c^2 = -\ow$.  Suppose that we have a
nonzero map $f : \M'' \rightarrow \M'_{X_1,X_2}/(L)$, let
$\overline{X}_1,\overline{X}_2$ denote the images of $X_1$ and
$X_2$ in $\overline{R}/(L)$, and fix $L'$ a non-zero nilpotent in
$\kk_E[X_1,X_2](L,X_1^2,X_1 X_2, X_2^2)$. Our map $f$ must send
$$ \e \mapsto \alpha h_1 + \beta h_2 .$$
Write $\alpha = \alpha_0 u^r + \alpha_L u^t L' $ and $\beta =
\beta_0 u^s + \beta_L u^v L'$ with
$\alpha_0,\alpha_L,\beta_0,\beta_L$ polynomials in $u^e$ over $\Fpp
\otimes \kk_E$ which either are zero or have nonzero constant term.
We wish to show that $\alpha_0 = \beta_0 = 0$.  We consider the
relation $f \p(\e) = \p f(\e)$, first paying attention only to the
terms not involving nilpotents:

$$\phi(\beta_0) u^{ps} c = -\alpha_0 u^{r+k},$$
$$ \phi(\alpha_0) u^{pr} = \beta_0 u^{e_2 - k + s} c.$$

If $\alpha_0$, $\beta_0$ are nonzero, we must therefore have
$r=p-j$ and $s=1+j$.  We turn next to the terms involving
nilpotents. The $g_1$-term in this relation is:
$$ w\phi(\beta_L) u^{pv} L' = - c\alpha_L u^{t+k} L' + c\beta_0
u^{1+j}(\overline{X}_2 - \overline{X}_1 u^{e_2 j}).$$ But if
$\beta_0 \overline{X}_2 \neq 0$, equality could not hold here,
because there can be no other terms of degree $1+j$ in $u$!  If
$\beta_0 \neq 0$ it follows that $\overline{X}_2 = 0$.  But
similar consideration of the $g_2$-term yields $\overline{X}_1=0$.
Since $\overline{X}_1$ and $\overline{X}_2$ cannot both be zero,
it follows that $\alpha_0 = \beta_0 = 0$, and we are done.

We note very briefly some of the features of this calculation for
the other parts of Theorems \ref{main-prin} and \ref{main-super}.  In part (2) of Theorem
\ref{main-prin}, the case $\tilde{x}_1^2 \equiv w \pmod{\mm_E}$
requires slightly more work (in most cases an $\alpha_0$ is forced
to be zero on its own, but in the more complicated case, one needs
to use $j \neq 1$ and consider $\beta_0$ as well to see that
$\alpha_0 = 0$). There is a similar feature in part (1) of Theorem
\ref{main-super} when $\tilde{b}^2 w \equiv \pm 1 \pmod{\mm_E}$; in
this case, there is a $\beta_0$ which satisfies $\phi(\beta_0) = \mp
\beta_0$, and then an equation of the form $\pm \beta_B =
\phi(\beta_0) - \phi(\beta_B)$ implies $\beta_0=0$. (Apply $\phi$ to
this equation again.)
\end{proof}

\begin{cor} The Breuil-M\'ezard conjecture
\cite[Conj. 2.3.1.1]{BreuilMezard} holds for $k=2$ and $\tau$ tame.
\end{cor}

\begin{proof} This is immediate from the computation of
$\mu_{\mathrm{aut}}(2,\rhobar,\tau)$ with $\tau$ tame.
\end{proof}

\begin{cor} Theorem \ref{modularity} holds.
\end{cor}

\subsection*{Acknowledgements}  The author is especially indebted to Christophe
Breuil for numerous discussions; without Breuil's generosity this
work would never have been completed. The author is grateful to
Richard Taylor for suggesting this problem, and to the anonymous
referee for his or her helpful comments. The author thanks the
Institut des Hautes \'Etudes Scientifiques, Bures-sur-Yvette,
France, for its warm hospitality during a year-long stay in
2002-2003.

\bibliographystyle{amsalpha}
\bibliography{savitt}

\providecommand{\bysame}{\leavevmode\hbox to3em{\hrulefill}\thinspace}
\providecommand{\MR}{\relax\ifhmode\unskip\space\fi MR }
\providecommand{\MRhref}[2]{%
  \href{http://www.ams.org/mathscinet-getitem?mr=#1}{#2}
}
\providecommand{\href}[2]{#2}
\begin{thebibliography}{BCDT01}

\bibitem[BCDT01]{BCDT}
Christophe Breuil, Brian Conrad, Fred Diamond, and Richard Taylor, \emph{On the
  modularity of elliptic curves over {${\mathbf{{{Q}}}}$}}, J.A.M.S.
  \textbf{14} (2001), 843--939.

\bibitem[BM02]{BreuilMezard}
Christophe Breuil and Ariane M\'{e}zard, \emph{Multiplicit\'{e}s modulaires et
  repr\'{e}sentations de ${{\rm {{G}{L}}}}_2(\mathbf{Z}_p)$ et de
  $\rm{{G}al}({\Qbar}_p/{\Q}_p)$ en $l=p$}, Duke Math. J. \textbf{115} (2002),
  no.~2, 205--310, With an appendix by Guy Henniart.

\bibitem[Bre]{BreuilBarcelona}
Christophe Breuil, \emph{p-adic {H}odge theory, deformations, and local
  {L}anglands}, {N}otes from a course given at {C}entre de {R}ecerca
  {M}atematica, {J}uly 18 to 28, 2001, available at {\sf
  http://www.crm.es/Publications/ps-pdf-preprints/Quadern20-1.pdf}.

\bibitem[Bre98]{BreuilENS}
\bysame, \emph{Construction de repr\'esentations {$p$}-adiques semi-stables},
  Ann. Sci. \'Ecole Norm. Sup. (4) \textbf{31} (1998), no.~3, 281--327.

\bibitem[Bre99]{BreuilStrong}
\bysame, \emph{Repr\'{e}sentations semi-stables et modules fortement
  divisibles}, Invent. math. \textbf{136} (1999), no.~1, 89--122.

\bibitem[Bre00]{Br}
\bysame, \emph{Groupes $p$-divisibles, groupes finis et modules filtr\'{e}s},
  Ann. of Math. (2) \textbf{152} (2000), no.~2, 489--549.

\bibitem[Car86]{Carayol}
Henri Carayol, \emph{Sur les repr\'esentations {$l$}-adiques associ\'ees aux
  formes modulaires de {H}ilbert}, Ann. Sci. \'Ecole Norm. Sup. (4) \textbf{19}
  (1986), no.~3, 409--468.

\bibitem[CDT99]{CDT}
Brian Conrad, Fred Diamond, and Richard Taylor, \emph{Modularity of certain
  potentially {B}arsotti-{T}ate {G}alois representations}, J.A.M.S. \textbf{12}
  (1999), 521--567.

\bibitem[CF00]{ColmezFontaine}
Pierre Colmez and Jean-Marc Fontaine, \emph{Construction des
  repr\'{e}sentations p-adiques semi-stables}, Invent. math. \textbf{140}
  (2000), no.~1, 1--43.

\bibitem[CI]{CI2}
Robert Coleman and Adrian Iovita, \emph{Hidden structures on semi-stable
  curves}, Preprint.

\bibitem[Edi92]{edix}
Bas Edixhoven, \emph{The weight in {S}erre's conjectures on modular forms},
  Invent. math. \textbf{109} (1992), no.~3, 563--594.

\bibitem[Fal87]{FaltingsHodgeTate}
Gerd Faltings, \emph{Hodge-{T}ate structures and modular forms}, Math. Ann.
  \textbf{278} (1987), no.~1-4, 133--149.

\bibitem[Fal97]{FaltingsCrystalline}
\bysame, \emph{Crystalline cohomology of semistable curve---the {${\bf Q}\sb
  p$}-theory}, J. Algebraic Geom. \textbf{6} (1997), no.~1, 1--18.

\bibitem[FM95]{FontaineMazur}
Jean-Marc Fontaine and Barry Mazur, \emph{Geometric {G}alois
  {R}epresentations}, Elliptic curves, modular forms, \& {F}ermat's last
  theorem (John Coates and Shing-Tung Yau, eds.), International Press, 1995,
  pp.~41--78.

\bibitem[Fon94]{AST223}
Jean-Marc Fontaine, \emph{Representations $p$-adiques semi-stables},
  Ast\'{e}risque \textbf{223} (1994), 113--184.

\bibitem[Gro90]{GrossTameness}
Benedict~H. Gross, \emph{A tameness criterion for {G}alois representations
  associated to modular forms (mod {$p$})}, Duke Math. J. \textbf{61} (1990),
  no.~2, 445--517.

\bibitem[Kis]{KisinModularity}
Mark Kisin, \emph{Moduli of finite flat group schemes and modularity},
  Preprint.

\bibitem[Ray74]{Raynaud}
Michel Raynaud, \emph{Sch\'{e}mas en groupes de type $(p,p,\ldots,p)$}, Bull.
  Soc. Math. France \textbf{102} (1974), 241--280.

\bibitem[Sai97]{SaitoModularForms}
Takeshi Saito, \emph{Modular forms and {$p$}-adic {H}odge theory}, Invent.
  Math. \textbf{129} (1997), no.~3, 607--620.

\bibitem[Sav04]{SavittCompositio}
David Savitt, \emph{Modularity of some potentially {B}arsotti-{T}ate {G}alois
  representations}, Compos. Math. \textbf{140} (2004), no.~1, 31--63.

\bibitem[Tat67]{tate}
John Tate, \emph{$p$-divisible groups}, Proc. Conf. Local Fields (Driebergen,
  1966), Springer-Verlag, 1967, pp.~158--183.

\bibitem[Tat97]{Tate2}
\bysame, \emph{Finite flat group schemes}, Modular forms and Fermat's last
  theorem (Boston, MA, 1995), Springer, New York, 1997, pp.~121--154.

\end{thebibliography}

\vskip0.3cm

\noindent {\sc Department of Mathematics, McGill University, 805
Sherbrooke St West, Montr\'eal, Canada, H3A 2K6;}

\vskip0.2cm

\noindent {\sf dsavitt@math.mcgill.ca}

\end{document}